\documentclass[final]{amsart}
\usepackage{amssymb}
\usepackage{amsaddr}
\usepackage[utf8]{inputenc}
\usepackage[T1]{fontenc}
\usepackage[english]{babel}
\usepackage{csquotes}
\usepackage{comment}
\usepackage{xargs}
\usepackage{graphicx, color}
\usepackage[figurewithin=none]{caption,subcaption} 
\graphicspath{{./}} \usepackage{enumitem}
 \usepackage[notcite,notref]{showkeys}
 \usepackage{hyperref}
 \hypersetup{
 	hidelinks,
   colorlinks  = true,    urlcolor    = blue,    linkcolor   = blue,    citecolor   = blue,      pdfauthor   = {Isabel Hubard, Elías Mochán, Antonio Montero },pdfsubject  = {Maniplexes},pdftitle    = {Voltage operations on maniplexes},  }

\usepackage[capitalise,noabbrev, nameinlink, poorman]{cleveref}
\usepackage{tikz-cd}
\tikzcdset{row sep=tiny,arrows={shorten=-4pt}
}

\usepackage[ shadow, colorinlistoftodos,textsize=tiny, obeyFinal]{todonotes}
\setuptodonotes{fancyline, backgroundcolor=gray!10,bordercolor=gray}

\newcommandx{\inline}[2][1=]{\todo[inline, #1]{#2}}

\makeatletter
    \providecommand\@dotsep{5}
\makeatother

\theoremstyle{plain}
\newtheorem{thm}{Theorem}[section]

\newtheorem{lemma}[thm]{Lemma}
\newtheorem{prop}[thm]{Proposition}

\newtheorem{coro}[thm]{Corollary}

\theoremstyle{definition}

\newtheorem{exam}[thm]{Example}
\newtheorem{eg}[thm]{Example}

\theoremstyle{remark}

\newtheorem{rem}[thm]{Remark}

\newtheorem{question}[thm]{Question}
\newtheorem{problem}{Problem}

\numberwithin{equation}{section}

\usepackage{xargs}

\renewcommand{\leq}{\leqslant} \renewcommand{\geq}{\geqslant}
\renewcommand{\epsilon}{\varepsilon} \renewcommand{\subset}{\subseteq}  
\renewcommand{\{}{\lbrace}
\renewcommand{\}}{\rbrace}

\renewcommand{\bar}{\overline}

\newcommand{\bZ}{\mathbb{Z}}
\newcommand{\bD}{\mathbb{D}}

\newcommand{\cC}{\mathcal{C}}
\newcommand{\cF}{\mathcal{F}}
\newcommand{\cG}{\mathcal{G}}

\newcommand{\cM}{\mathcal{M}}
\newcommand{\cN}{\mathcal{N}}
\newcommand{\cP}{\mathcal{P}}

\newcommand{\cS}{\mathcal{S}}
\newcommand{\cT}{\mathcal{T}}
\newcommand{\cU}{\mathcal{U}}
\newcommand{\cX}{\mathcal{X}}
\newcommand{\cY}{\mathcal{Y}}

\newcommand{\pman}{\mathrm{pMpx}}

\DeclareMathOperator{\aut}{\Gamma}

\DeclareMathOperator{\rk}{rk}

\DeclareMathOperator{\mon}{Mon}

\DeclareMathOperator{\fg}{\Pi}

\DeclareMathOperator{\oo}{\mathcal{O}}
\DeclareMathOperator{\pyr}{Pyr}
\DeclareMathOperator{\prism}{Pri}
\DeclareMathOperator{\med}{Med}

\newcommand{\mix}{\diamondsuit}

\newcommandx{\Proot}[2][1=\cP, 2= \baseFlag, usedefault]{\left( #1,#2 \right)}

\newcommandx{\pth}[2][1=i, usedefault]{ {}^{#1} #2}
\DeclareMathOperator{\vr}{{V}}

\newcommandx{\ertimes}[1][1=\eta, usedefault]{\rtimes_{#1}}
\newcommandx{\eltimes}[1][1=\eta, usedefault]{\rlimes_{#1}}
\newcommand{\mixer}{\mu}
\newcommand{\mixtimes}[1]{\ertimes[\mixer]{#1}}
\DeclareMathOperator{\ot}{{ot}}
\newcommand{\1}{{\mathbf{1}}}
\newcommand{\2}{{\mathbf{2}}}
\DeclareMathOperator{\trp}{{Trp}}

\newcommandx{\mlink}[6][1,3,5,6]{\arrow[dash, #1]{#2}[font=\tiny, description, #3]{#4}[outer sep=3pt, #5]{#6}}
\newcommandx{\mdart}[6][1,3,5,6]{\arrow[#1]{#2}[font=\tiny, description, #3]{#4}[outer sep=3pt, #5]{#6}}
\newcommandx{\msemi}[6][1,3,5,6]{\arrow[dash, #1]{#2}[font=\tiny, description, near start, #3]{#4}[near end, #5]{#6}}

\newcommand{\gen}[1]{\left\langle #1 \right\rangle}

\usepackage[style=numeric, sorting=nyt,  maxbibnames=9, backend=biber]{biblatex}
\begin{document}

\title{Voltage operations on maniplexes}
\author{Isabel Hubard}
\address{Institute of Mathematics, National Autonomous University of Mexico (IM UNAM), 04510 Mexico City, Mexico}
\email{isahubard@im.unam.mx}
\author{Elías Mochán} 
\address{College of Science, Northeastern University, 02115 Boston, USA}
\email{j.mochanquesnel@northeastern.edu}
\author{Antonio Montero}
\address{Institute of Mathematics, National Autonomous University of Mexico (IM UNAM), 04510 Mexico City, Mexico}
\email{amontero@im.unam.mx}

\maketitle

\listoftodos\relax

\begin{abstract}
  Classical geometric and topological operations on polyhedra, maps and polytopes often give rise to structures with the same symmetry group as the original one, but with more flags.
    In this paper we introduce the notion of \emph{voltage operations} on maniplexes, as a way to generalize such operations.
Our technique  provides
     a way to study classical operations in a graph theory setting.
     In fact, voltage operations can be applied to symmetry type graphs, and more generally to $n$-valent properly $n$-edge colored graphs.
We focus on studying the interactions between voltage operations and the symmetries of the operated object, and show that they can be potentially used to build maniplexes with prescribed symmetry type graphs. Moreover, a complete characterization of when an operation can be seen as a voltage operation is given.
\end{abstract}
 
\section{Introduction}

Traditionally, operations on polyhedra  introduce local transformations on the vertices, edges or faces. Classical examples include the Wythoffian ones, such as truncation and medial.
Operations in maps and polytopes have been studied previously.
For example, in \cite{OrbanicPellicerWeiss_2010_MapOperations$k$}, Orbanič, Pellicer and Weiss explore map operations to build $k$-orbit maps.
They take a combinatorial approach and divide each flag  (incident triplets vertex-edge-face) of a map into several new ones.
In a very recent paper \cite{CunninghamPellicerWilliams_StratifiedOperationsManiplexes_preprint} Cunningham, Pellicer and Williams introduce \emph{stratified operations} to study monodromy groups of some important families of non-reflexive maniplexes (a generalization of -the flag graph of- polytopes, where some conditions are relaxed).
This concept is closely related to our concept of voltage operation and many of their results can be written in our language and vice versa.

Classical operations often give rise to maps or polytopes with the same symmetry group as the original one, but with more flags.
Thus, by applying one of these operations to a polytope we get a new polytope with more flag orbits. Furthermore, it is not difficult to see that, often, the ``local configuration'' of the flags does not depend on the original polytope but on the operation per se.
Could we use operations to get \emph{any} local configuration of the flag orbits?
We shall explain this question more precisely.

The flag graph of a polytope $\cP$ is a properly edge-colored $n$-valent graph whose vertices are the flags of $\cP$ and two of them are adjacent if the corresponding flags differ in exactly one element.
The colors of the edges are given by the type of element they differ in.
The quotient of the flag graph of $\cP$ by its symmetry group is the \emph{symmetry type graph} of $\cP$.
The symmetry type graph can be thought as a way to represent the local configuration of the flag orbits.

We can notice that when we apply a classical operation to different regular polyhedra, the resulting polyhedra often have isomorphic symmetry type graphs.
(For example, if $\med$ represents the medial operation and $\cP$ is a regular polyhedron, the symmetry type graph of $\med(\cP)$ will have either one single vertex, or two vertices and an edge of color $2$ joining them, depending on whether or not $\cP$ is self-dual.)
Moreover, if we apply an operation $\oo$ to two polyhedra, both with the same symmetry type graph $\cT$, it is very likely that both resulting polyhedra have the same  symmetry type graph $\cT'$.

A natural question to ask is: given a properly edge-colored $n$-valent graph $\cT$, is there a maniplex whose symmetry type graph is $\cT$?
Going back to operations, can we find an operation $\oo$ such that if $\cP$ is a reflexible maniplex, then the symmetry type graph of $\oo(\cP)$ is precisely $\cT$? 

These questions give rise to \emph{voltage operations} as a potential useful tool to find answers.
Voltage operations define a technique that uses voltage graphs to generalize the above mentioned classical operations, as well as many others, even in higher dimensions (or ranks).
These operations can also be defined for maniplexes.
But not only that, they can also be defined for \emph{premaniplexes}, a more general concept that includes both maniplexes and their symmetry type graphs.

It is important to remark that voltage graphs have been used before in constructing polytope-like structures.
Notably, in \cite{PellicerPotocnikToledo_2019_ExistenceResultTwo} the authors use voltage graphs to build $2$-orbit $n$-maniplexes (for $n \geq 4$) with prescribed symmetry type graphs.
In \cite{CunninghamDelRioFrancosHubardToledo_2015_SymmetryTypeGraphs} the authors use voltage graphs, without explicitly mentioning them, to find generators and relations for the automorphism group of $3$-orbit polytopes.
In \cite{Mochan_2021_AbstractPolytopesTheir_PhDThesis} the second author of this manuscript exploited extensively the use of voltage graphs to solve some relevant problems in the area, namely, problems 1 and 2 in    \cite{CunninghamPellicer_2018_OpenProblems$k$}.

We start this paper by giving some basic definitions that we shall need throughout the manuscript. In \cref{sec:operations} we give the definition of a voltage operation and look in detail into the mix as a voltage operation; we also study some connectivity properties of the voltage operations. In Section\nobreakspace \ref {sec:ClassicalExamples}  we explore some examples of classical operations on polytopes and polyhedra.
In particular we show how pyramids, prisms, Wythoffian operations, among others, can be seen as voltage operations.
In Section\nobreakspace \ref {sec:automorphisms} we describe how automorphisms and voltage operations interact.
We show that if $\cM$ is a maniplex with symmetry type $\cT$ and $\oo$ is a voltage operation, then every automorphism of $\cM$ induces an automorphism of $\oo(\cM)$.
That is, the automorphism group of $\cM$ is a subgroup of the automorphism group of $\oo(\cM)$.
Moreover, we show that the symmetry type graph of $\oo(\cM)$ with respect to the automorphism group of $\cM$ can be obtained y applying the voltage operation to $\cT$. In Section\nobreakspace \ref {sec:composition} we show that voltage operations are closed under composition and we find a simple way to describe the composition of two voltage operations.
Finally, in Section\nobreakspace \ref {sec:universal} we give conditions for two voltage operations to be equivalent.
We also characterize voltage operations in terms of what they do to the universal maniplex.
More precisely, in Theorem\nobreakspace \ref {thm:covers} we prove that if $\oo$ is an application that assigns a premaniplex $\oo(\cX)$ to each premaniplex $\cX$, then we prove that $\oo$ is a voltage operation if and only if $\oo(\cU/\Gamma) \cong \oo(\cU)/\Gamma$ for every group $\Gamma \leq \aut(\cU)$.
 \section{Preliminaries}\label{sec:defs}

\inline{Pasar un corrector y dejar todo en un sólo tipo de inglés}

\subsection{Graphs}

In this work we use the definition of graph used in \cite{MalnicNedelaSkoviera_2000_LiftingGraphAutomorphisms}, which is slightly more general that the usual definition.
A \emph{graph} $X$ is a quadruple $(D,V,I,(\cdot)^{-1})$ where $D$ and $V$ are disjoint sets, $I:D \to V$ is a mapping and $(\cdot)^{-1}: D \to D$ is an involutory permutation of $D$.
The set $V$ is the set of \emph{vertices} of $X$, the set $D$ is the set of \emph{darts} of $X$.
For a dart $d$, the vertex $I(d)$ is the \emph{initial vertex} or \emph{starting point} of $d$ and $d^{-1}$ is the \emph{inverse} of $d$. The \emph{terminal vertex} or \emph{endpoint} of $d$ is the starting point of $d^{-1}$.
The \emph{edges} are the orbits of $D$ under the action of $(\cdot)^{-1}$. If an edge consists only of a single dart, that is, a dart $d$ that satisfies $d^{-1}=d$, then it is called a \emph{semiedge}.
A \emph{loop} is an edge consisting of two darts whose initial vertex is the same.
A \emph{link} is an edge that is not a loop or a semiedge. That is, an edge whose two darts are different and have different starting points.

Usually a graph $X$ is defined by their vertex set $V(X)$ and the set $E(X)$ of edges.
We can recover the set of darts as the set of formal pairs \[D(X)= \left\{ \pth[e]{v} : v \in V(X), e \in E(X) \text{ and $v$ is incident to $e$}  \right\}. \]
The function $I$ is given by $I(\pth[e]{v})=v$ and if an edge $e$ is incident to the vertices $u$ and $v$, then $\left( \pth[e]{v} \right)^{-1} = \pth[e]{u}$.
Formally speaking, this definition turns loops into semiedges.
In this work we will not use graphs with loops hence we allow this small abuse.

If $X=(V(X),D(X), I_{X}, (\cdot)^{-1}_{X})$ and $Y=(V(Y),D(Y), I_{Y}, (\cdot)^{-1}_{Y})$ are graphs, a \emph{graph homomorphism} $f: X \to Y$ is a pair of functions $f=(f_V, f_D)$ such that $f_V: V(X) \to V(Y)$, $f_D: D(X) \to D(Y)$ that satisfy
\[\begin{aligned}
I_{Y}\left( d f_{D} \right) &= \left( I_{X} (d) \right)f_V && \text{and} \\
\left( (d)f_D \right)_{Y}^{-1} &= ( d^{-1}_{X} )f_D
\end{aligned}
\]
for every dart $d \in D(X)$.
Note that we evaluate graph homomorphisms on the right. If $v\in V(X)$ and $d \in D(X)$, we shall write $vf$ and $df$ instead of $vf_{V}$ and $df_{D}$.
If both $f_{V}$ and $f_{D}$ are bijective and $f^{-1}:=(f_{V}^{-1},f_{D}^{-1})$ is also a graph homomorphism, then we say that $f$ (and $f^{-1})$ is an \emph{isomorphism} and that $X$ and $Y$ are \emph{isomorphic} (and write $X \cong Y$).
Naturally, for a graph $X$, an \emph{automorphism} is an isomorphism $f:X \to X$.

A \emph{path} is a finite sequence  $W=(d_{1}, \dots, d_{k})$ of darts such that the endpoint of $d_{i}$ is the starting point of $d_{i+1}$ for $i \in \{0, \dots, k-1\}$.
We usually omit commas and parentheses and simply write $W=d_{1} \dots d_{k}$.
The number $k$ is the \emph{length} of $W$. The starting point of $d_1$ is the \emph{starting point} of $W$; we also say that $W$ starts at $I(d_{1})$. Similarly, the endpoint of $W$ is the endpoint of $d_{k}$ and we say that $W$ ends at this vertex. A single vertex is a path of length $0$. A path is \emph{closed} if its starting point and its endpoint are the same vertex.

Notions such as degree of a vertex, cycle, subgraphs, connectivity, connected components and trees extend naturally from the classic definition of graphs and paths.

Given a path $W = d_{1} \dots d_{k}$, an \emph{elementary (graph) move} consist in inserting or removing a pair of consecutive inverse darts at any point in the sequence $d_{1}, \dots, d_{k}$.
If a path $V$ can be obtained from a path $W$ by applying a series of elementary moves then we say that the paths are \emph{graph-homotopic} (and write $W \sim V$).
Clearly graph-homotopy is an equivalence relation and we often identify a path with its homotopy class.
Observe that if $W$ is graph-homotopic to a path of length $0$ then $W$ must be closed.

Given two paths $W$ and $V$ we say that they are \emph{compatible} if the endpoint of $W$ is the starting point of $V$.
We can operate compatible paths by concatenation.
This is, if $W=d_{1} \dots d_{k}$ and $V=a_{1} \dots a_{\ell}$ then $WV= d_{1}  \dots  d_{k}  a_{1}  \dots  a_{\ell}$.
If $W \sim W'$ and $V \sim V'$ then $WV \sim W'V'$, which implies that we can operate not only compatible paths but homotopy classes of compatible paths.

The \emph{fundamental groupoid} of a graph $X$, denoted by $\fg(X)$ is the set of graph-homotopy classes of $X$ with the partial operation defined above.
If $u$ is a vertex in $X$, then the fundamental group of $X$ at $u$, denoted by $\fg^{u}(X)$ is the set of graph-homotopy classes of closed paths at $u$ with concatenation as operation.
Observe that $\fg^{u}(X)$ is actually a group and that if $V$ is a path from $v$ to $u$ then  $\fg^{v}(X) = V \fg^{u}(X) V^{-1} $.

If $X$ is a graph and $G$ is a group, a \emph{voltage assignment} is a function $\xi: \fg(X) \to G$ that satisfies $\xi(WV) = \xi(V)\xi(W)$  for any two compatible paths $W$ and $V$
\footnote{Usually a voltage assignment $\xi$ is defined such that $\xi(WV)=\xi(W)\xi(V)$ (cf. \cite{MalnicNedelaSkoviera_2000_LiftingGraphAutomorphisms}). The reason we do it the other way is because we are considering right actions, as is customary in the polytopes and maniplexes literature.}.
The group $G$ is called the \emph{voltage group} of $\xi$ and the pair $(X,\xi)$ is called a \emph{voltage graph}.

Since every path can be though as the product of its darts, we can regard a voltage assignment as a function $\xi:D(X) \to G$ such that for every dart $d$, $\xi(d^{-1})= \xi(d)^{-1}$. The voltage of a path $W= d_{1}, \dots, d_{k}$ is simply $\xi(d_{k}) \cdots \xi(d_{1})$.

If $(X,\xi)$ is a voltage graph, the \emph{derived graph} is the graph $X^{\xi}$ whose vertices and darts are the elements in $V(X) \times G$ and $D(X) \times G$, respectively.
The initial vertex of a dart $(d,g)$ is $(I(d),g)$ and its inverse is $(d^{-1}, \xi(d)g)$.
Informally speaking, the adjacent vertices of a given vertex $(x,g)$ are the vertices $(x, \xi(d)g)$ for each dart $d$ starting at $x$.

If $X$ is a graph and $\xi: \fg(X) \to \Gamma$ and $\zeta: \fg(X) \to \Gamma$ are voltage assignments, then $\xi$ and $\zeta$ are \emph{equivalent} if there is an isomorphism $X^{\xi} \to X^{\zeta}$ such that the following diagram commutes:
\begin{equation}\label{eq:EquivVolts}
\begin{tikzcd}  [column sep = small]
X^\xi \arrow[dashed]{rr} \arrow{rd}{} & & X^{\zeta} \arrow{dl}{} \\
& X &
\end{tikzcd}
\end{equation}
where the arrows pointing to $X$ are the projection to the first coordinate, which is a homomorphism.

The following Theorem is well known (see \cite{MalnicNedelaSkoviera_2000_LiftingGraphAutomorphisms}):

\begin{thm}\label{thm:VoltsEquiv}
    If $(\cX,\xi)$ is a connected voltage graph and $x_0$ is a vertex in $\cX$, there exists a voltage assignment $\xi'$ on $\cX$ satisfying that:
    \begin{itemize}
        \item $\xi'$ is equivalent to $\xi$;
        \item for every dart $d$ in $\cX$, $\xi'(d)\in\xi(\fg^{x_0}(\cX))$; and
        \item there is a spanning tree $T$ of $\cX$ such that for every dart $d$ in $T$, $\xi'(d)=1$.
    \end{itemize}
\end{thm}

Note that if $\cX$ is connected, without loss of generality, we can assume that a voltage assignment $\xi$ on $\cX$ satisfies the last two conditions of Theorem\nobreakspace \ref {thm:VoltsEquiv}.

We can think of the spanning tree $T$ as a \emph{fundamental region} for the group $G$.

\subsection{Premaniplexes and maniplexes}
A \emph{properly $n$-edge-colored graph} is a graph $X$ and a function $c: D(x) \to \left\{ 0, \dots, n-1 \right\} $ such that $c(d) = c(d^{-1})$ for every dart $d$ and if $d_{1}$ and $d_{2}$ are such that $I(d_{1}) = I(d_{2})$, then $c(d_{1})\neq c(d_{2})$.
Observe that $c$ induces a proper edge-coloring in the classical sense.

A \emph{premaniplex of rank $n$} or simply $n$-premaniplex is a properly $n$-edge-colored graph $\cX$ such that every vertex is the starting point of one dart of each color and if $i,j \in \left\{0, \dots, n-1 \right\} $ are such that $\left| i-j \right| \geq 2 $, then the alternating paths of length 4 with colors $i,j$ are closed.

If a premaniplex $\cX$ is connected and simple, that is, if there are no semiedges or parallel edges, then we say that $\cX$ is a maniplex.
Maniplexes were introduced in \cite{Wilson_2012_ManiplexesPart1} as a combinatorial generalization of maps to higher ranks.
Constructions of maniplexes with given symmetry properties can be found in \cite{PellicerPotocnikToledo_2019_ExistenceResultTwo}.
Natural examples of maniplexes are the flag graphs of polytopes.
In fact, if $\cX$ is a maniplex, we usually call its vertices \emph{flags}.
Moreover, if $\cM$ is an $n$-maniplex and $i \in \{0, \dots, n-1\}$, the \emph{$i$-faces} of $\cM$ are the connected components of $\cM$ after removing the edges of color $i$. More generally, if $0 \leq k < \ell \leq n-1$ the \emph{$(k,\ell)$-sections} of $\cM$ are the connected components of $\cM$ after removing the edges of color $i$ if $i < k$ or $i > \ell$.
Note that when $\cM$ is the flag graph of a polytope $\cP$, the faces and sections of $\cM$ are in correspondence with the faces and sections of $\cP$.
These close relation between flag graphs of polytopes and maniplexes allows us to think of polytopes in a graph-theoretical approach.
In this paper, we slightly abuse of language and whenever we talk of a polytope $\cP$ we actually refer to the flag graph $\cP$.
In \cite{GarzaVargasHubard_2018_PolytopalityManiplexes} Garza-Vargas and Hubard characterize when a maniplex is the flag graph of a polytope.
Finally, observe that the notions of flags, sections and faces extend naturally to premaniplexes.

If $\cX$ is a premaniplex, whenever we write $x \in \cX$ we mean that that $x$ is a flag (vertex) in $\cX$.
If $x$ is a vertex in an $n$-premaniplex, and $i \in \left\{ 0, \dots, n-1 \right\} $ we denote by  $\pth{x}$ the dart of color $i$ whose starting point is $x$. We denote by $x^{i}$ the \emph{$i$-adjacent vertex} of $x$, that is, the endpoint of the dart $\pth{x}$.

Given two $n$-premaniplexes $\cX$ and $\cX'$, a \emph{(premaniplex) homomorphism} from $\cX$ to $\cX'$ is a function that preserves $i$-adjacencies, for $i\in\{0, \dots, n-1\}$.
In the particular case where $\cX$ and $\cX'$ are prepolytopes, these homomorphisms are called rap-maps.
It is easy to prove that if $\cX'$ is connected, every homomorphism with codomain $\cX'$ is surjective (see  \cite[Lemma 2.5]{MonsonPellicerWilliams_2014_MixingMonodromyAbstract} for a proof for rap-maps that naturally extends to homomorphism of premaniplexes). In general, every connected component of $\cX'$ is either contained completely in the image of an homomorphism or in its complement.
The class of all $n$-premaniplexes together with the premaniplex homomorphisms as arrows forms a category which we will call $\pman^n$.
A surjective homomorphism is called a \emph{covering}. If there is a covering from $\cX$ to $\cX'$ we say that \emph{$\cX$ covers $\cX'$}.
It is easy to see that if $\cX'$ is connected, any homomorphism from $\cX$ to $\cX'$ is a covering.
The notions of \emph{isomorphism} and \emph{automorphism} of premaniplexes are defined in the usual way.
The automorphism group of $\cX$ is denoted by $\aut(\cX)$, and acts naturally on flags, $i$-faces and $(k,l)$-sections. Such actions will be consider as right actions.
A premaniplex $\cX$ is \emph{regular} if $\aut(\cX)$ acts transitively on vertices.
We say that $\cX$ is a \emph{$k$-orbit} premaniplex if $\aut(\cX)$ induces $k$ orbits on vertices.
These notions coincide with the analogue notions for polytopes.
For example, the flag graph of a polytope is a regular maniplex if and only if the polytope is regular itself.

If $\cX$ is an $n$-premaniplex and $\Gamma \leq \aut(\cX)$ the \emph{quotient} $\cX/\Gamma$ is the $n$-colored graph whose vertices are the orbits $\{x\Gamma  :  x \in V(\cX)\}$ and for $i \in \{0, \dots, n-1\}$ $(x\Gamma)^{i} = (x^{i})\Gamma$.
Observe that the adjacencies are well-defined since $(x^{i}) \gamma = (x \gamma)^{i}$ for every $\gamma \in \aut(\cX)$.
Moreover, it is straightforward that $\cX/\Gamma$ is an $n$-premaniplex.

If $\cM$ is a maniplex and $\Gamma \leq \aut(\cM)$, the \emph{symmetry type graph} of $\cM$ with respect to $\Gamma$  is the quotient $\cM/\Gamma$. In particular, when $\Gamma = \aut(\cM)$ then the quotient $\cM/\Gamma$ is the \emph{symmetry type graph} of $\cM$ (see \cite{CunninghamDelRioFrancosHubardToledo_2015_SymmetryTypeGraphs, Mochan_2021_AbstractPolytopesTheir_PhDThesis}).

The \emph{universal rank-$n$ Coxeter group} (with string diagram) is the group $\cC^{n}=\gen{\rho_{0}, \dots \rho_{n-1}}$ defined by the following relations:
\begin{equation}\label{eq:relsReg}
\begin{aligned}
\rho_{i}^{2} &= \epsilon && \text{for all } i \in \{0, \dots, n-1\}, \\
(\rho_{i}\rho_{j})^{2} &= \epsilon && \text{if } |i-j| \geq 2.
\end{aligned}
\end{equation}

The \emph{universal $n$-maniplex} $\cU^{n}$ is the Cayley graph associated with the universal Coxeter group $\cC^{n}$.
That is, the vertex set of $\cU^{n}$ is $\cC^{n}$ and for $\gamma \in \cC^{n}$, the $i$-adjacent vertex of $\gamma$ is $\rho_{i}\gamma$.
The maniplex $\cU^{n}$ is in fact the flag graph of the universal polytope (see \cite[Theorem 5.2]{Hartley_1999_AllPolytopesAre}).
We omit the rank of both the  universal Coxeter group and the universal maniplex whenever it is implicit.

Since the universal Coxeter group $\cC$ acts transitively by automorphisms on the universal maniplex $\cU$, then $\cU$ is a regular maniplex. In fact, $\aut(\cU)$ is precisely $\cC$.

We denote by $r_{i}$ the permutation of $V(\cU)$ that swaps every vertex $x$ with $x^{i}$.
More precisely, for $\gamma \in \cC = V(\cU)$, $r_i: \gamma \mapsto \rho_i\gamma $.

The \emph{monodromy group} of $\cU$, denoted by $\mon(\cU)$, is the permutation group of $V(\cU)$ generated by $\{r_{0}, \dots, r_{n-1}\}$. We call \emph{monodromies} the elements of $\mon(\cU)$
The permutations $r_0, \dots, r_{n-1}$ admit a (left) action on $\cX$ that is compatible with the (right) action of $\aut(\cU)$, that is \[(r_{i}x) \gamma = r_i(x\gamma)\] for every  $x \in V(\cU)$, $\gamma \in \aut(\cU)$ and $i \in \{0,\dots, n-1\}$. This induces an action of $\mon(\cU)$ on $V(\cU)$.

Observe that the elements $r_0, \dots, r_{n-1}$ satisfy the relations
\begin{equation}\label{eq:monU}
\begin{aligned}
r_{i}^{2} &= 1 && \text{for all } i \in \{0, \dots, n-1\}, \\
(r_{i}r_{j})^{2} &= 1 && \text{if } |i-j| \geq 2.
\end{aligned}
\end{equation}

In fact, there exists an isomorphism from $\cC$ to $\mon(\cU)$ mapping $\rho_{i}$ to $r_{i}$ (see \cite[Theorem 3.9]{MonsonPellicerWilliams_2014_MixingMonodromyAbstract}).

The relations in Equation\nobreakspace \textup {(\ref {eq:monU})} imply that if $\cX$ is any $n$-premaniplex, then $\mon(\cU)$ acts on $\cX$ by $r_i: x \mapsto x^{i}$ for every vertex $x$ and $i \in \{0, \dots, n-1 \}$.
The \emph{monodromy group} of $\cX$ is the permutation group induced by this action.
Naturally, we denote this group by $\mon(\cX)$.
It is important to remark that the term \emph{connection group} has been used for what we call the monodromy group and has gained popularity in the last few years.

Recall that for a vertex $x$ and $i \in \{0, \dots, n-1\}$, $\pth{x}$ denotes the dart of color $i$ whose starting point is $x$.
We can generalize this notation.
If $W$ is a path starting at $x$ and following the colors $i_{1}, i_{2}, \dots, i_{k}$ then we write
$W=\pth[i_{k}, i_{k-1}, \dots, i_{1} ]{x}$.
For every vertex $x$ and $i_{1}, i_{2}, \dots, i_{k} \in \left\{ 0, \dots, n-1 \right\} $, the path
$\pth[i_{k}, i_{k-1}, \dots, i_{1} ]{x}$ ends at the vertex $r_{i_{k}} \cdots r_{i_{1}}x$.

Given a vertex $u$ of $\cU$ and a path $W=\pth[i_k,i_{k-1}\ldots,i_1]{x}$ in a premaniplex $\cX$, the \emph{lift of $W$ starting at $u$} is simply the path $\pth[i_k,i_{k-1},\ldots,i_1]{u}$ in $\cU$.

An \emph{elementary (maniplex) move} on a path $\pth[w]{x}$ with $w=i_{k}, i_{k-1}, \dots, i_{1}$ consist of either adding or removing the same color two times at any two consecutive positions (i.e. if $v=i_k,\ldots,i_\ell,j,j,i_{\ell-1},\ldots,i_1$ then  $\pth[w]{x}\mapsto \pth[v]{x}$ and $\pth[v]{x}\mapsto \pth[w]{x}$ are elementary moves)
or swapping two non-consecutive colors in consecutive positions, more precisely, if $|i_\ell-i_{\ell-1}|>1$, and $u=i_k,\ldots,i_{\ell+1},i_{\ell-1},i_\ell,i_{\ell-2},\ldots,i_1$ then $\pth[w]{x}\mapsto \pth[u]{x}$ is an elementary move.
We say that two paths in a premaniplex are \emph{maniplex-homotopic} if we can turn one into the other by a finite sequence of elementary (maniplex) moves.
Observe that two paths in an premaniplex are homotopic if their lifts to the universal maniplex $\cU$ starting at a the same vertex also end at the same vertex.

The homotopy class of a path in a premaniplex is uniquely determined by its starting vertex and a monodromy of the universal maniplex.
If $w \in \mon(\cU^{n})$ and $w=r_{i_{k}} \cdots r_{i_{1}}$ is a word representing $w$, then we denote by $\pth[w]x$ the homotopy class of the path $\pth[i_{k}, i_{k-1}, \dots, i_{1} ]{x}$.
Observe that any two words representing $w$ yield homotopic paths and if $W_{1}, W_{2} \in \pth[w]{x}$, then both paths end at $wx$.

Let $\cX$ be a premaniplex and let $\fg(\cX)$ be its fundamental groupoid. If $\xi:\fg(\cX)\to \Gamma$ is a voltage assignment such that $\xi(W)$ is the identity in $\Gamma$ whenever $W$ is a path of length 4 alternating between two non-consecutive colors, we say that the pair $(\cX,\xi)$ is a \emph{voltage premaniplex}.
In other words, a voltage premaniplex is a voltage graph where the voltage of the maniplex homotopy class of a  path is well defined.
 \section{Voltage operations} \label{sec:operations}

Let $\cX$ be an $n$-premaniplex and let $\cY$ be an $m$-premaniplex.
Consider the voltage assignment $\eta: \fg(\cY) \to \mon(\cU^{n})$. We define the $m$-colored graph $\cX \ertimes \cY$ in the following way: the vertex set of $\cX \ertimes \cY$ is $\vr(\cX) \times \vr(\cY)$ and, for each $i\in \{0,1,\dots, m\}$, there is an edge of color $i$ from  $(x,y)$ to $\left(\eta(\pth{y})x,r_i y\right)$.
Since, in $\cY$, there is an $i$-edge from $y$ to $r_iy$, then the darts $\pth{y}$ and $\pth{r_iy}$ are inverse, implying that $\eta(\pth{r_iy})\eta(\pth{y})$ is the identity. Thus, each vertex of the graph $\cX \ertimes \cY$ has exactly one $i$-edge and, hence, it has degree $m$.

In the terminology of \cite{MalnicNedelaSkoviera_2000_LiftingGraphAutomorphisms}, $\cX\ertimes \cY$ is the derived graph from the \emph{voltage space} $(\cX,\mon(\cU);\eta)$ with $\cX$ being the abstract fiber.

A pair $\left( \cY,\eta \right)$ is called a \emph{voltage operator} if it is a voltage premaniplex with $\eta: \fg(\cY) \to \mon(\cU^{n})$ (or, if the rank of $\cY$ is $m$, an $\left( n,m \right)$-\emph{voltage operator}).
Similarly, we say that the mapping $\cX \mapsto \cX \ertimes \cY$, where $\cX$ runs over all $n$-premaniplexes and $(\cY,\eta)$ is a fixed voltage operator, is a \emph{voltage operation} (or an $(n,m)$-\emph{voltage operation}).

Given premaniplexes $\cX$ and $\cY$ and a voltage assignment $\eta$, we shall see that $\cX \ertimes \cY$ is a premaniplex itself, although need not be a maniplex (even if $\cX$ is a maniplex).
Before showing this, we give some straightforward examples of voltage operations.

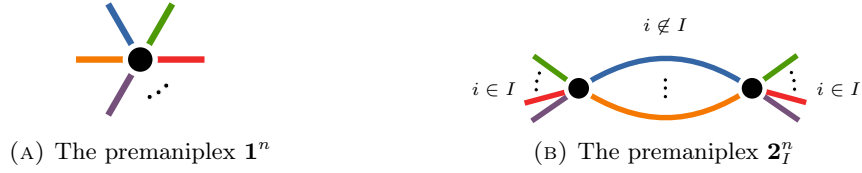
\begin{figure}\label{fig:EjemplosRango3}
	\centering
	\begin{subfigure}[b]{0.45\textwidth}
	\centering
        \includegraphics[width=.3\textwidth]{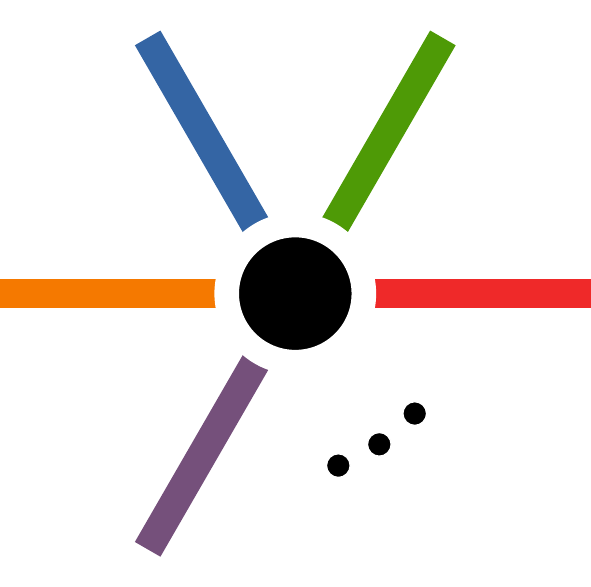}
 \caption{The premaniplex $\1^{n}$}
 \label{fig:1ton}
	\end{subfigure}
	\hfill \begin{subfigure}[b]{0.45\textwidth}
		\centering
			\begin{scriptsize}
\def\svgwidth{\textwidth}
		\begingroup \makeatletter \providecommand\color[2][]{\errmessage{(Inkscape) Color is used for the text in Inkscape, but the package 'color.sty' is not loaded}\renewcommand\color[2][]{}}\providecommand\transparent[1]{\errmessage{(Inkscape) Transparency is used (non-zero) for the text in Inkscape, but the package 'transparent.sty' is not loaded}\renewcommand\transparent[1]{}}\providecommand\rotatebox[2]{#2}\newcommand*\fsize{\dimexpr\f@size pt\relax}\newcommand*\lineheight[1]{\fontsize{\fsize}{#1\fsize}\selectfont}\ifx\svgwidth\undefined \setlength{\unitlength}{671.8930318bp}\ifx\svgscale\undefined \relax \else \setlength{\unitlength}{\unitlength * \real{\svgscale}}\fi \else \setlength{\unitlength}{\svgwidth}\fi \global\let\svgwidth\undefined \global\let\svgscale\undefined \makeatother \begin{picture}(1,0.24878184)\lineheight{1}\setlength\tabcolsep{0pt}\put(0,0){\includegraphics[width=\unitlength,page=1]{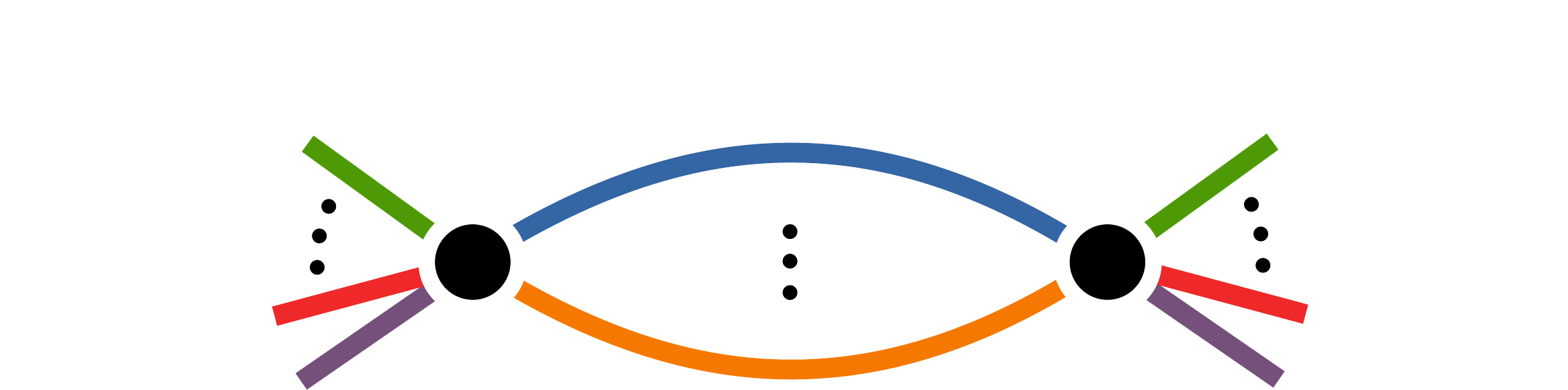}}\put(0.50397153,0.21573051){\color[rgb]{0,0,0}\makebox(0,0)[t]{\lineheight{1.25}\smash{\begin{tabular}[t]{c}$i \not\in I$\end{tabular}}}}\put(0.15588195,0.06644261){\color[rgb]{0,0,0}\makebox(0,0)[rt]{\lineheight{1.25}\smash{\begin{tabular}[t]{r}$i  \in I$\end{tabular}}}}\put(0.85676304,0.06644261){\color[rgb]{0,0,0}\makebox(0,0)[lt]{\lineheight{1.25}\smash{\begin{tabular}[t]{l}$i \in I$\end{tabular}}}}\end{picture}\endgroup  		\end{scriptsize}
		\caption{ The premaniplex $\2^n_{I}$}
		\label{fig:2_I}
	\end{subfigure}
	\caption{Premaniplexes with $1$ and $2$ vertices.}
\end{figure}

Let us denote by $\1^{n}$ the symmetry type graph of a regular $n$-maniplex.
That is, the premaniplex whit only one vertex and $n$ semiedges (see Figure\nobreakspace \ref {fig:1ton}).
Moreover, if $G$ is a group and $\eta: \1^{n} \to G$ is a voltage assignment that assigns $g_{i}$ to the semiedge of color $i$, then we denote this voltage graph by $\left( \1^{n}, [g_{0}, \dots, g_{n-1}] \right)$.

\begin{exam}\label{eg:d-auto} \leavevmode
The concept of a $d$-automorphism of a polytope is defined in \cite{HubardOrbanicIvicWeiss_2009_MonodromyGroupsSelf}, and can be generalized to maniplexes in a straightforward way.
Let ${\mathcal M}$ be an $n$-maniplex
and consider $\mon({\mathcal U}^n) =\langle r_0, r_1, \dots, r_{n-1} \rangle$ and $d$ an automorphism of $\mon(\cU^n)$.
We shall keep the vertices of $\cM$ and the action of $\mon(\cU^n)$ on them, but we choose a new set of labeled generators of the same permutation group acting on the same set of flags to obtain a new maniplex, ${\mathcal M}^d$.
More precisely, the maniplex ${\mathcal M}^d$ is defined as follows: the vertices of $\cM^d$ are the vertices of $\cM$ and given a vertex $x$, the dart of color $i$ starting at $x$ ends at $d(r_i)x$.
In other words, in $\cM^d$, $x^i = d(r_i)x$. If $\cM \cong \cM^d$, we call each isomorphism $\varphi:\cM\to\cM^d$ a \emph{$d$-automorphism} of $\cM$, and we say that $\cM$ is \emph{$d$-automorphic}.

The maniplex $\cM^d$ can easily be seen as an $(n,n)$-voltage operation:
\begin{center}
if $(\cY,\eta) = (\1^{n},[d(r_0), d(r_1) \dots, d(r_{n-1})])$, then $\cM \ertimes \cY = \cM^d$.
\end{center}
Classical examples of the above operation are the dual and the Petrial of an $n$-maniplex $\cX$:

	\begin{itemize}
	    \item If $(\cY,\eta) = (\1^{n},[r_0, \dots, r_{n-1}])$, then $\cX \ertimes \cY$ is $\cX$ itself.
		\item If $(\cY,\eta) = (\1^{n},[r_{n-1}, \dots, r_{0}])$, then $\cX \ertimes \cY$ is the dual of $\cX$.
		\item If $(\cY,\eta) = (\1^{n},[r_0, \dots, r_{n-4}, r_{n-3}r_{n-1}, r_{n-2}, r_{n-1}])$, then $\cX \ertimes \cY$ is the generalized Petrial of $\cX$.
\end{itemize}
In the above examples we know that whenever $\cX$ is in fact a maniplex, then so is $\cX \ertimes \cY$.
\end{exam}

If $(\cY, \eta)$ is a voltage operator, we want to show that $\cX \ertimes \cY$ is a premaniplex, for every premaniplex $\cX$.
Thus, we need to show that the alternating paths in $\cX \ertimes \cY$ of length $4$ with colors $i$ and $j$, with $\left| i-j \right| \geq 2 $, are closed.
In other words, given a premaniplex $\cX$ and a voltage operator $(\cY, \eta)$ we want $\pth[i,j,i,j]{(x,y)}=(x,y)$, whenever $\left| i-j \right| \geq 2 $, for every $x\in\cX$ and $y \in \cY$.

More generally, we are interested in obtaining properties for the paths in $\cX \ertimes \cY$ from the paths in $\cX$ or $\cY$.
The following result is a straightforward but useful observation towards this direction.

\begin{rem}\label{rem:endsPaths}
Let $(\cY, \eta)$ a voltage operator. Assume that $W=\pth[i_{k},\dots, i_{1}]{y}$ is a path in $\cY$ starting at $y$.
For every $x \in \cX$ the path $\pth[i_{k},\dots, i_{1}]{(x,y)}$ in $\cX \ertimes \cY$ that starts at $(x,y)$ and follows the same colors as $W$ finishes at $\left(\eta(W)x, r_{i_{k}} \cdots r_{i_{1}}y\right)$.
\end{rem}

By Remak\nobreakspace \ref {rem:endsPaths} the path $\pth[i,j,i,j]{(x,y)}$ starts at $(x,y)$ and finishes at
$\left(\eta(\pth[i,j,i,j]{y})x, r_ir_jr_ir_j y\right)$.
Thus, the path $\pth[i,j,i,j]{(x,y)}$ is closed if and only if
\[\begin{aligned}
    \eta(\pth[i,j,i,j]{y})x &= x && \text{and} \\
	r_{i}r_{j}r_{i}r_{j}y &= y,
  \end{aligned}\]
Note now that since $\cY$ is a premaniplex, then $r_{i}r_{j}r_{i}r_{j}y = y$ for every $y\in \cY$, which implies that $\pth[i,j,i,j]{(x,y)}$ is closed if and only if $\eta(\pth[i,j,i,j]{y})$ fixes $x$. But $\eta(\pth[i,j,i,j]{y})$ is the identity, since $(\cY, \eta)$ is a voltage premaniplex. Therefore, we have the following proposition.

\begin{prop}
Given a voltage premaniplex $(\cY, \eta)$ and a premaniplex $\cX$, the voltage operation $\cX \ertimes \cY$ is a premaniplex.
\end{prop}

Since we are working with voltage premaniplexes, for the rest of the paper, whenever we refer to paths being ``homotopic'' or to the ``homotopy class'' of a path, we are thinking in terms of maniplex homotopy. In the same way, the notation $\fg(\cY)$ will denote the fundamental groupoid consisting of paths in $\cY$ considered up to maniplex homotopy.

Although we have special interest on voltage operations that give rise to connected premaniplexes, there are interesting examples of voltage operations in which we (often) obtain disconnected objects.

A \emph{rooted} premaniplex is a pair $(\cX,x)$ where $\cX$ is a connected premaniplex and the \emph{root} $x$ is a vertex of $\cX$.
If $\cX$ is not connected, then $(\cX,x)$ denotes the rooted premaniplex with the connected component of $x$ in $\cX$ as its underlying premaniplex.
If $(\cX,x)$ and $(\cY,y)$ are rooted premaniplexes and $(\cY,\eta)$ is a voltage operator, then $(\cX,x)\ertimes(\cY,y)$ denotes the rooted premaniplex $(\cX\ertimes \cY, (x,y))$.

\begin{exam}\label{eg:sections} \leavevmode
Let $\cX$ be a $n$-premaniplex, $-1 \leq k < \ell \leq n$ and let $(\cY,\eta) = (\1^{\ell - k - 1},[r_{k+1}, \dots, r_{\ell-1}])$. Then $\cX \ertimes \cY$ is a graph whose connected components are the $(k,\ell)$-sections of $\cX$.
    In particular, if $k=-1$ and $\cX$ is a maniplex, then $\cX \ertimes \cY$ determines the set of all $\ell$-faces of $\cX$, and if $(\cX,x)$ is a rooted premaniplex.
    Moreover, if $y$ denotes the only vertex in $\cY$, for each $x \in \cX$, then $(\cX,x)\ertimes (\cY,y)$  is the $\ell$-face of $\cX$ that contains $x$.
\end{exam}

Note that if $\cY$ is disconnected, then $\cX\ertimes\cY$  consists of disjoint copies of $\cX\rtimes_{\eta_i}\cY_i$, where $\cY_i$ runs over the connected components of $\cY$ and $\eta_i$ is the restriction of $\eta$ to $\fg(\cY_i$).

In a recent paper \cite{CunninghamPellicerWilliams_StratifiedOperationsManiplexes_preprint} the term \emph{stratified operations} is introduced. They are operations on maniplexes defined in terms of a set $A$ of \emph{strata} admitting a left action of $\mon(\cU)$ and a function $\varphi$ that assigns a monodromy to each pair $(a,r_i)$ where $r_i$ is a generator of $\mon(\cU)$.
We can describe stratified operations in terms of voltage operations. If $\oo$ is a stratified operation with strata set $A$, we can define the voltage operator $(\cY,\eta)$ by taking the vertices of $\cY$ to be the set $A$, defining the adjacencies by $a^i=r_i a$, and defining the voltage assignment $\eta$ by $\eta(\pth{a})=\varphi(a,r_i)$.
There may be a subtle difference between the stratified operation $\oo$ and the voltage operation defined by $(\cY,\eta)$. Mainly, if $\cM$ is a maniplex $\oo(\cM)$ must be a maniplex too, and therefore connected, while $\cM\ertimes\cY$ might have more than one connected component.
One can check that all stratified operations are the result of applying a voltage operation and then choosing one connected component. The converse, however is not true. There are voltage operations that do not define stratified operations when one chooses a connected component of the result.
This is because the definition of stratified operation requires that the natural projection $\oo(\cM)\to\cM$ must be surjective.
 The snub operation described in Section\nobreakspace \ref {sec:ClassicalExamples} and discussed further in Section\nobreakspace \ref {sec:automorphisms} is a voltage operation, but not a stratified operation, since each connected component only uses half of the flags of the original maniplex when this is orientable.

\subsection{The mix as a voltage operation on premaniplexes}\label{sec:mix}

The mix of two regular polytopes was defined in terms of their automorphism groups (see \cite{McMullenSchulte_2002_MixRegularPolytope, MonsonPellicerWilliams_2014_MixingMonodromyAbstract}).
In \cite{CunninghamPellicer_2018_OpenProblems$k$} the mix is defined for rooted maniplexes as a natural generalization of the parallel product of rooted maps (see \cite{Wilson_1994_ParallelProductsGroups}), and hence it is also defined of rooted polytopes.
If $(\cM, \Phi)$ and $(\cN, \Psi)$ are rooted maniplexes, then $(\cM,\Phi) \mix (\cN,\Psi)$ is the smallest maniplex that covers both $(\cM, \Phi)$ and $ (\cN, \Psi)$. In particular, it is connected.
We shall see that the mix can be defined as a voltage operation on premaniplexes.

Let $\cY$ be an $n$-premaniplex. We denote by $\mixer$ the voltage assignment that maps each dart of color $i$ to $r_{i}$.
In other words $\mixer(\pth[\omega]y):=\omega$.
We call $\mixer$ the \emph{mixing voltage (for $\cY$)} and we say that $(\cY,\mixer)$ is a \emph{mix operator}.
Then $\cX \mixtimes \cY$ is the \emph{mix} $\cX \mix \cY$, and it is easy to see that this generalizes the same concept previously defined for regular abstract polytopes. By rooting the premaniplexes, this voltage operation generalizes the one for rooted maniplexes.

With our definition of the mix as a voltage operation we allow the resulting graph to be disconnected.
However, each connected component of $\cX \mixtimes \cY$ covers both $\cX$ and $\cY$, so if at least one of them is simple (that is, a maniplex) the mentioned components are maniplexes as well. Moreover, the connected component of $\cX\mix\cY$ containing the vertex $(x,y)$ is the smallest premaniplex that covers both $(\cX,x)$ and $(\cY,y)$ (cf. \cite[Proposition 3.10]{CunninghamPellicer_2018_OpenProblems$k$}).

We can use this definition of the mix to find the smallest cover of a maniplex satisfying some property, as we shall see in the next example.

\begin{eg}\label{eg:DoubleCovers}
 Let $I\subset\{0,1,\ldots,n-1\}$. We denote by $\2^n_I$ the $n$-premaniplex with 2 vertices with semiedges of colors in $I$ at each vertex and links of colors not in $I$ between the $2$ vertices (see Figure\nobreakspace \ref {fig:2_I}).
 We can use these premaniplexes to describe certain properties of maniplexes: \emph{orientable} $n$-maniplexes (those that are bipartite) are those that cover $\2^n_\emptyset$ and \emph{vertex-bipartite} maniplexes (those whose $1$-skeleton is bipartite) are those that cover $\2^n_{\{1,2,\ldots,n-1\}}$, for example.
 In general, if $\cX$ covers $\2^n_I$ it means that there is a coloring of the vertices (or flags) of $\cX$ with two colors such that $i$-adjacent flags are of the same color if and only if $i\in I$. See \cite{KoikePellicerRaggiWilson_2017_FlagBicoloringsPseudo, PellicerPotocnikToledo_2019_ExistenceResultTwo} for a detailed discussion on flag colorings.

 If a maniplex $\cX$ does not cover $\2^n_I$, then $\cX\mix \2^n_I$ is a maniplex that covers $\cX$ but also covers $\2^n_I$. We call this the \emph{double cover of $\cX$ with respect to $I$}.
 If $\cX$ is non-orientable, $\cX\mix \2^n_\emptyset$ is the so called \emph{orientable double cover of $\cX$}.

Note that if $\cX$ does cover $\2^n_I$, then $\cX\mix \2^n_I = \cX \mixtimes \2^n_I$ (with $\mixer$ the mixing voltage) consists of two isomorphic copies of $\cX$, but the flags of these copies will be colored by the vertices of $\2^n_I$.
More precisely, if we colored the vertices of $\2^n_I$ with white and black, then each vertex (flag) of $\cX\mix \2^n_I$ is colored white or black according to its second coordinate.
So $\cX\mix \2^n_I$ consists of two copies of $\cX$ but together with an $I$-compatible coloring, and the two copies have opposite colorings.
In particular, if $X$ is orientable and $I=\emptyset$, the two copies of $\cX \mix \2^{n}_{\emptyset}$ are mirror images; if $\cX$ is a chiral maniplex, then $\cX \mixtimes \2^{n}_{\emptyset}$ consists of its two enantiomorphic forms.
\end{eg}

      Before ending this section let us note that the mixing voltages are the only voltages such that the product is naturally commutative. More precisely:
      \begin{prop}\label{prop:MixConmuta}
        Let $(\cX,\eta_\cX)$ and $(\cY,\eta_\cY)$ be $(n,n)$-voltage operators.
        Then the function $(x,y)\mapsto (y,x)$ is an isomorphism between $\cX\ertimes[\eta_\cY]\cY$ and $\cY\ertimes[\eta_\cX]\cX$ if and only if both $\eta_\cX$ and $\eta_\cY$ are mixing voltages.
      \end{prop}
      \begin{proof}
        The alleged function is an isomorphism if and only if it maps $(x,y)^i = (\eta_\cY(\pth{y})x,y^i)$ to $(y,x)^i = (\eta_\cX(\pth{x})y,x^i)$, for all $x\in\cX$, $y\in\cY$ and $i\in\{0,1,\ldots,n-1\}$.
          This means that $\eta_\cY(\pth{y}) = r_i$ and $r_i = \eta_\cX(\pth{x})$.
      \end{proof}

      \begin{question}
        Is it possible that $\cX\ertimes[\eta_\cY]\cY$ and $\cY\ertimes[\eta_\cX]\cX$ are isomorphic without $\eta_\cX$ and $\eta_\cY$ being mixing voltages?
      \end{question}

\subsection{Connectivity}\label{sec:connectivity}
We now turn our attention to the connectivity of $\cX\ertimes \cY$.
We say that an $(n,m)$-voltage operator $(\cY,\eta)$  \emph{preserves connectivity} if whenever $\cX$ is connected $\cX\ertimes \cY$ is connected as well (in the context of \cite{CunninghamPellicerWilliams_StratifiedOperationsManiplexes_preprint}, these are called \emph{fully stratified operations}).
We first analyze when we can find a path between two vertices of $\cX \ertimes \cY$ and use this to determine when a voltage operator preserves connectivity.

\begin{lemma}\label{prop:paths}
	Let $(x,y)$ and $(x',y')$ two vertices of $\cX \ertimes \cY$. There is a path in $\cX \ertimes \cY$ that starts at $(x,y)$ and ends at $(x',y')$ if and only if there exists a path $W$ in $\cY$ that connects $y$ with $y'$ and such that $\eta(W)x= x'$.
\end{lemma}
\begin{proof}
	Assume there is a path $\pth[i_{k},\dots,i_{1}]{(x,y)}$ that ends at $(x',y')$. But by Remak\nobreakspace \ref {rem:endsPaths} the path $\pth[i_{k}, \dots, i_{1}]{(x,y)}$ ends at $\left( \eta(W)x, r_{i_{k}}\cdots r_{i_{1}}y\right)$, where
	 $W$ denotes the path $\pth[i_{k}, \dots, i_{1}]{y}$.
	Hence, $\eta(W)x = x'$ and $r_{i_{k}}\cdots r_{i_{1}}y = y'$.
	Moreover, $r_{i_{k}}\cdots r_{i_{1}}y$ is precisely the endpoint of $W$, which implies that $W$ connects $y$ with $y'$.

	Conversely, assume that there exists a path $W$ that connects $y$ with $y'$ and satisfies that $\eta(W)x = x'$.
	Since $W$ starts at $y$, it can be written as $W=\pth[i_{k}, \dots, i_{1}]{y}$, for some
	 $i_{k}, \dots, i_{1} \in \left\{ 0, \dots, m-1 \right\}$.
	This implies that $y' = r_{i_{k}}\cdots r_{i_{1}}y$. By Remak\nobreakspace \ref {rem:endsPaths}, the path $\pth[i_{k}, \dots, i_{1}]{(x,y)}$ ends at $\left( \eta(W)x, r_{i_{k}}\cdots r_{i_{1}}y\right) = (x',y')$.

\end{proof}

\begin{prop} \label{prop:connectedness}
The graph $\cX \ertimes \cY$ is connected if and only if $\cY$ is connected and $\eta(\fg^{y_0}(\cY))$ acts transitively on $\cX$ for some vertex $y_0 \in \cY$.
\end{prop}

\begin{proof}
    We start by assuming that $\cX \ertimes\cY$ is connected, and let $x,x'\in \cX$.
    First note that $\cX\ertimes\cY$ covers $\cY$, so $\cY$ must be connected as well.
    Now we want to prove that $\eta(\fg^{y_0}(\cY))$ acts transitively on $\cX$.
    Since $\cX \ertimes\cY$ is connected, there is a path from $(x,y_0)$ to $(x',y_0)$.
    By Lemma~\ref{prop:paths} there exists $W\in \fg^{y_0}(\cY)$ such that $\eta(W)x=x'$.
    Since $x$ and $x'$ were arbitrary, we have proved that $\eta(\fg^{y_0}(\cY))$ acts transitively on $\cX$.

    Now we assume that $\cY$ is connected and $\eta(\fg^{y_0}(\cY))$ acts transitively on $\cX$. We will prove that $\cX \ertimes \cY$ is connected.
    Let $x_0 \in \cX$ be fixed, and let $x\in\cX$ and $y\in \cY$.
    Since $\cY$ is connected, there is a path $V$ from $y_0$ to $y$ in $\cY$. Let $\sigma=\eta(V)$.
    Since $\eta(\fg^{y_0}(\cY))$ acts transitively on $\cX$, there exists $W \in \fg^{y_0}(\cY)$ such that $\eta(W)x_0=\sigma^{-1}x$. Then the path $WV$ goes from $y_0$ to $y$ and satisfies that $\eta(WV)x_0=x$.
    By Lemma\nobreakspace \ref {prop:paths} the existence of the path $WV$ implies that there exists a path in $\cX \ertimes \cY$ from $(x_0,y_0)$ to $(x,y)$. That is, $(x_0,y_0)$ is connected to every  vertex of $\cX \ertimes \cY$, implying that $\cX \ertimes \cY$ is connected.
\end{proof}

\begin{coro}\label{coro:connected}
    A voltage operator $(\cY,\eta)$ preserves connectivity if and only if $\cY$ is connected and $\eta(\fg^{y_0}(\cY)) = \mon(\cU)$ for some $y_0\in \cY$.
\end{coro}

In \cite{CunninghamDelRioFrancosHubardToledo_2015_SymmetryTypeGraphs} the authors give a way to find generators of the automorphism group of a maniplex, given its symmetry type graph  (STG). They do so they in terms of a spanning tree of the STG with trivial voltage in all its darts. Hence, it will be helpful to see Corollary\nobreakspace \ref {coro:connected} in the light of such a spanning tree of $\cY$.

\begin{coro}
    Let $(\cY,\eta)$ be a voltage operator. Assume that $\cY$ is connected and has a spanning tree $T$ such that all arcs in $T$ have trivial voltage.
    Then $(\cY,\eta)$ preserves connectivity if and only if $\mon(\cU)$ is generated by the voltages of the darts in $\cY$ not in $T$.
\end{coro}
\begin{proof}
    The group $\fg^{y_0}(\cY)$ is generated by the paths $C_d$ where $d$ varies among darts in $\cY$ not in $T$ and $C_d$ denotes the path that starts in $y_0$, goes to the initial vertex of $d$ through $T$, then uses the dart $d$, and then goes back to $y_0$ through $T$.
    This implies that $\eta(\fg^{y_0}(\cY))$ is generated by the voltages of those paths, but $\eta(C_d)=\eta(d)$.
\end{proof}

\section{More examples of voltage operations} \label{sec:ClassicalExamples}
Here, we shall see some examples of operations on maps and polytopes that can be seen as voltage operations.
In fact, if we have an operation on polytopes (or maps) that can be seen as a voltage operation, by defining the corresponding voltage operator we have defined the operation on maniplexes and premaniplexes.  

\subsubsection*{Wythoffian constructions}
	  We start by recalling the Wythoffian constructions from convex regular polytopes (see for example \cite{coxeter1999beauty, schulte2016wythoffian}).
	  Given a regular convex polytope, its symmetry group is generated by reflections $\{\rho_0, \dots, \rho_{n-1}\}$.
	  We call these reflections the {\em generating reflections} of the group and note that their mirrors form a cone that, intersected with the polytope, results on a fundamental region for the polytope with respect to its symmetry group.
	  For a Wythoffian construction, first choose a non-empty set $A$ of generating reflections and pick a point $v$ of the fundamental region that is not fixed by any of the reflections in $A$, but fixed by any generating reflection not in $A$ (if they exist).
	  The vertices of the new convex polytope are the orbit of $v$ under $\gen{\rho_0, \dots \rho_{n-1}}$ and the polytope is the convex hull of such vertices.
	  This can be represented on the Coxeter diagram (a graph with $n$ nodes, each representing a generating reflection, and an edge  labeled by the order of the product of their nodes, whenever it is bigger than $2$), and marking the nodes representing the generating reflections in $A$.

	  If $n=3$, then $A$ is a non-empty subset of a $3$ element set, so there are seven possibilities for $A$, each giving a different construction. In particular, for $A = \{\rho_0\}$ and $A = \{\rho_2\}$ the Wythoffian constructions give rise to the identity and dual operations, respectively, that have been given as voltage operations in Example\nobreakspace \ref {eg:d-auto}.
	  It is not difficult to generalize the  Wythoffian constructions of regular polyhedra to maps on surfaces, by investigating the flag adjacencies of the constructed polyhedra.
	  In particular, the medial and truncation of a map have been studied in \cite{ hubard2013medial, HubardOrbanicIvicWeiss_2009_MonodromyGroupsSelf} and \cite{OrbanicPellicerWeiss_2010_MapOperations$k$}, respectively.

\begin{figure}
	\centering
	\begin{subfigure}[b]{0.45\textwidth}
		\centering
			\begin{scriptsize}
\def\svgwidth{\textwidth}
		\begingroup \makeatletter \providecommand\color[2][]{\errmessage{(Inkscape) Color is used for the text in Inkscape, but the package 'color.sty' is not loaded}\renewcommand\color[2][]{}}\providecommand\transparent[1]{\errmessage{(Inkscape) Transparency is used (non-zero) for the text in Inkscape, but the package 'transparent.sty' is not loaded}\renewcommand\transparent[1]{}}\providecommand\rotatebox[2]{#2}\newcommand*\fsize{\dimexpr\f@size pt\relax}\newcommand*\lineheight[1]{\fontsize{\fsize}{#1\fsize}\selectfont}\ifx\svgwidth\undefined \setlength{\unitlength}{793.821464bp}\ifx\svgscale\undefined \relax \else \setlength{\unitlength}{\unitlength * \real{\svgscale}}\fi \else \setlength{\unitlength}{\svgwidth}\fi \global\let\svgwidth\undefined \global\let\svgscale\undefined \makeatother \begin{picture}(1,0.66631666)\lineheight{1}\setlength\tabcolsep{0pt}\put(0,0){\includegraphics[width=\unitlength,page=1]{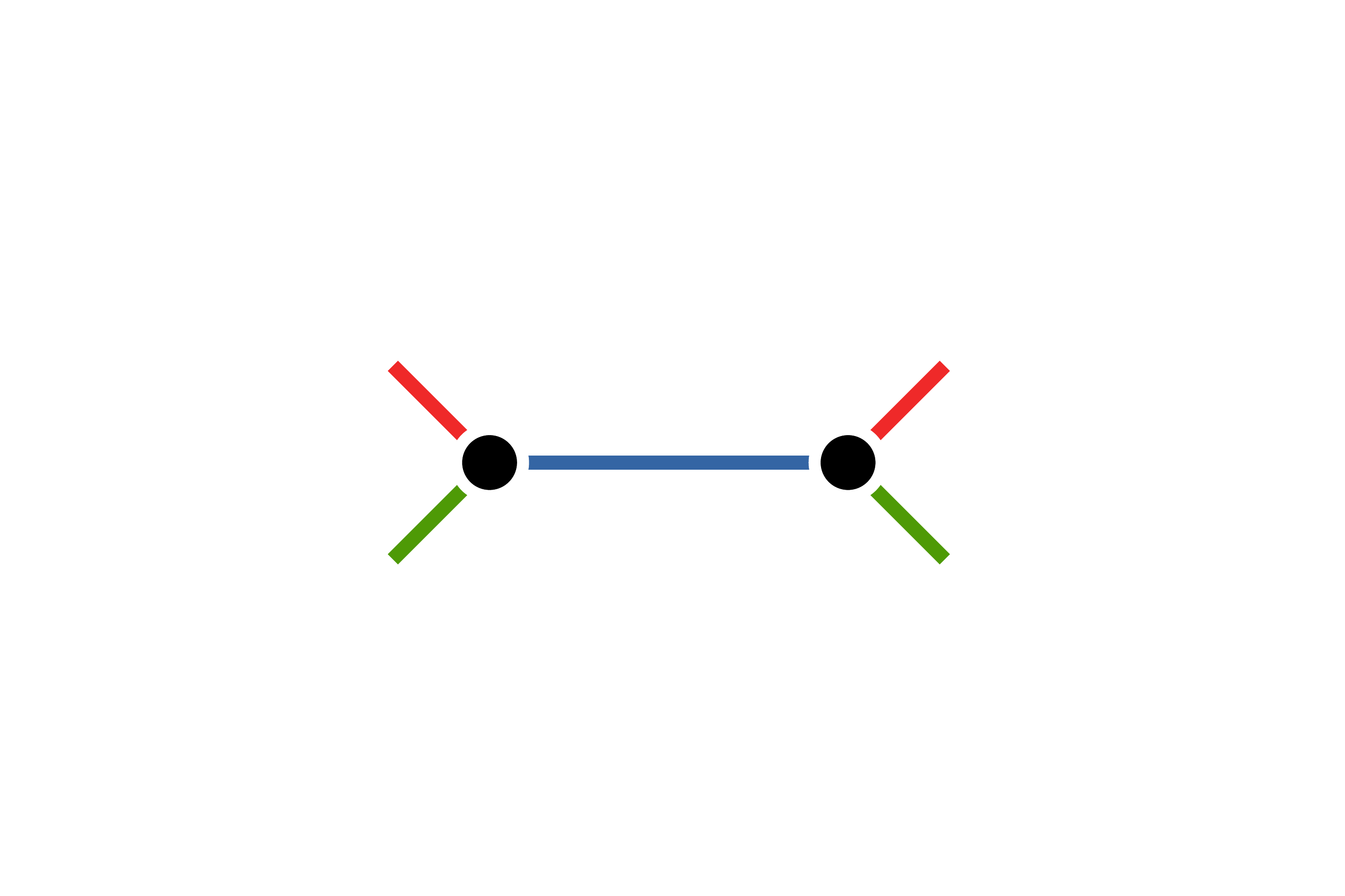}}\put(0.71790213,0.39940349){\color[rgb]{0,0,0}\makebox(0,0)[lt]{\lineheight{1.25}\smash{\begin{tabular}[t]{l}$r_1$\end{tabular}}}}\put(0.29127845,0.40281965){\color[rgb]{0,0,0}\makebox(0,0)[rt]{\lineheight{1.25}\smash{\begin{tabular}[t]{r}$r_1$\end{tabular}}}}\put(0.28027205,0.24267557){\color[rgb]{0,0,0}\makebox(0,0)[rt]{\lineheight{1.25}\smash{\begin{tabular}[t]{r}$r_0$\end{tabular}}}}\put(0.71993091,0.24593291){\color[rgb]{0,0,0}\makebox(0,0)[lt]{\lineheight{1.25}\smash{\begin{tabular}[t]{l}$r_2$\end{tabular}}}}\end{picture}\endgroup  		\end{scriptsize}
		\caption{\label{fig:medial} Medial operator $A=\left\{ \rho_{1} \right\} $}
	\end{subfigure}
	\hfill \begin{subfigure}[b]{0.45\textwidth}
		\centering
			\begin{scriptsize}
\def\svgwidth{\textwidth}
		\begingroup \makeatletter \providecommand\color[2][]{\errmessage{(Inkscape) Color is used for the text in Inkscape, but the package 'color.sty' is not loaded}\renewcommand\color[2][]{}}\providecommand\transparent[1]{\errmessage{(Inkscape) Transparency is used (non-zero) for the text in Inkscape, but the package 'transparent.sty' is not loaded}\renewcommand\transparent[1]{}}\providecommand\rotatebox[2]{#2}\newcommand*\fsize{\dimexpr\f@size pt\relax}\newcommand*\lineheight[1]{\fontsize{\fsize}{#1\fsize}\selectfont}\ifx\svgwidth\undefined \setlength{\unitlength}{793.821464bp}\ifx\svgscale\undefined \relax \else \setlength{\unitlength}{\unitlength * \real{\svgscale}}\fi \else \setlength{\unitlength}{\svgwidth}\fi \global\let\svgwidth\undefined \global\let\svgscale\undefined \makeatother \begin{picture}(1,0.66631666)\lineheight{1}\setlength\tabcolsep{0pt}\put(0,0){\includegraphics[width=\unitlength,page=1]{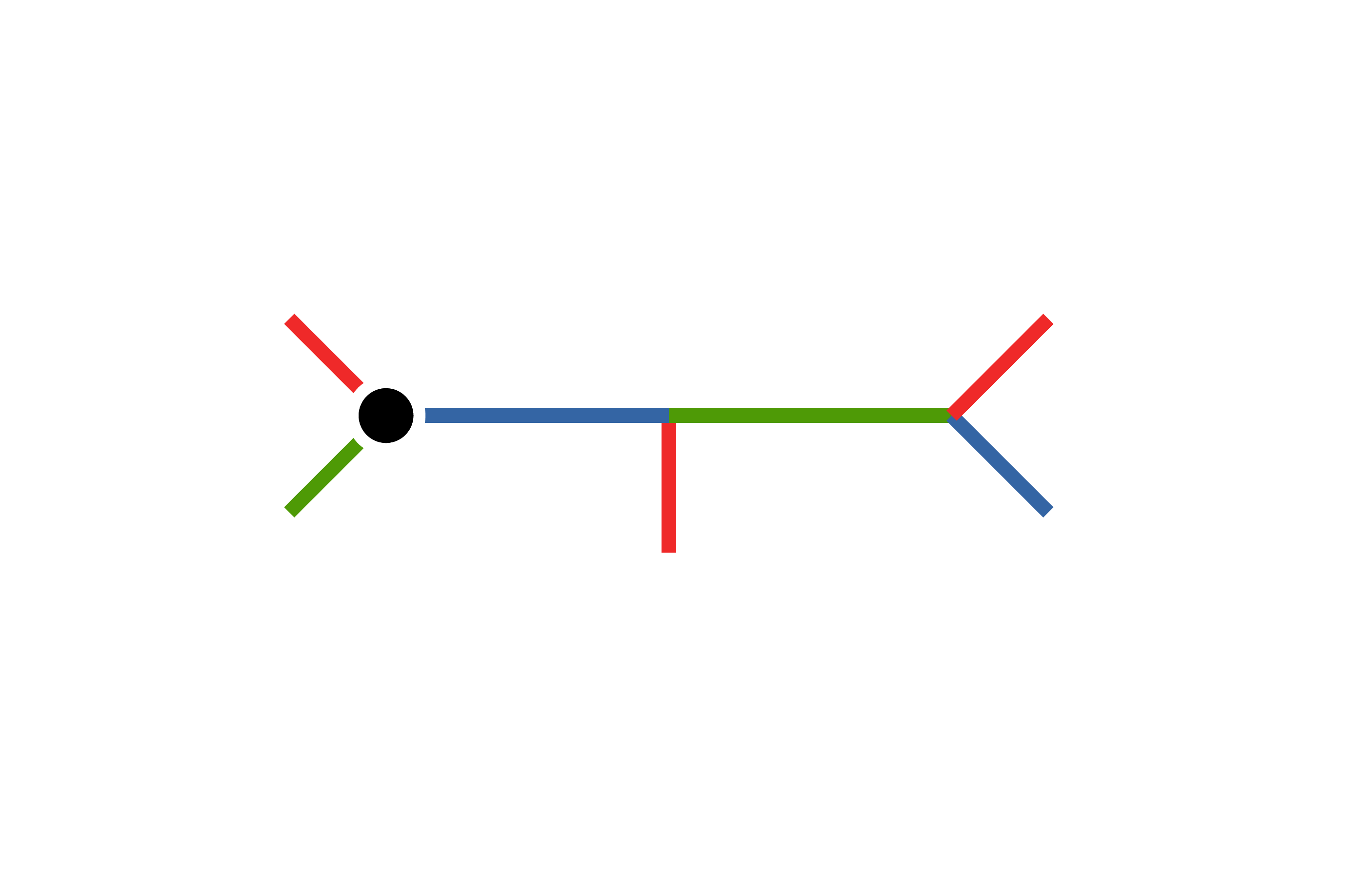}}\put(0.49682375,0.2059241){\color[rgb]{0,0,0}\makebox(0,0)[t]{\lineheight{1.25}\smash{\begin{tabular}[t]{c}$r_1$\end{tabular}}}}\put(0,0){\includegraphics[width=\unitlength,page=2]{trunc.pdf}}\put(0.80012247,0.43957113){\color[rgb]{0,0,0}\makebox(0,0)[lt]{\lineheight{1.25}\smash{\begin{tabular}[t]{l}$r_0$\end{tabular}}}}\put(0.20008051,0.43957113){\color[rgb]{0,0,0}\makebox(0,0)[rt]{\lineheight{1.25}\smash{\begin{tabular}[t]{r}$r_1$\end{tabular}}}}\put(0.80012247,0.25157233){\color[rgb]{0,0,0}\makebox(0,0)[lt]{\lineheight{1.25}\smash{\begin{tabular}[t]{l}$r_2$\end{tabular}}}}\put(0.20008051,0.25157233){\color[rgb]{0,0,0}\makebox(0,0)[rt]{\lineheight{1.25}\smash{\begin{tabular}[t]{r}$r_2$\end{tabular}}}}\end{picture}\endgroup  		\end{scriptsize}
		\caption{\label{fig:truncation} Truncation operator $A=\left\{ \rho_{0}, \rho_{1} \right\}$}
	\end{subfigure}

		\begin{subfigure}[b]{0.45\textwidth}
		\centering
			\begin{scriptsize}
\def\svgwidth{\textwidth}
		\begingroup \makeatletter \providecommand\color[2][]{\errmessage{(Inkscape) Color is used for the text in Inkscape, but the package 'color.sty' is not loaded}\renewcommand\color[2][]{}}\providecommand\transparent[1]{\errmessage{(Inkscape) Transparency is used (non-zero) for the text in Inkscape, but the package 'transparent.sty' is not loaded}\renewcommand\transparent[1]{}}\providecommand\rotatebox[2]{#2}\newcommand*\fsize{\dimexpr\f@size pt\relax}\newcommand*\lineheight[1]{\fontsize{\fsize}{#1\fsize}\selectfont}\ifx\svgwidth\undefined \setlength{\unitlength}{793.821464bp}\ifx\svgscale\undefined \relax \else \setlength{\unitlength}{\unitlength * \real{\svgscale}}\fi \else \setlength{\unitlength}{\svgwidth}\fi \global\let\svgwidth\undefined \global\let\svgscale\undefined \makeatother \begin{picture}(1,0.66631666)\lineheight{1}\setlength\tabcolsep{0pt}\put(0,0){\includegraphics[width=\unitlength,page=1]{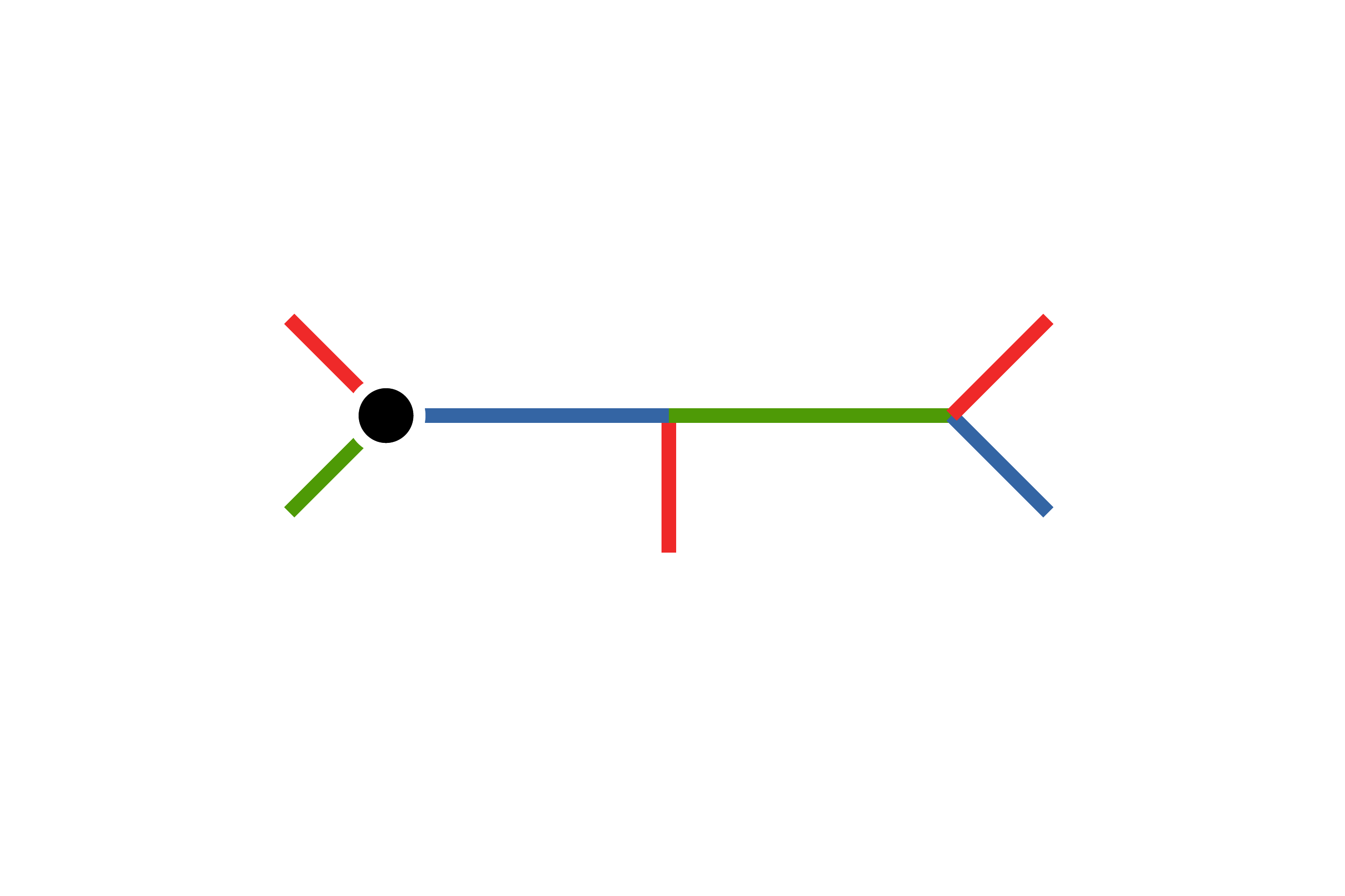}}\put(0.49682377,0.20592411){\color[rgb]{0,0,0}\makebox(0,0)[t]{\lineheight{1.25}\smash{\begin{tabular}[t]{c}$r_1$\end{tabular}}}}\put(0,0){\includegraphics[width=\unitlength,page=2]{trunc_dual.pdf}}\put(0.80012249,0.43957111){\color[rgb]{0,0,0}\makebox(0,0)[lt]{\lineheight{1.25}\smash{\begin{tabular}[t]{l}$r_2$\end{tabular}}}}\put(0.20008048,0.43957111){\color[rgb]{0,0,0}\makebox(0,0)[rt]{\lineheight{1.25}\smash{\begin{tabular}[t]{r}$r_1$\end{tabular}}}}\put(0.80012249,0.25157234){\color[rgb]{0,0,0}\makebox(0,0)[lt]{\lineheight{1.25}\smash{\begin{tabular}[t]{l}$r_0$\end{tabular}}}}\put(0.20008048,0.25157234){\color[rgb]{0,0,0}\makebox(0,0)[rt]{\lineheight{1.25}\smash{\begin{tabular}[t]{r}$r_0$\end{tabular}}}}\end{picture}\endgroup  		\end{scriptsize}
		\caption{\label{fig:trunc_dual} $A=\left\{ \rho_{1}, \rho_{2} \right\}$}
	\end{subfigure}
	\hfill \begin{subfigure}[b]{0.45\textwidth}
		\centering
			\begin{scriptsize}
\def\svgwidth{\textwidth}
		\begingroup \makeatletter \providecommand\color[2][]{\errmessage{(Inkscape) Color is used for the text in Inkscape, but the package 'color.sty' is not loaded}\renewcommand\color[2][]{}}\providecommand\transparent[1]{\errmessage{(Inkscape) Transparency is used (non-zero) for the text in Inkscape, but the package 'transparent.sty' is not loaded}\renewcommand\transparent[1]{}}\providecommand\rotatebox[2]{#2}\newcommand*\fsize{\dimexpr\f@size pt\relax}\newcommand*\lineheight[1]{\fontsize{\fsize}{#1\fsize}\selectfont}\ifx\svgwidth\undefined \setlength{\unitlength}{793.821464bp}\ifx\svgscale\undefined \relax \else \setlength{\unitlength}{\unitlength * \real{\svgscale}}\fi \else \setlength{\unitlength}{\svgwidth}\fi \global\let\svgwidth\undefined \global\let\svgscale\undefined \makeatother \begin{picture}(1,0.66631666)\lineheight{1}\setlength\tabcolsep{0pt}\put(0,0){\includegraphics[width=\unitlength,page=1]{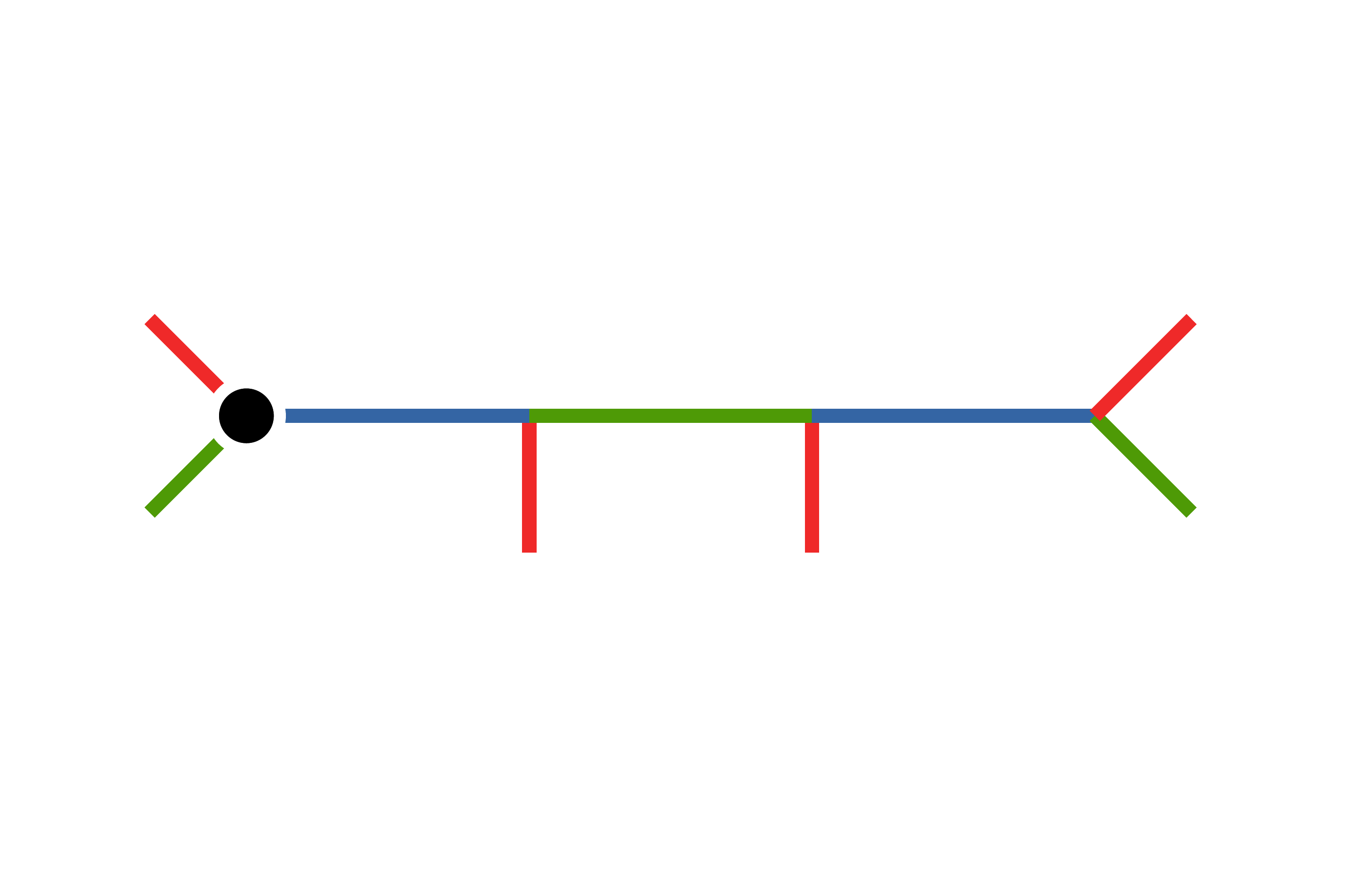}}\put(0.39299411,0.20574558){\color[rgb]{0,0,0}\makebox(0,0)[t]{\lineheight{1.25}\smash{\begin{tabular}[t]{c}$r_0$\end{tabular}}}}\put(0,0){\includegraphics[width=\unitlength,page=2]{whyt_02.pdf}}\put(0.90395212,0.43939261){\color[rgb]{0,0,0}\makebox(0,0)[lt]{\lineheight{1.25}\smash{\begin{tabular}[t]{l}$r_2$\end{tabular}}}}\put(0.09625087,0.43939261){\color[rgb]{0,0,0}\makebox(0,0)[rt]{\lineheight{1.25}\smash{\begin{tabular}[t]{r}$r_0$\end{tabular}}}}\put(0.89614996,0.25139386){\color[rgb]{0,0,0}\makebox(0,0)[lt]{\lineheight{1.25}\smash{\begin{tabular}[t]{l}$r_1$\end{tabular}}}}\put(0.09625087,0.25139386){\color[rgb]{0,0,0}\makebox(0,0)[rt]{\lineheight{1.25}\smash{\begin{tabular}[t]{r}$r_1$\end{tabular}}}}\put(0,0){\includegraphics[width=\unitlength,page=3]{whyt_02.pdf}}\put(0.60330377,0.20574553){\color[rgb]{0,0,0}\makebox(0,0)[t]{\lineheight{1.25}\smash{\begin{tabular}[t]{c}$r_2$\end{tabular}}}}\end{picture}\endgroup  		\end{scriptsize}
		\caption{ $A=\left\{ \rho_{0},\rho_{2} \right\}$}\label{fig:whyt_02}
	\end{subfigure}

			\begin{subfigure}[b]{0.45\textwidth}
		\centering
			\begin{scriptsize}
\def\svgwidth{\textwidth}
		\begingroup \makeatletter \providecommand\color[2][]{\errmessage{(Inkscape) Color is used for the text in Inkscape, but the package 'color.sty' is not loaded}\renewcommand\color[2][]{}}\providecommand\transparent[1]{\errmessage{(Inkscape) Transparency is used (non-zero) for the text in Inkscape, but the package 'transparent.sty' is not loaded}\renewcommand\transparent[1]{}}\providecommand\rotatebox[2]{#2}\newcommand*\fsize{\dimexpr\f@size pt\relax}\newcommand*\lineheight[1]{\fontsize{\fsize}{#1\fsize}\selectfont}\ifx\svgwidth\undefined \setlength{\unitlength}{793.821464bp}\ifx\svgscale\undefined \relax \else \setlength{\unitlength}{\unitlength * \real{\svgscale}}\fi \else \setlength{\unitlength}{\svgwidth}\fi \global\let\svgwidth\undefined \global\let\svgscale\undefined \makeatother \begin{picture}(1,0.66631666)\lineheight{1}\setlength\tabcolsep{0pt}\put(0,0){\includegraphics[width=\unitlength,page=1]{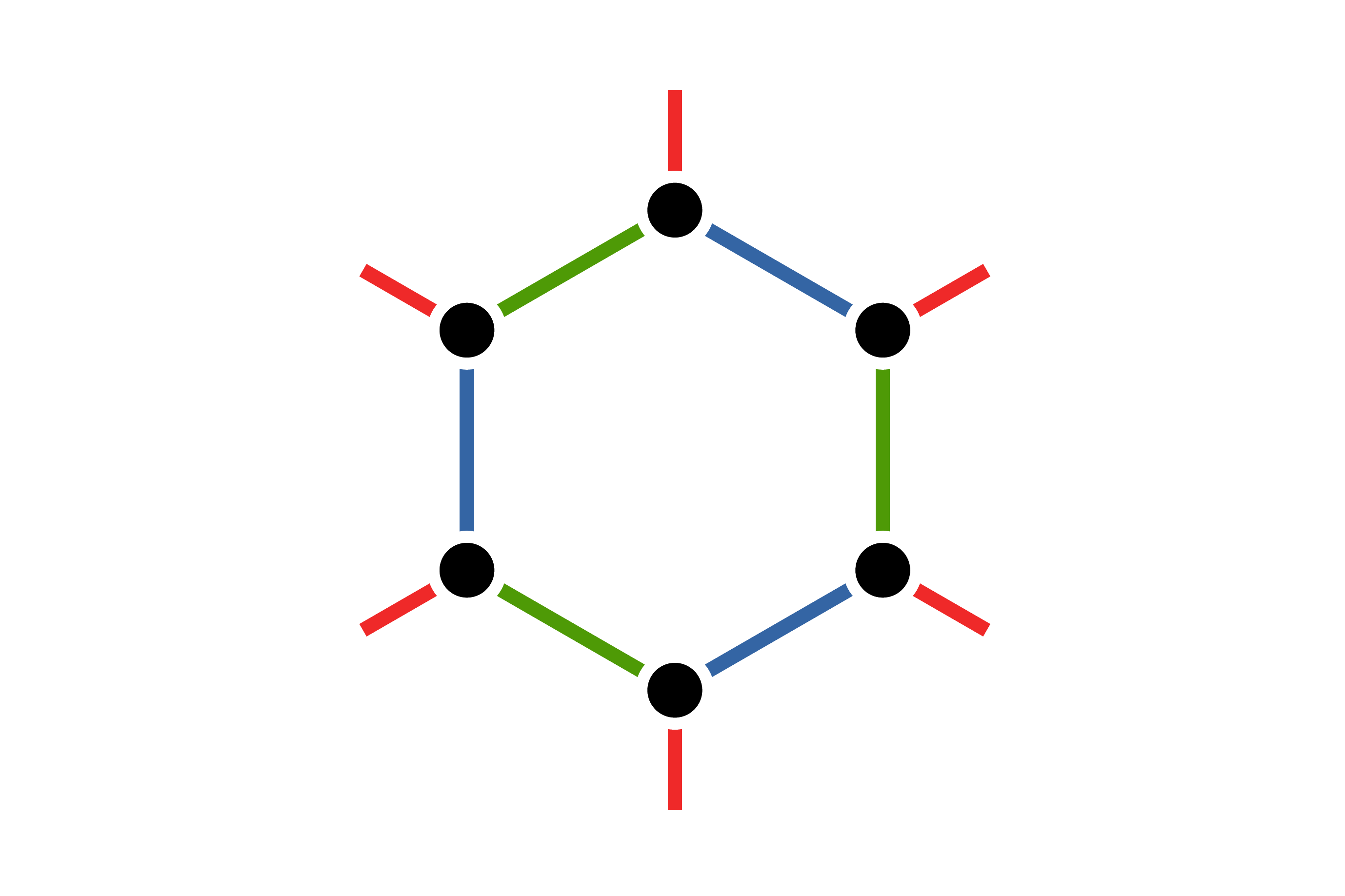}}\put(0.24850542,0.18195712){\color[rgb]{0,0,0}\makebox(0,0)[rt]{\lineheight{1.25}\smash{\begin{tabular}[t]{r}$r_0$\end{tabular}}}}\put(0.24850542,0.47151129){\color[rgb]{0,0,0}\makebox(0,0)[rt]{\lineheight{1.25}\smash{\begin{tabular}[t]{r}$r_0$\end{tabular}}}}\put(0.50358095,0.62727216){\color[rgb]{0,0,0}\makebox(0,0)[t]{\lineheight{1.25}\smash{\begin{tabular}[t]{c}$r_1$\end{tabular}}}}\put(0.75169756,0.47151129){\color[rgb]{0,0,0}\makebox(0,0)[lt]{\lineheight{1.25}\smash{\begin{tabular}[t]{l}$r_1$\end{tabular}}}}\put(0.75169756,0.18195712){\color[rgb]{0,0,0}\makebox(0,0)[lt]{\lineheight{1.25}\smash{\begin{tabular}[t]{l}$r_2$\end{tabular}}}}\put(0.50183757,0.01822305){\color[rgb]{0,0,0}\makebox(0,0)[t]{\lineheight{1.25}\smash{\begin{tabular}[t]{c}$r_2$\end{tabular}}}}\end{picture}\endgroup  		\end{scriptsize}
		\caption{ $A=\left\{ \rho_{0}, \rho_{1},\rho_{2} \right\}$}\label{fig:whyt_012}
	\end{subfigure}
	\caption{Whytoffian operators for rank $3$. The edges in red, green and blue represent $0$-, $1$- and $2$- adjacencies, respectively.} \label{fig:whytoffian}
\end{figure}

	  Figure\nobreakspace \ref {fig:whytoffian} shows voltage operators $(\cY, \eta)$ for five of the Wythoffian constructions to be seen as a voltage operation. (As pointed out above, the two remaining ones correspond to when $A = \{\rho_0\}$ and $A = \{\rho_2\}$ that are the identity and dual operations, respectively.)

\begin{figure}
	\centering
    \begin{scriptsize}
\def\svgwidth{.5\textwidth}
	\begingroup \makeatletter \providecommand\color[2][]{\errmessage{(Inkscape) Color is used for the text in Inkscape, but the package 'color.sty' is not loaded}\renewcommand\color[2][]{}}\providecommand\transparent[1]{\errmessage{(Inkscape) Transparency is used (non-zero) for the text in Inkscape, but the package 'transparent.sty' is not loaded}\renewcommand\transparent[1]{}}\providecommand\rotatebox[2]{#2}\newcommand*\fsize{\dimexpr\f@size pt\relax}\newcommand*\lineheight[1]{\fontsize{\fsize}{#1\fsize}\selectfont}\ifx\svgwidth\undefined \setlength{\unitlength}{966.00767601bp}\ifx\svgscale\undefined \relax \else \setlength{\unitlength}{\unitlength * \real{\svgscale}}\fi \else \setlength{\unitlength}{\svgwidth}\fi \global\let\svgwidth\undefined \global\let\svgscale\undefined \makeatother \begin{picture}(1,0.73907862)\lineheight{1}\setlength\tabcolsep{0pt}\put(0,0){\includegraphics[width=\unitlength,page=1]{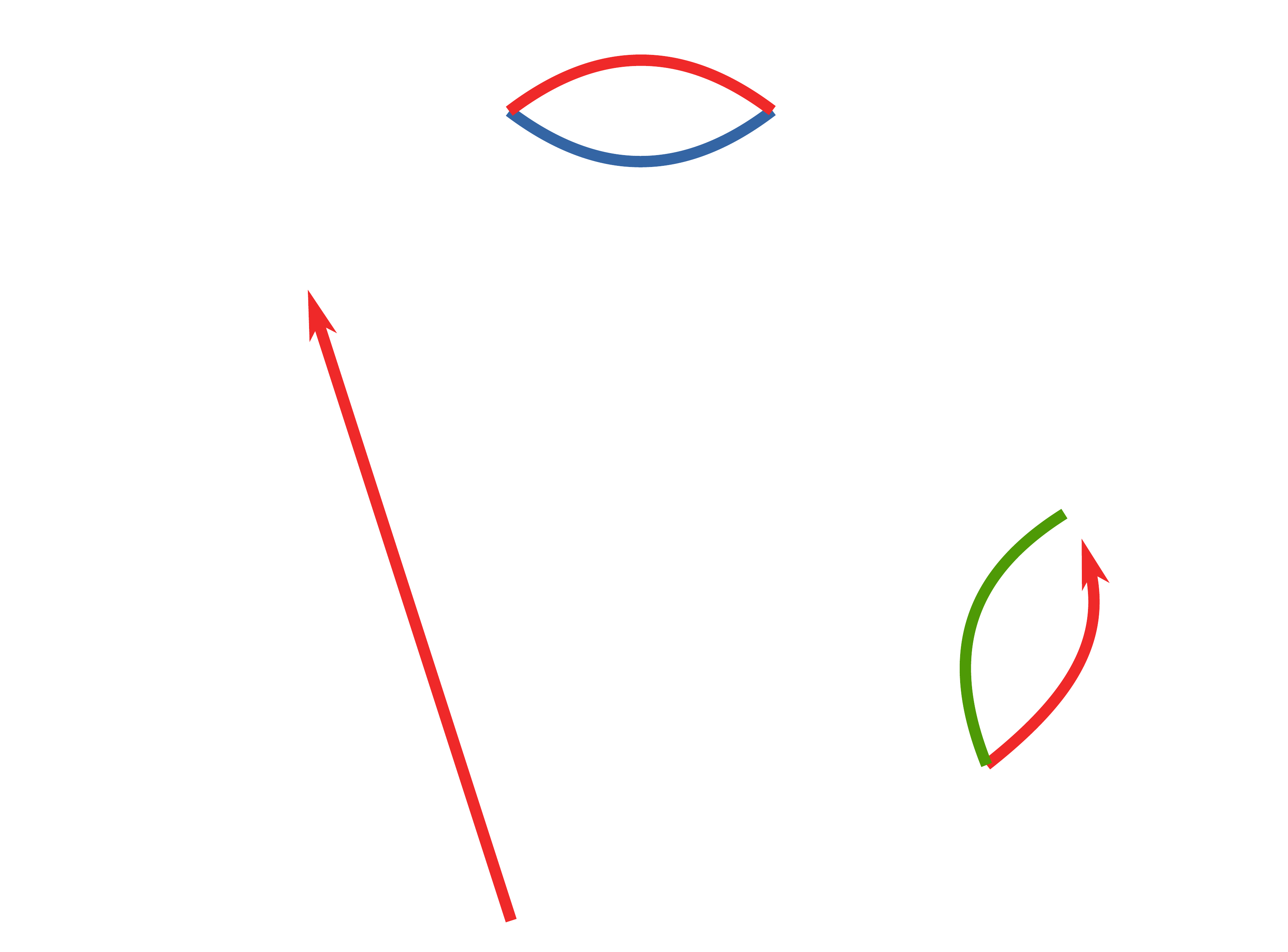}}\put(0.69752058,0.38578242){\color[rgb]{0,0,0}\makebox(0,0)[rt]{\lineheight{1.25}\smash{\begin{tabular}[t]{r}$r_2 r_1$\end{tabular}}}}\put(0,0){\includegraphics[width=\unitlength,page=2]{snub.pdf}}\put(0.31855147,0.38578242){\color[rgb]{0,0,0}\makebox(0,0)[lt]{\lineheight{1.25}\smash{\begin{tabular}[t]{l}$r_0 r_1$\end{tabular}}}}\put(0.86755895,0.21347074){\color[rgb]{0,0,0}\makebox(0,0)[lt]{\lineheight{1.25}\smash{\begin{tabular}[t]{l}$r_2 r_1$\end{tabular}}}}\put(0.13260783,0.21347074){\color[rgb]{0,0,0}\makebox(0,0)[rt]{\lineheight{1.25}\smash{\begin{tabular}[t]{r}$r_0 r_1$\end{tabular}}}}\put(0.49749165,0.71609023){\color[rgb]{0,0,0}\makebox(0,0)[t]{\lineheight{1.25}\smash{\begin{tabular}[t]{c}$r_0 r_2$\end{tabular}}}}\end{picture}\endgroup  	\end{scriptsize}
	\caption{ Snub operator}\label{fig:snub}
\end{figure}

    The \emph{omnitruncation} in rank $3$ (Figure\nobreakspace \ref {fig:whyt_012}) can be easily generalized to higher dimensions for a regular convex polytope $\cP$: it corresponds to the Whytoffian construction where $A$ is the set of all reflections $\{\rho_0, \rho_1, \dots, \rho_{n-1}\}$.
    We can define this operation when $\cP$ is an abstract polytope (regular or not) in the following way.
    For $i\in\{0,1,\ldots,n-1\}$, the $i$-faces of $\ot(\cP)$ are chains of $\cP$ that do not contain the greatest nor least elements and have $n-i$ elements,
    and given two faces $C$ and $C'$ of $\ot(\cP)$ we say that $C\leq_{\ot(\cP)} C'$ if and only if $C'\subset C$.
    In particular, the vertices of $\ot(\cP)$ are the flags of $\cP$ and an edge of the form $\Phi\setminus\{\Phi_i\}$ joins the vertices $\Phi$ and $\Phi^i$.
    This means that the 1-skeleton of $\ot(\cP)$ is the flag graph of $\cP$.
    Another way to construct $\ot(\cP)$ is as the colorful polytope of the flag graph of $\cP$ \cite{AraujoPardoHubardOliverosSchulte_2013_ColorfulPolytopesGraphs}.

    The omnitruncation on rank $n$ can be thought as a voltage operation. To see this, use the following construction for the voltage operator. Let $\cS$ be  the flag graph of an $(n-1)$-simplex. First, change the color of each edge by increasing it by 1 (i.e. if two flags of $\cS$ are $i$-adjacent, they will be joined by an edge of color $i+1$). All these edges will have the identity voltage. Note that now each vertex has an edge of color $i$, for $i\in \{1,2,\dots, n-1\}$.
    Now, label each connected component of the graph with edges of colors $ 1, \dots, n-2$ with a different number from 0 to $n-2$.
    Then, add a semiedge of color $0$ at each vertex with voltage $r_j$ where $j$ is the label of connected component of the vertex.
    The resulting voltage operator $(\cY,\eta)$ satisfies that the flag graph of $\ot(\cP)$ is $\cP\ertimes \cY$.

    As can be expected all Wythoffian constructions for higher dimensions can be seen as voltage operations. We do not give here all the voltage operators, as for rank $n$, there can be up to $n!$ vertices in it. However, Theorem~Theorem\nobreakspace \ref {thm:derivedAsProduct} gives a way to find each of them.

   The {\em snub} of a polyhedron can be thought as a generalized Wythoffian construction, where the rotational subgroup is used. We shall not give details of this kind of constructions, but rather note some interesting differences between the snub and the Wythoffian described above.

   The {\em snub operator} given in Figure\nobreakspace \ref {fig:snub}. In contrast with the Wythoffian operations that preserve connectivity of any maniplex, the snub operator only preserves connectivity for non-orientable premaniplexes.
	    If we apply this snub operation to an orientable premaniplex $\cX$ (in particular to a regular convex polyhedron) we get two copies of what is usually referred to as \emph{the snub $\cX$}, each copy having all the ``rotational'' symmetry of $\cX$ but not the reflection symmetry of $\cX$.
	    We will explore this phenomenon deeper in Section\nobreakspace \ref {sec:automorphisms}.

\subsubsection*{Prisms and pyramids over polytopes}

Let $\cP$ be an abstract polytope. The \emph{pyramid over $\cP$}, denoted by $\pyr(\cP)$ is the poset $\cP\times \{0,1\}$ with the product order, that is, $(F,\ell)\leq(G,\ell')$ if and only if $F\leq G$ and $\ell\leq \ell'$.
Similarly, the \emph{prism over $\cP$}, denoted by $\prism(\cP)$ is defined as the poset $((\cP\setminus\cP_{-1})\times \Lambda) \cup\{F_{-1}\}$ where $\Lambda$ is the poset $\{\{\lambda_0\},\{\lambda_1\},\{\lambda_0,\lambda_1\}\}$ ordered by inclusion and $F_{-1}$ is the unique minimum element of the prism.

Both the prism and the pyramid over an $n$-polytope are $(n+1)$-polytopes.
They are particular cases of products of polytopes, in the sense of \cite{GleasonHubard_2018_ProductsAbstractPolytopes}: the pyramid is the join product by a vertex, while the prism is the direct product by an edge.

One can see that each flag $\Psi$ in $\pyr(\cP)$ is of the form
\[
    \Psi = \{(\Phi_{-1},0),(\Phi_0,0),\ldots,(\Phi_t,0),(\Phi_t,1),(\Phi_{t+1},1),\ldots,(\Phi_n,1)\},
\]
where $\Phi$ is a flag in $\cP$ and $t\in\{-1,0,1,\ldots,n\}$. Hence, we identify the set of flags of $\pyr(\cP)$ with $\cF(\cP)\times \{-1,0,\ldots,n\}$.

In a similar way, one can identify the set of flags of $\prism(\cP)$ with $\cF\times \{0,1,\ldots,n\}\times\{\lambda_0,\lambda_1\}$, where the last coordinate tells us if the 0-face is of the form $(\Phi_0,\{\lambda_0\})$ or $(\Phi_0,\{\lambda_1\})$.

Let us first observe the pyramid. The $i$-adjacencies of the flags (see \cite{GleasonHubard_2018_ProductsAbstractPolytopes}), for $i \in \{0,1,\dots, n\}$, are given by
\[
    (\Phi,t)^i=\begin{cases}
      (\Phi^i,t) & \text{if } 0\leq i<t,\\
      (\Phi,t-1) & \text{if } i=t,\\
      (\Phi,t+1) & \text{if } i=t+1,\\
      (\Phi^{i-1},t) & \text{if } t+1<i.
    \end{cases}
\]

\begin{figure}
	\centering
    \begin{scriptsize}
\def\svgwidth{\textwidth}
	\begingroup \makeatletter \providecommand\color[2][]{\errmessage{(Inkscape) Color is used for the text in Inkscape, but the package 'color.sty' is not loaded}\renewcommand\color[2][]{}}\providecommand\transparent[1]{\errmessage{(Inkscape) Transparency is used (non-zero) for the text in Inkscape, but the package 'transparent.sty' is not loaded}\renewcommand\transparent[1]{}}\providecommand\rotatebox[2]{#2}\newcommand*\fsize{\dimexpr\f@size pt\relax}\newcommand*\lineheight[1]{\fontsize{\fsize}{#1\fsize}\selectfont}\ifx\svgwidth\undefined \setlength{\unitlength}{1357.23324273bp}\ifx\svgscale\undefined \relax \else \setlength{\unitlength}{\unitlength * \real{\svgscale}}\fi \else \setlength{\unitlength}{\svgwidth}\fi \global\let\svgwidth\undefined \global\let\svgscale\undefined \makeatother \begin{picture}(1,0.23111161)\lineheight{1}\setlength\tabcolsep{0pt}\put(0,0){\includegraphics[width=\unitlength,page=1]{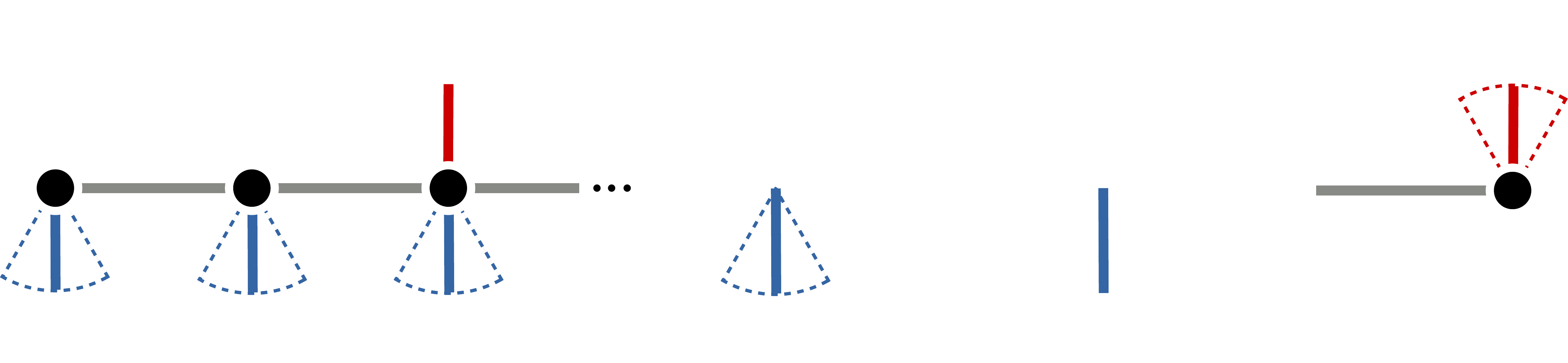}}\put(0.09799002,0.12170555){\color[rgb]{0,0,0}\makebox(0,0)[t]{\lineheight{1.25}\smash{\begin{tabular}[t]{c}$0$\end{tabular}}}}\put(0.22330285,0.12170555){\color[rgb]{0,0,0}\makebox(0,0)[t]{\lineheight{1.25}\smash{\begin{tabular}[t]{c}$1$\end{tabular}}}}\put(0.0353336,0.02772093){\color[rgb]{0,0,0}\makebox(0,0)[t]{\lineheight{1.25}\smash{\begin{tabular}[t]{c}$r_{i-1}$\end{tabular}}}}\put(0.16067806,0.02727338){\color[rgb]{0,0,0}\makebox(0,0)[t]{\lineheight{1.25}\smash{\begin{tabular}[t]{c}$r_{i-1}$\end{tabular}}}}\put(0.28652188,0.19466808){\color[rgb]{0,0,0}\makebox(0,0)[t]{\lineheight{1.25}\smash{\begin{tabular}[t]{c}$r_i$\end{tabular}}}}\put(0,0){\includegraphics[width=\unitlength,page=2]{pyr.pdf}}\put(0.83959674,0.19391411){\color[rgb]{0,0,0}\makebox(0,0)[t]{\lineheight{1.25}\smash{\begin{tabular}[t]{c}$r_i$\end{tabular}}}}\put(0,0){\includegraphics[width=\unitlength,page=3]{pyr.pdf}}\put(0.90208062,0.12170555){\color[rgb]{0,0,0}\makebox(0,0)[t]{\lineheight{1.25}\smash{\begin{tabular}[t]{c}$n$\end{tabular}}}}\put(0,0){\includegraphics[width=\unitlength,page=4]{pyr.pdf}}\put(0.49481395,0.1948047){\color[rgb]{0,0,0}\makebox(0,0)[t]{\lineheight{1.25}\smash{\begin{tabular}[t]{c}$r_i$\end{tabular}}}}\put(0.46348574,0.12170555){\color[rgb]{0,0,0}\makebox(0,0)[rt]{\lineheight{1.25}\smash{\begin{tabular}[t]{r}$t$\end{tabular}}}}\put(0.52614216,0.12170555){\color[rgb]{0,0,0}\makebox(0,0)[lt]{\lineheight{1.25}\smash{\begin{tabular}[t]{l}$t+1$\end{tabular}}}}\put(0.70384118,0.19391411){\color[rgb]{0,0,0}\makebox(0,0)[t]{\lineheight{1.25}\smash{\begin{tabular}[t]{c}$r_i$\end{tabular}}}}\put(0.03543615,0.00703192){\color[rgb]{0,0,0}\makebox(0,0)[t]{\lineheight{1.25}\smash{\begin{tabular}[t]{c}$i > 0 $\end{tabular}}}}\put(0.16086687,0.00796255){\color[rgb]{0,0,0}\makebox(0,0)[t]{\lineheight{1.25}\smash{\begin{tabular}[t]{c}$i > 1 $\end{tabular}}}}\put(0.28607651,0.21631953){\color[rgb]{0,0,0}\makebox(0,0)[t]{\lineheight{1.25}\smash{\begin{tabular}[t]{c}$i < 1$\end{tabular}}}}\put(0.49436862,0.21645615){\color[rgb]{0,0,0}\makebox(0,0)[t]{\lineheight{1.25}\smash{\begin{tabular}[t]{c}$i < t $\end{tabular}}}}\put(0.7036687,0.21569017){\color[rgb]{0,0,0}\makebox(0,0)[t]{\lineheight{1.25}\smash{\begin{tabular}[t]{c}$i < n-2$\end{tabular}}}}\put(0.28617968,0.00796249){\color[rgb]{0,0,0}\makebox(0,0)[t]{\lineheight{1.25}\smash{\begin{tabular}[t]{c}$i > 2 $\end{tabular}}}}\put(0.49474174,0.00727451){\color[rgb]{0,0,0}\makebox(0,0)[t]{\lineheight{1.25}\smash{\begin{tabular}[t]{c}$i > t+1 $\end{tabular}}}}\put(0.70385749,0.00727451){\color[rgb]{0,0,0}\makebox(0,0)[t]{\lineheight{1.25}\smash{\begin{tabular}[t]{c}$i > n-1 $\end{tabular}}}}\put(0.83942433,0.21569022){\color[rgb]{0,0,0}\makebox(0,0)[t]{\lineheight{1.25}\smash{\begin{tabular}[t]{c}$i < n-1$\end{tabular}}}}\put(0.96473716,0.21569022){\color[rgb]{0,0,0}\makebox(0,0)[t]{\lineheight{1.25}\smash{\begin{tabular}[t]{c}$i < n$\end{tabular}}}}\put(0.28599088,0.02727331){\color[rgb]{0,0,0}\makebox(0,0)[t]{\lineheight{1.25}\smash{\begin{tabular}[t]{c}$r_{i-1}$\end{tabular}}}}\put(0.33790244,0.12170555){\color[rgb]{0,0,0}\makebox(0,0)[t]{\lineheight{1.25}\smash{\begin{tabular}[t]{c}$2$\end{tabular}}}}\put(0.49455292,0.02658535){\color[rgb]{0,0,0}\makebox(0,0)[t]{\lineheight{1.25}\smash{\begin{tabular}[t]{c}$r_{i-1}$\end{tabular}}}}\put(0.67234049,0.12170555){\color[rgb]{0,0,0}\makebox(0,0)[rt]{\lineheight{1.25}\smash{\begin{tabular}[t]{r}$n-2$\end{tabular}}}}\put(0.76632512,0.12170555){\color[rgb]{0,0,0}\makebox(0,0)[t]{\lineheight{1.25}\smash{\begin{tabular}[t]{c}$n-1$\end{tabular}}}}\put(0.7036687,0.02658535){\color[rgb]{0,0,0}\makebox(0,0)[t]{\lineheight{1.25}\smash{\begin{tabular}[t]{c}$r_{i-1}$\end{tabular}}}}\put(0.96473716,0.19391411){\color[rgb]{0,0,0}\makebox(0,0)[t]{\lineheight{1.25}\smash{\begin{tabular}[t]{c}$r_i$\end{tabular}}}}\end{picture}\endgroup  	\end{scriptsize}
	\caption{Voltage operator for the pyramid over an $n$-polytope}\label{fig:pyr}
\end{figure}

Let $(\cY,\eta)$ be the $(n,n+1)$-voltage operator whose vertices are the numbers $\{-1,0,\ldots,n\}$ with the $i$-adjacencies (for $i \in \{0,1,\dots, n\}$) given by:
\[
    t^i=\begin{cases}
      t & \text{if }  i<t \text{ or } t+1< i, \\
      t-1 & \text{if }  i=t,\\
      t+1 & \text{if } i=t+1,
    \end{cases}
\]
and the voltage assignment given by:
\[
    \xi(\pth{t})=\begin{cases}
      r_i & \text{if }  i<t,\\
      1 & \text{if }  i=t \text{ or } i=t+1,\\
      r_{i-1} & \text{if } t+1<i.
    \end{cases}
\]
(see Figure\nobreakspace \ref {fig:pyr}). Then the flag graph of $\pyr(\cP)$ is $\cG(\cP)\ertimes \cY$ where $(\cY,\eta)$.

We now turn our attention to the prism, where the $i$-adjacencies of the flags (see \cite{GleasonHubard_2018_ProductsAbstractPolytopes}), for $i \in \{0,1,\dots, n\}$, are given by
\[
    (\Phi,t,\lambda)^i=\begin{cases}
      (\Phi^i,t,\lambda) & \text{if } i<t,\\
       (\Phi,t,\lambda') & \text{if } 0=i=t, \\
      (\Phi,t-1,\lambda) & \text{if } 0<i=t,\\
      (\Phi,t+1,\lambda) & \text{if } i=t+1,\\
      (\Phi^{i-1},t,\lambda) & \text{if } t+1<i,\\
    \end{cases}
\]
where $\lambda'$ is $\lambda_1$ if $\lambda=\lambda_0$ and vice versa.

\begin{figure}
	\centering
    \begin{scriptsize}
\def\svgwidth{\textwidth}
	\begingroup \makeatletter \providecommand\color[2][]{\errmessage{(Inkscape) Color is used for the text in Inkscape, but the package 'color.sty' is not loaded}\renewcommand\color[2][]{}}\providecommand\transparent[1]{\errmessage{(Inkscape) Transparency is used (non-zero) for the text in Inkscape, but the package 'transparent.sty' is not loaded}\renewcommand\transparent[1]{}}\providecommand\rotatebox[2]{#2}\newcommand*\fsize{\dimexpr\f@size pt\relax}\newcommand*\lineheight[1]{\fontsize{\fsize}{#1\fsize}\selectfont}\ifx\svgwidth\undefined \setlength{\unitlength}{1021.91286102bp}\ifx\svgscale\undefined \relax \else \setlength{\unitlength}{\unitlength * \real{\svgscale}}\fi \else \setlength{\unitlength}{\svgwidth}\fi \global\let\svgwidth\undefined \global\let\svgscale\undefined \makeatother \begin{picture}(1,0.56703244)\lineheight{1}\setlength\tabcolsep{0pt}\put(0,0){\includegraphics[width=\unitlength,page=1]{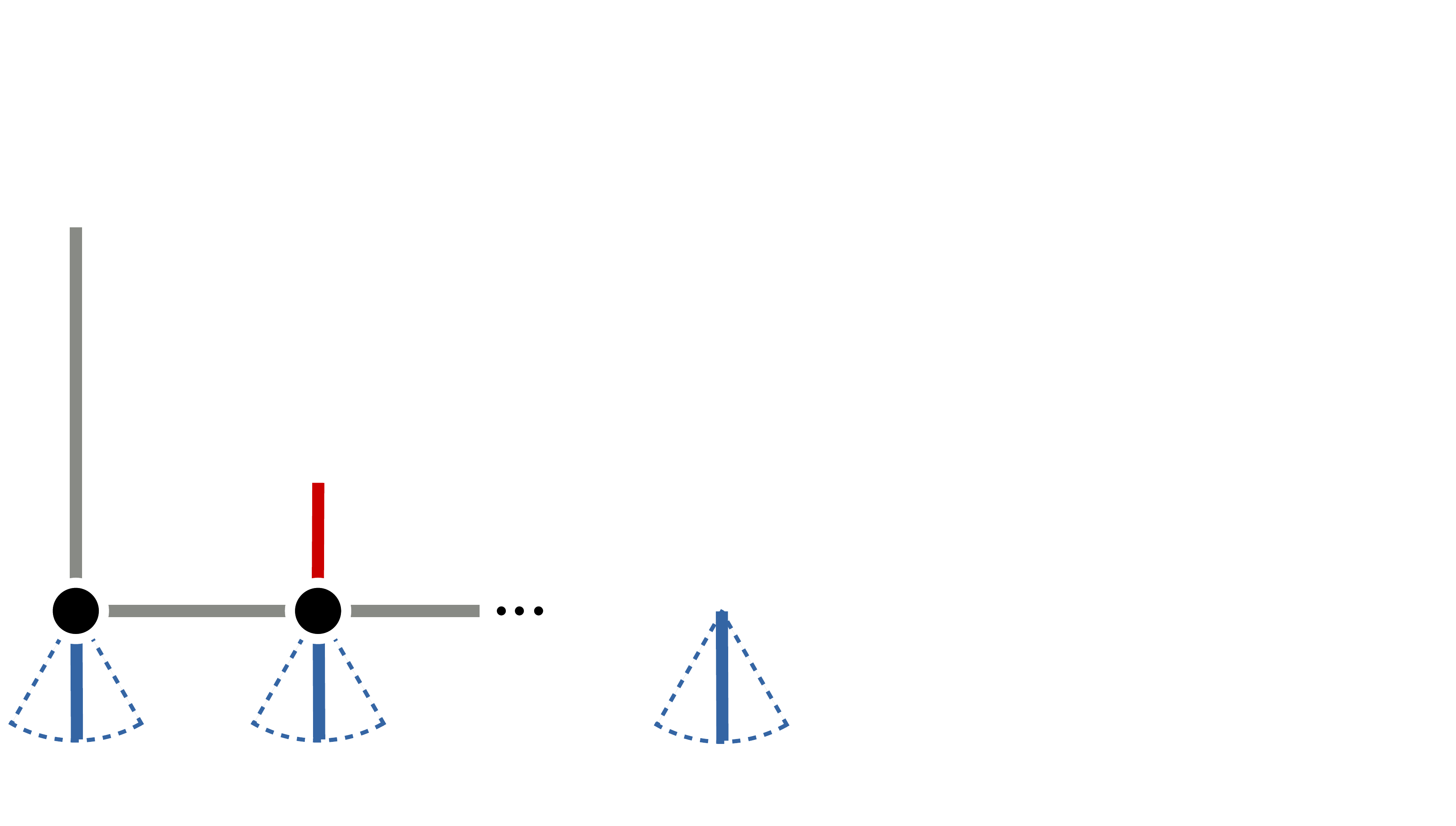}}\put(0.13528679,0.16131864){\color[rgb]{0,0,0}\makebox(0,0)[t]{\lineheight{1.25}\smash{\begin{tabular}[t]{c}$1$\end{tabular}}}}\put(0.05211292,0.03590043){\color[rgb]{0,0,0}\makebox(0,0)[t]{\lineheight{1.25}\smash{\begin{tabular}[t]{c}$r_{i-1}$\end{tabular}}}}\put(0,0){\includegraphics[width=\unitlength,page=2]{pri.pdf}}\put(0.45428092,0.16131864){\color[rgb]{0,0,0}\makebox(0,0)[rt]{\lineheight{1.25}\smash{\begin{tabular}[t]{r}$t$\end{tabular}}}}\put(0.53749679,0.16131864){\color[rgb]{0,0,0}\makebox(0,0)[lt]{\lineheight{1.25}\smash{\begin{tabular}[t]{l}$t+1$\end{tabular}}}}\put(0.05236368,0.01025314){\color[rgb]{0,0,0}\makebox(0,0)[t]{\lineheight{1.25}\smash{\begin{tabular}[t]{c}$i > 1 $\end{tabular}}}}\put(0.21879541,0.01025306){\color[rgb]{0,0,0}\makebox(0,0)[t]{\lineheight{1.25}\smash{\begin{tabular}[t]{c}$i > 2 $\end{tabular}}}}\put(0.49579295,0.00933933){\color[rgb]{0,0,0}\makebox(0,0)[t]{\lineheight{1.25}\smash{\begin{tabular}[t]{c}$i > t+1 $\end{tabular}}}}\put(0.05207107,0.52192079){\color[rgb]{0,0,0}\makebox(0,0)[t]{\lineheight{1.25}\smash{\begin{tabular}[t]{c}$i > 1 $\end{tabular}}}}\put(0.21850281,0.52192071){\color[rgb]{0,0,0}\makebox(0,0)[t]{\lineheight{1.25}\smash{\begin{tabular}[t]{c}$i > 2 $\end{tabular}}}}\put(0.49550035,0.52100698){\color[rgb]{0,0,0}\makebox(0,0)[t]{\lineheight{1.25}\smash{\begin{tabular}[t]{c}$i > t+1 $\end{tabular}}}}\put(0.21854465,0.03590034){\color[rgb]{0,0,0}\makebox(0,0)[t]{\lineheight{1.25}\smash{\begin{tabular}[t]{c}$r_{i-1}$\end{tabular}}}}\put(0.28748996,0.16131864){\color[rgb]{0,0,0}\makebox(0,0)[t]{\lineheight{1.25}\smash{\begin{tabular}[t]{c}$2$\end{tabular}}}}\put(0.49554217,0.03498663){\color[rgb]{0,0,0}\makebox(0,0)[t]{\lineheight{1.25}\smash{\begin{tabular}[t]{c}$r_{i-1}$\end{tabular}}}}\put(0.7316672,0.16131864){\color[rgb]{0,0,0}\makebox(0,0)[rt]{\lineheight{1.25}\smash{\begin{tabular}[t]{r}$n-1$\end{tabular}}}}\put(0.85649101,0.16131864){\color[rgb]{0,0,0}\makebox(0,0)[t]{\lineheight{1.25}\smash{\begin{tabular}[t]{c}$n$\end{tabular}}}}\put(0.21909418,0.25738655){\color[rgb]{0,0,0}\makebox(0,0)[t]{\lineheight{1.25}\smash{\begin{tabular}[t]{c}$r_i$\end{tabular}}}}\put(0.95364967,0.25638519){\color[rgb]{0,0,0}\makebox(0,0)[t]{\lineheight{1.25}\smash{\begin{tabular}[t]{c}$r_i$\end{tabular}}}}\put(0.49573327,0.257568){\color[rgb]{0,0,0}\makebox(0,0)[t]{\lineheight{1.25}\smash{\begin{tabular}[t]{c}$r_i$\end{tabular}}}}\put(0.77334862,0.25638519){\color[rgb]{0,0,0}\makebox(0,0)[t]{\lineheight{1.25}\smash{\begin{tabular}[t]{c}$r_i$\end{tabular}}}}\put(0.21850281,0.2861425){\color[rgb]{0,0,0}\makebox(0,0)[t]{\lineheight{1.25}\smash{\begin{tabular}[t]{c}$i < 1$\end{tabular}}}}\put(0.49514172,0.28632395){\color[rgb]{0,0,0}\makebox(0,0)[t]{\lineheight{1.25}\smash{\begin{tabular}[t]{c}$i < t $\end{tabular}}}}\put(0.77311955,0.28530663){\color[rgb]{0,0,0}\makebox(0,0)[t]{\lineheight{1.25}\smash{\begin{tabular}[t]{c}$i < n-1$\end{tabular}}}}\put(0.95342061,0.28530669){\color[rgb]{0,0,0}\makebox(0,0)[t]{\lineheight{1.25}\smash{\begin{tabular}[t]{c}$i < n$\end{tabular}}}}\put(0,0){\includegraphics[width=\unitlength,page=3]{pri.pdf}}\put(0.13526558,0.42804769){\color[rgb]{0,0,0}\makebox(0,0)[t]{\lineheight{1.25}\smash{\begin{tabular}[t]{c}$1$\end{tabular}}}}\put(0,0){\includegraphics[width=\unitlength,page=4]{pri.pdf}}\put(0.45425972,0.42804769){\color[rgb]{0,0,0}\makebox(0,0)[rt]{\lineheight{1.25}\smash{\begin{tabular}[t]{r}$t$\end{tabular}}}}\put(0.53747559,0.42804769){\color[rgb]{0,0,0}\makebox(0,0)[lt]{\lineheight{1.25}\smash{\begin{tabular}[t]{l}$t+1$\end{tabular}}}}\put(0.28746875,0.42804769){\color[rgb]{0,0,0}\makebox(0,0)[t]{\lineheight{1.25}\smash{\begin{tabular}[t]{c}$2$\end{tabular}}}}\put(0.73164596,0.42804769){\color[rgb]{0,0,0}\makebox(0,0)[rt]{\lineheight{1.25}\smash{\begin{tabular}[t]{r}$n-1$\end{tabular}}}}\put(0.85646977,0.42804769){\color[rgb]{0,0,0}\makebox(0,0)[t]{\lineheight{1.25}\smash{\begin{tabular}[t]{c}$n$\end{tabular}}}}\put(0.02433245,0.27227319){\color[rgb]{0,0,0}\makebox(0,0)[rt]{\lineheight{1.25}\smash{\begin{tabular}[t]{r}$0$\end{tabular}}}}\put(0.05182031,0.54756808){\color[rgb]{0,0,0}\makebox(0,0)[t]{\lineheight{1.25}\smash{\begin{tabular}[t]{c}$r_{i-1}$\end{tabular}}}}\put(0.21825205,0.54756799){\color[rgb]{0,0,0}\makebox(0,0)[t]{\lineheight{1.25}\smash{\begin{tabular}[t]{c}$r_{i-1}$\end{tabular}}}}\put(0.49524957,0.54665429){\color[rgb]{0,0,0}\makebox(0,0)[t]{\lineheight{1.25}\smash{\begin{tabular}[t]{c}$r_{i-1}$\end{tabular}}}}\end{picture}\endgroup  	\end{scriptsize}
	\caption{Voltage operator for the prism over an $n$-polytope}\label{fig:prism}
\end{figure}

Let $(\cY',\eta')$ be the $(n,n+1)$-voltage operator whose vertices are $\{0,\ldots,n\}\times\{\lambda_0,\lambda_1\}$ with the adjacencies given by:
\[
    (t,\lambda)^i=\begin{cases}
      (t,\lambda) & \text{if } i<t \text{ or } i>t+1,\\
      (t,\lambda') & \text{if }  0=i=t, \\
      (t-1,\lambda) & \text{if } i=t>0,\\
      (t+1,\lambda) & \text{if } i=t+1,\\
    \end{cases}
\]
and the voltage assignment given by:
\[
    \eta'(\pth{t})=\begin{cases}
      r_i & \text{if } i<t,\\
      1 & \text{if }  i=t \text{ or } i=t+1,\\
      r_{i-1} & \text{if } t+1<i.
    \end{cases}
\]
We show this operator in Figure\nobreakspace \ref {fig:prism}. Then the flag graph of $\pyr(\cP)$ is $\cG(\cP)\ertimes \cY$ where $(\cY,\eta)$.

The flag graph of $\prism(\cP)$ is $\cG(\cP)\ertimes[\eta']\cY'$.

\subsubsection*{The trapezotope}
Given a polytope $\cP$, we define the \emph{trapezotope over $\cP$} denoted by $\trp(\cP)$ as follows:
the faces of rank $i$ of the trapezotope are the ordered pairs $(F,G)$ where $F$ and $G$ are faces of $\cP$, $F\leq G$ and $\rk(G)-\rk(F)=i$.
In particular, the vertices of $\trp(\cP)$ are the pairs $(F,F)$, so in correspondence to the faces of $\cP$, and the edges of $\trp(\cP)$ are in correspondence to the edges of the Hasse diagram of $\cP$. Thus, the $1$-skeleton of the trapezotope of $\cP$ is the Hasse diagram of $\cP$.
The trapezotope over $\cP$ is in fact the dual of the \emph{antiprism over $\cP$}, defined in \cite{GleasonHubard_2021_AntiprismAbstractPolytope}.

In \cite{GleasonHubard_2021_AntiprismAbstractPolytope} the authors showed that the trapezotope of a polytope $\cP$ (or rather its dual) is in fact a polytope and that the flags are in one-to-one correspondence with the set $\cF(\cP)\times \bZ_2^{n+1}$.
When we denote a flag of $\trp(\cP)$ by $(\Phi,v)$, the $i$-th coordinate of the vector $v$ tells us if the faces $(\Phi,v)_i=(F_i,G_i)$ and $(\Phi,v)_{i-1}=(F_{i-1},G_{i-1})$ differ in the first coordinate (when $v_i=0$) or in the second one (when $v_i=1$).
Using this natural correspondence, one can prove that the $i$-adjacent flag of the flag $(\Phi,v)$ is:
\[
    (\Phi,v)^i=
    \begin{cases}
       (\Phi,v+e_0) &   \text{if } i=0,\\
       (\Phi^j,v) & \text{if } v_{i+1}=v_i,\\
        (\Phi,v+e_i+e_{i+1}) & \text{if } v_{i+1} \neq v_i \text{ and } i\neq 0,
    \end{cases}
\]
where $j$ is the rank of $F_i$ if $v_i=0$ or the rank of $G_i$ if $v_i=1$ and $e_{i}$ denotes the vector in $\bZ_2^{n+1}$ whose all but the $i$-th coordinate are zero.
Let $\cY$ be the premaniplex whose vertices are $\bZ_2^n$ and the adjacencies are given by:
\[
    v^i=
    \begin{cases}
       v+e_0 &   \text{if } i=0,\\
       v & \text{if } v_{i+1}=v_i\\
        v+e_i+e_{i+1} & \text{if } v_{i+1} \neq v_i \text{   and } i\neq 0.
    \end{cases}
\]
And let $\eta$ be the voltage assignment given by
\[
    \eta(\pth{v})=
    \begin{cases}
       1 &   \text{if } i=0 \text{ or } v_{i+1} \neq v_i,\\
       r_j & \text{if } v_{i+1} = v_i,
    \end{cases}
\]
where $j$ is defined the same way as before.
One can confirm that the voltage operator $(\cY,\eta)$ satisfies that for every polytope $\cP$, the flag graph of $\trp(\cP)$ is $\cG(\cP)\ertimes \cY$.

\begin{figure}
\def\svgwidth{.6\textwidth}
\begin{scriptsize}
\begingroup \makeatletter \providecommand\color[2][]{\errmessage{(Inkscape) Color is used for the text in Inkscape, but the package 'color.sty' is not loaded}\renewcommand\color[2][]{}}\providecommand\transparent[1]{\errmessage{(Inkscape) Transparency is used (non-zero) for the text in Inkscape, but the package 'transparent.sty' is not loaded}\renewcommand\transparent[1]{}}\providecommand\rotatebox[2]{#2}\newcommand*\fsize{\dimexpr\f@size pt\relax}\newcommand*\lineheight[1]{\fontsize{\fsize}{#1\fsize}\selectfont}\ifx\svgwidth\undefined \setlength{\unitlength}{1173.12323142bp}\ifx\svgscale\undefined \relax \else \setlength{\unitlength}{\unitlength * \real{\svgscale}}\fi \else \setlength{\unitlength}{\svgwidth}\fi \global\let\svgwidth\undefined \global\let\svgscale\undefined \makeatother \begin{picture}(1,0.86030148)\lineheight{1}\setlength\tabcolsep{0pt}\put(0,0){\includegraphics[width=\unitlength,page=1]{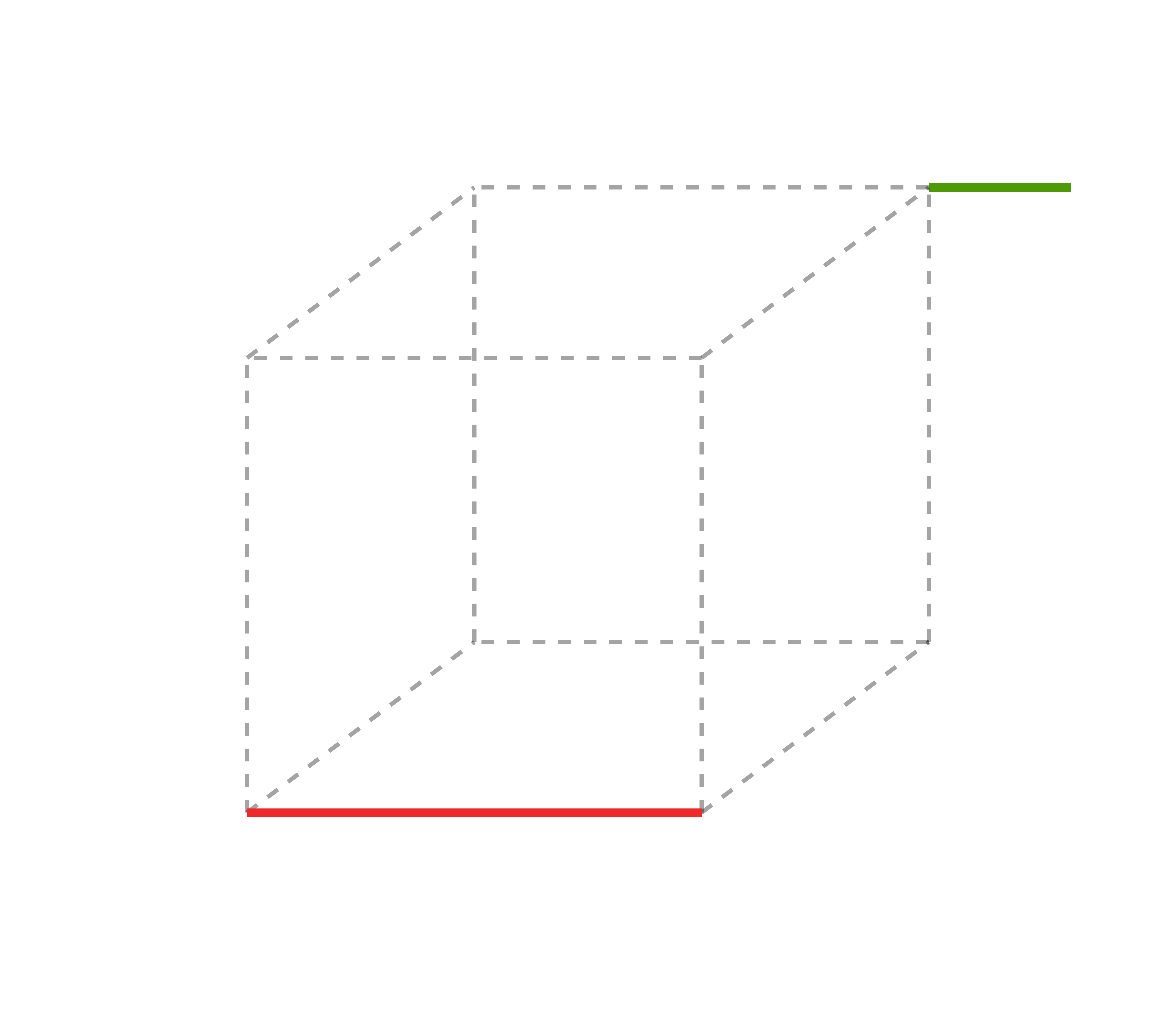}}\put(0.93500697,0.7009833){\color[rgb]{0,0,0}\makebox(0,0)[lt]{\lineheight{1.25}\smash{\begin{tabular}[t]{l}$r_0$\end{tabular}}}}\put(0.79002753,0.8338811){\color[rgb]{0,0,0}\makebox(0,0)[t]{\lineheight{1.25}\smash{\begin{tabular}[t]{c}$r_1$\end{tabular}}}}\put(0.25843632,0.7009833){\color[rgb]{0,0,0}\makebox(0,0)[rt]{\lineheight{1.25}\smash{\begin{tabular}[t]{r}$r_1$\end{tabular}}}}\put(0.74170106,0.16939205){\color[rgb]{0,0,0}\makebox(0,0)[lt]{\lineheight{1.25}\smash{\begin{tabular}[t]{l}$r_0$\end{tabular}}}}\put(0.21010984,0.01233103){\color[rgb]{0,0,0}\makebox(0,0)[t]{\lineheight{1.25}\smash{\begin{tabular}[t]{c}$r_0$\end{tabular}}}}\put(0.93500697,0.31437148){\color[rgb]{0,0,0}\makebox(0,0)[lt]{\lineheight{1.25}\smash{\begin{tabular}[t]{l}$r_1$\end{tabular}}}}\put(0.06513042,0.16939208){\color[rgb]{0,0,0}\makebox(0,0)[rt]{\lineheight{1.25}\smash{\begin{tabular}[t]{r}$r_1$\end{tabular}}}}\put(0.06513042,0.55600389){\color[rgb]{0,0,0}\makebox(0,0)[rt]{\lineheight{1.25}\smash{\begin{tabular}[t]{r}$r_0$\end{tabular}}}}\put(0,0){\includegraphics[width=\unitlength,page=2]{tpr.pdf}}\put(0.18594661,0.12106561){\color[rgb]{0,0,0}\makebox(0,0)[rt]{\lineheight{1.25}\smash{\begin{tabular}[t]{r}$(0,0,0)$\end{tabular}}}}\put(0.60880328,0.12106561){\color[rgb]{0,0,0}\makebox(0,0)[lt]{\lineheight{1.25}\smash{\begin{tabular}[t]{l}$(1,0,0)$\end{tabular}}}}\put(0.25843632,0.30228986){\color[rgb]{0,0,0}\makebox(0,0)[lt]{\lineheight{1.25}\smash{\begin{tabular}[t]{l}$(0,1,0)$\end{tabular}}}}\put(0.18594661,0.58016711){\color[rgb]{0,0,0}\makebox(0,0)[rt]{\lineheight{1.25}\smash{\begin{tabular}[t]{r}$(0,0,1)$\end{tabular}}}}\put(0,0){\includegraphics[width=\unitlength,page=3]{tpr.pdf}}\end{picture}\endgroup  \end{scriptsize}
\caption{The trapezotope operator in rank $3$. Edges in red, green and blue represent $0$-, $1$- and $2$- adjacencies, respectively.}
\label{fig:tpr}
\end{figure}
\subsubsection*{The $k$-bubble}

The last example that we shall see in this section is the
 $k$-bubble  of a polytope. This operation was introduced by Helfand in \cite{Helfand_2013_ConstructionsKOrbit_PhDThesis} as a generalization of the truncation of the vertices.

Let $\cP$ be an abstract $n$-polytope. The $k$-bubble $[\cP]_k$ of $\cP$ is defined as follows.
Let $\cP_i$ denote the set of $i$-faces of $\cP$.
The set $([\cP]_k)_i$ of $i$-faces of $[\cP]_{k}$ is defined by
\[([\cP]_k)_i=
\begin{cases}
  \cP_i & \text{ if } i<k, \\
  \{(F,G):F\in \cP_k,G\in \cP_{k+1}, F<G\} & \text{ if } i=k,\\
  \cP_i \cup \{(F,G):F\in \cP_k,G\in \cP_{i+1}, F<G\} & \text{ if } i>k.
\end{cases}
\]

The order in $[\cP]_k$ is defined by the following rules: Let $H$ and $H'$ be faces in $\cP$ with ranks different than $k$, let $F$ be a $k$-face of $\cP$ and let $G$ and $G'$ be faces of $\cP$ properly containing $F$
\begin{itemize}
    \item $H<H'$ in $[\cP]_k$ if and only if $H<H'$ in $\cP$.
    \item $H<(F,G)$ in $[\cP]_k$ if and only if $H<F$ in $\cP$.
    \item $(F,G)<H$ in $[\cP]_k$ if and only if $G\leq H$ in $\cP$.
    \item $(F,G)<(F,G')$ in $[\cP]_k$ if and only if $G<G'$ in $\cP$.
\end{itemize}

One can see that the flags of $[\cP]_k$ are of the form \[\{F_0,\ldots,F_{k-1},(F_k,F_{k+1}),\ldots,(F_k,F_\ell),F_\ell,\ldots,F_{n-1}\}\] where $\{F_0,\ldots,F_{n-1}\}$ is a flag in $\cP$.
So every flag of $[\cP]_k$ is determined by a flag $\Phi$ of $\cP$ and a number $\ell\in\{k+1,\ldots,n\}$. If we denote this flag by $(\Phi,\ell)$ then the flag adjacencies on $[\cP]_k$ are described by the following rules:

\[
(\Phi,\ell)^i=
\begin{cases}
  (r_i\Phi,\ell) & \text{ if } i\leq k+1\\
  (r_{i+1}\Phi,\ell) & \text{ if } k+1\leq i \leq \ell-2, \\
  (\Phi,\ell-1) & \text{ if } k+1\leq i =\ell-1,\\
  (\Phi,\ell+1) & \text{ if } i=\ell,\\
  (r_i\Phi,\ell) & \text{ if } i \geq \ell+1.
\end{cases}
\]

\begin{figure}
\def\svgwidth{.9\textwidth}
\begin{scriptsize}
\begingroup \makeatletter \providecommand\color[2][]{\errmessage{(Inkscape) Color is used for the text in Inkscape, but the package 'color.sty' is not loaded}\renewcommand\color[2][]{}}\providecommand\transparent[1]{\errmessage{(Inkscape) Transparency is used (non-zero) for the text in Inkscape, but the package 'transparent.sty' is not loaded}\renewcommand\transparent[1]{}}\providecommand\rotatebox[2]{#2}\newcommand*\fsize{\dimexpr\f@size pt\relax}\newcommand*\lineheight[1]{\fontsize{\fsize}{#1\fsize}\selectfont}\ifx\svgwidth\undefined \setlength{\unitlength}{1187.36678957bp}\ifx\svgscale\undefined \relax \else \setlength{\unitlength}{\unitlength * \real{\svgscale}}\fi \else \setlength{\unitlength}{\svgwidth}\fi \global\let\svgwidth\undefined \global\let\svgscale\undefined \makeatother \begin{picture}(1,0.20397796)\lineheight{1}\setlength\tabcolsep{0pt}\put(0,0){\includegraphics[width=\unitlength,page=1]{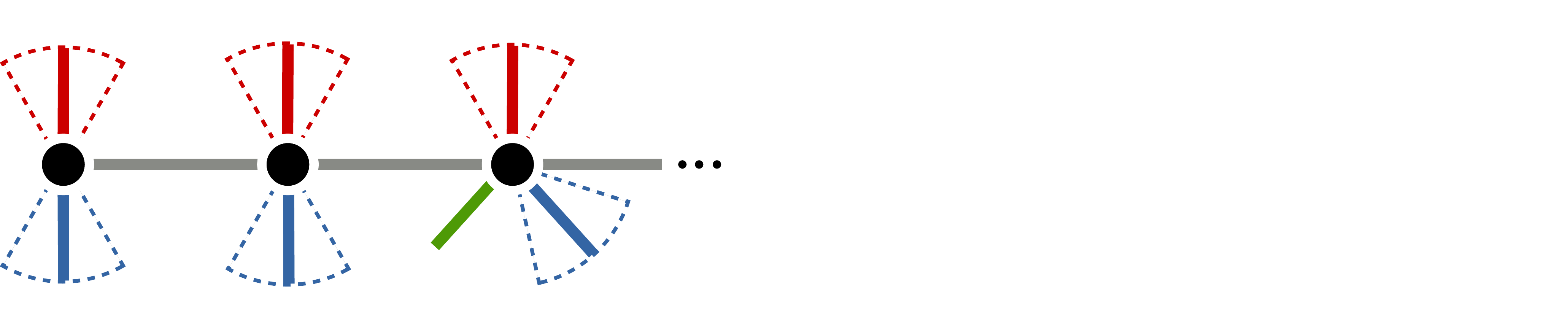}}\put(0.11201694,0.11242072){\color[rgb]{0,0,0}\makebox(0,0)[t]{\lineheight{1.25}\smash{\begin{tabular}[t]{c}$k+1$\end{tabular}}}}\put(0.25534291,0.11262173){\color[rgb]{0,0,0}\makebox(0,0)[t]{\lineheight{1.25}\smash{\begin{tabular}[t]{c}$k+2$\end{tabular}}}}\put(0.18379892,0.19447043){\color[rgb]{0,0,0}\makebox(0,0)[t]{\lineheight{1.25}\smash{\begin{tabular}[t]{c}$r_i$\end{tabular}}}}\put(0.04038848,0.00363958){\color[rgb]{0,0,0}\makebox(0,0)[t]{\lineheight{1.25}\smash{\begin{tabular}[t]{c}$r_i$\end{tabular}}}}\put(0.18366491,0.003128){\color[rgb]{0,0,0}\makebox(0,0)[t]{\lineheight{1.25}\smash{\begin{tabular}[t]{c}$r_i$\end{tabular}}}}\put(0.32751212,0.19447043){\color[rgb]{0,0,0}\makebox(0,0)[t]{\lineheight{1.25}\smash{\begin{tabular}[t]{c}$r_i$\end{tabular}}}}\put(0.37461578,0.01557628){\color[rgb]{0,0,0}\makebox(0,0)[t]{\lineheight{1.25}\smash{\begin{tabular}[t]{c}$r_i$\end{tabular}}}}\put(0.26840897,0.02846313){\color[rgb]{0,0,0}\makebox(0,0)[t]{\lineheight{1.25}\smash{\begin{tabular}[t]{c}$r_{k+2}$\end{tabular}}}}\put(0.04038848,0.1946266){\color[rgb]{0,0,0}\makebox(0,0)[t]{\lineheight{1.25}\smash{\begin{tabular}[t]{c}$r_i$\end{tabular}}}}\put(0,0){\includegraphics[width=\unitlength,page=2]{kBubble.pdf}}\put(0.9597107,0.19360859){\color[rgb]{0,0,0}\makebox(0,0)[t]{\lineheight{1.25}\smash{\begin{tabular}[t]{c}$r_i$\end{tabular}}}}\put(0,0){\includegraphics[width=\unitlength,page=3]{kBubble.pdf}}\put(0.88809056,0.11005178){\color[rgb]{0,0,0}\makebox(0,0)[t]{\lineheight{1.25}\smash{\begin{tabular}[t]{c}$n-1$\end{tabular}}}}\put(0,0){\includegraphics[width=\unitlength,page=4]{kBubble.pdf}}\put(0.56560277,0.1946266){\color[rgb]{0,0,0}\makebox(0,0)[t]{\lineheight{1.25}\smash{\begin{tabular}[t]{c}$r_i$\end{tabular}}}}\put(0.5297927,0.11106977){\color[rgb]{0,0,0}\makebox(0,0)[rt]{\lineheight{1.25}\smash{\begin{tabular}[t]{r}$\ell-1$\end{tabular}}}}\put(0.60141284,0.11106977){\color[rgb]{0,0,0}\makebox(0,0)[lt]{\lineheight{1.25}\smash{\begin{tabular}[t]{l}$\ell$\end{tabular}}}}\put(0.61334953,0.01557626){\color[rgb]{0,0,0}\makebox(0,0)[t]{\lineheight{1.25}\smash{\begin{tabular}[t]{c}$r_i$\end{tabular}}}}\put(0.50591932,0.01557626){\color[rgb]{0,0,0}\makebox(0,0)[t]{\lineheight{1.25}\smash{\begin{tabular}[t]{c}$r_{i+1}$\end{tabular}}}}\put(0.80453374,0.19360859){\color[rgb]{0,0,0}\makebox(0,0)[t]{\lineheight{1.25}\smash{\begin{tabular}[t]{c}$r_i$\end{tabular}}}}\put(0.80453374,0.00262159){\color[rgb]{0,0,0}\makebox(0,0)[t]{\lineheight{1.25}\smash{\begin{tabular}[t]{c}$r_{i+1}$\end{tabular}}}}\put(0.95974458,0.00262159){\color[rgb]{0,0,0}\makebox(0,0)[t]{\lineheight{1.25}\smash{\begin{tabular}[t]{c}$r_{i+1}$\end{tabular}}}}\end{picture}\endgroup  \end{scriptsize}
\caption{The $k$-bubble operator. The vertices (from left to right)  can be labeled with $\ell \in \{k+1, \dots, n\}$ and the voltage of the dart $\pth{\ell}$ is $r_i$ if $i \leq k$ (red), $r_{i+1}$ if $k+1 \leq i \leq \ell-2$ (green), $r_{i}$ if $i > \ell$ (blue)  and trivial otherwise (gray).}
\label{fig:kBubble}
\end{figure}

With this information we can see that the $k$-bubble is described by the voltage operator in Figure\nobreakspace \ref {fig:kBubble}.

\subsubsection*{The $2^\cM$ and $\hat{2}^\cM$ operations}
Given a finite maniplex $\cM$ of rank $n$, the maniplex $2^\cM$ was defined in \cite{douglas2015twist} as follows:
    Label the $0$-faces of $\cM$ with the set $\{1, \dots, \ell \}$.
    The flag set of $2^\cM$ is $\cM\times \bZ_2^{\ell}$.
    Then adjacencies of a flag $(\Psi, v)$, with $\Psi\in\cM$ and $v \in \bZ_2^{\ell}$, are given by the following formula, where $e_{k} \in \bZ_2^{\ell}$ is the vector that has $0$ in all its entries except the $k$-th one:
    \[
    (\Psi,v)^i= \begin{cases}
        (\Psi,v+e_k) & \text{if } i=0, \text{ and the $0$-face of $\Psi$ is labeled with $k$,}\\
        (\Psi^{i-1},v) &\text{otherwise}.
    \end{cases}
    \]

    It is well-known that if $\cM$ is the flag graph of a polytope $\cP$, then $2^\cM$ is the flag graph of the polytope $2^\cP$ defined by Danzer in \cite{danzer1984regular} (see \cite{Mochan_2021_AbstractPolytopesTheir_PhDThesis} for details).

    The maniplex $\hat{2}^\cM$ is defined as the dual of $2^{\cM^{\ast}}$, where $\cM^{\ast}$ is the dual of $\cM$. Thus, the set of flags of $\hat{2}^\cM$ is $\cM \times \bZ_2^{m}$, where the set of facets of $\cM$ is labeled with the set $\{1, \dots, m\}$.
    The adjacencies in $\hat{2}^\cM$ are given by:
     \[
    (\Psi,v)^i= \begin{cases}
        (\Psi^i,v)& \text{if } i<n,\\
        (\Psi,v + e_k) & i=n, \text{ and the facet of $\Psi$ is labeled with $k$.}
    \end{cases}
    \]
     where $e_{k}\in\bZ_2^{m}$ is the vector that has $0$ in all its entries except the $k$-th one.

     These operations on maniplexes are very interesting.
     In particular in the context of voltage operations, Example\nobreakspace \ref {eg:2^maniplex} will show that one cannot find a premaniplex $\cY$ that satisfies that $2^\cM = \cM \ertimes \cY$, for every maniplex $\cM$ (not even if we fix the rank of $\cM$, or ask $\cM$ to be regular).
     However we shall see that if $\cM$ is a regular premaniplex, then there exists a premaniplex $\cY_\cM$ and a voltage assignment   $\eta$ such that $\hat{2}^\cM = \cM \ertimes \cY_\cM$.

     Let $\cM$ be a regular $n$-premaniplex and let $\rho_0, \rho_1, \dots, \rho_{n-1}$ be the distinguished generators of $\aut(\cM)$ with respect to the base flag $\Phi$.
     We define the premaniplex $\cY:=\cY_\cM$ as follows.
     The set of vertices of $\cY$ is the set $\bZ_2^{m}$.
     Given $v \in \cY$ the entries of $v$ correspond to facets of $\cM$. Assume that the facets of $\cM$ are labeled with the set $\{1, \dots, m\}$, as above, in such a way that the facet corresponding to the base flag is labeled with $1$.
     Notice that every element in $\aut(\cM)$ permutes the facets of $\cM$, which in turn induces a action of $\aut(\cM)$ on $\bZ_2^{m}$ by permutation of coordinates. More precisely, for $\alpha \in \aut(\cM)$ and $v=(v_1, \dots v_{m}) \in \bZ_2^{m}$, let $v \alpha$ denote the vector $(v_{1\alpha^{-1}}, \dots, v_{m\alpha^{-1}})$, that is, the vector resulting from $v$ after permuting the coordinates according to $\alpha$.

     The adjacencies in $\cY$ are given as follows: for $i\in \{0,1,\dots, n-1\}$ there is an edge of color $i$ between $v$ and $v^i:=v \rho_{i}$.
Further, there is an edge of  color $n$ between $v$ and $v+e_{1}$.
     We now define the voltage assignment on $\cY$ as:
        \[
    \eta(^iv)= \begin{cases}
        r_i& \text{if } i<n,\\
        1 &i=n.
    \end{cases}
    \]
    We shall show that $\hat{2}^\cM \cong\cM \ertimes \cY$.
     For this, we recall that $\cM$ is regular, and hence any flag of $\cM$ can be written in a unique way as $\Phi\alpha$, where $\Phi$ is the base flag and $\alpha \in \aut(\cM)$.

    Consider $\varphi:\hat{2}^\cM \to \cM \ertimes \cY$ defined as \[ \varphi(\Phi\alpha, v)=(\Phi\alpha, v \alpha^{-1}).\] Then, $\varphi$ preserves incidences.
    In fact, we observe that if $k= 1\alpha$:
    \[\begin{aligned}
        \Big(\varphi(\Phi\alpha, v)\Big)^n
        &=\Big(\Phi\alpha, v \alpha^{-1} \Big)^n \\
        &= \Big( \eta(^n(v \alpha^{-1}))(\Phi\alpha), ( v \alpha^{-1} )^n \Big) \\
        &= \Big( \Phi\alpha, v \alpha^{-1}  + e_{1} \Big) \\
        &= \Big( \Phi\alpha, (v + e_{1 \alpha}) \alpha^{-1}  \Big) \\
        &= \Big( \Phi\alpha, (v+e_{k}) \alpha^{-1}\Big)\\
        &=  \varphi \Big( \Phi\alpha, v+e_{k}\Big) \\
        &= \varphi\Big( (\Phi\alpha,v)^n \Big).
    \end{aligned}\]
On the other hand, for $i<n$ we have that:
\[\begin{aligned}
    \Big(\varphi( \Phi\alpha, v)\Big)^i
    &= \Big( \Phi\alpha, v \alpha^{-1}\Big)^i\\
    &= \Big( \eta\big(^i(v \alpha^{-1})\big)(\Phi\alpha), (v\alpha^{-1})^i\Big)\\
    &= \Big( r_i\Phi\alpha, (v\alpha^{-1})^i\Big)\\
    &= \Big( \Phi^i\alpha,  v\alpha^{-1}\rho_i \Big)\\
    &= \Big( \Phi^i\alpha, v (\rho_i \alpha)^{-1}  )\Big)\\
    &= \Big( \Phi \rho_{i} \alpha, v (\rho_i \alpha)^{-1}  )\Big)\\
    &= \varphi \Big( \Phi\rho_i\alpha, v\Big)\\
    &= \varphi \Big( \big((\Phi\alpha)^i, v\big)\Big) \\
    &= \varphi\Big( (\Phi\alpha,v)^i \Big).
\end{aligned}\]

The computation above shows that $\hat{2}^\cM$ and $\cM \ertimes \cY$ are isomorphic.

It is clear that if we drop the regularity condition on $\cM$ the given voltage operator $(\cY, \eta)$ is not well-defined as it depends heavily on the existence of the distinguished generators of the automorphism group.
Thus, we can ask ourselves the following questions. 

\begin{question}\label{Q1}
Is it true that for any maniplex $\cM$, there exists a voltage operator $(\cY_\cM,\eta)$ such that $\hat{2}^\cM \cong \cM \ertimes \cY_\cM$?
\end{question}

For example, we have not been able to determine if there exists such voltage operator for the case when $\cM$ is the pyramid over a digon.
If the answer for the above question is negative, then we could ask:

\begin{question}
 Give necessary and sufficient conditions on a maniplex $\cM$ in such a way that there exists a voltage operator $(\cY_\cM,\eta)$ such that $\hat{2}^\cM \cong \cM \ertimes \cY_\cM$.
\end{question}

Note that a sufficient condition for the last question is to have $\cM$ to be regular, as shown above. However chances are this condition is not necessary. For example, we believe that if the symmetry type graph of $\cM$ is itself a regular pre-maniplex, the above voltage operator might exist.
 \section{Automorphisms}
\label{sec:automorphisms}

In this section we will study the interplay between the automorphism groups of a premaniplex and the resulting premaniplex after a voltage operation, as well as their symmetry type graphs.

Let $\cX$ be a premaniplex and let $(\cY, \eta)$ be a voltage operator.
Observe that the elements of $\aut(\cX)$ act as automorphisms on $\cX \ertimes \cY$.
Indeed, if $\gamma \in \aut(\cX)$, then the mapping $\bar{\gamma}: (x,y) \mapsto (x\gamma, y)$ induces an automorphism of $\cX \ertimes \cY$.
To see this, just note that $\bar{\gamma}$ commutes with the monodromies\[\begin{aligned}
  r_{i} \left( (x,y) \bar{\gamma} \right)
	&= r_{i}\left( x\gamma,y \right) \\
	&= \left( \eta (^{i}y)x \gamma, r_{i}y \right) \\
	&= \left( \eta(^{i} y)x, r_{i}y \right)\bar{\gamma}\\
	&= \left( r_{i}(x,y) \right) \bar{\gamma}.
\end{aligned}\]

Thus, $\aut(\cX) \leq \aut(\cX \ertimes \cY)$ for every voltage operator $(\cY, \eta)$ and we may replace $\bar{\gamma}$ by $\gamma$ to simplify notation.
The latter result implies that $\aut(\cX)$ acts by automorphism on $\cX \ertimes \cY$.
Note that $\cX$ naturally covers $\cX/\Gamma$ for every group $\Gamma \leq \aut(\cX)$.
As expected $\cX \ertimes Y$ naturally covers $(\cX \ertimes Y)/\Gamma$,
In fact, we have slightly more general result:

\begin{thm}\label{thm:covers}
Let $\cX$ and $\cX'$ be $n$-premaniplexes, $(\cY,\eta)$.
If $\cX$ covers $\cX'$, then $\cX\ertimes \cY$ covers $\cX'\ertimes\cY$.

In particular, when $\cX'$ is obtained from $\cX$ as a quotient by a group $\Gamma \leq \aut(\cX)$,
we get that \[\cX' \ertimes Y \cong \left( \cX / \Gamma \right) \ertimes \cY \cong \left( \cX \ertimes \cY \right)/ \Gamma.\]
\end{thm}
\begin{proof}
    Let $\varphi:\cX\to\cX'$ be a covering homomorphism. If $\cX'=\cX/\Gamma$ we chose $\varphi$ as the natural projection $x\mapsto x\Gamma$.
    Define $\varphi' : \cX \ertimes \cY  \to \cX'\ertimes \cY$ by $\varphi': (x,y) \mapsto (\varphi(x), y)$.

    Now, using the fact that $\varphi$ is a homomorphism we simply check that
    \[
    \begin{aligned}
        r_i\left(\varphi'(x,y)\right) &= r_i\left(\varphi(x),y\right)\\
        &= \left(\eta(\pth{y})\varphi(x), r_i y\right)\\
        &= \left(\varphi(\eta(\pth{y})x), r_i y \right)\\
        &= \varphi'\left( \eta(\pth{y})x, r_i y \right)\\
        &= \varphi' \left(r_i(x,y)\right).
    \end{aligned}
    \]

    This proves that $\varphi'$ is a homomorphism, and it is trivial to see that it is surjective. Therefore it is a covering.
    
    In the particular case when $\cX'=\cX/\Gamma$ we get that $\varphi'(x,y) = (x\Gamma,y)$, but by using the natural action of $\Gamma$ on $\cX\ertimes \cY$ we get that $(x\Gamma,y) = (x,y)\Gamma$.
\end{proof}

The above theorem is of particular interest when $\cX$ is a maniplex and $\Gamma =\aut(\cX)$. In this case, $\cT:=\cX/\Gamma$ is the symmetry type graph of $\cX$. Theorem\nobreakspace \ref {thm:covers} tells us that the symmetry type graph of $\cX \ertimes \cY$ with respect to $\Gamma$ is precisely $\cT \ertimes \cY$.
However it is important to remark that $\Gamma$ might be a proper subgroup of $\aut(\cX \ertimes \cY)$, that is, $\cT \ertimes \cY$ might not be the symmetry type graph of $\cX \ertimes \cY$, even when $\cX \ertimes \cY$ is a maniplex.
An example of such situation will be given in Example\nobreakspace \ref {eg:PyramidMedialVolts}.
This phenomenon will be explored deeply in \cite{HubardMochanMontero_MoreVoltageOperations_preprint}.

\emph{Words of caution}: the fact that $\aut(\cX)$ acts by automorphisms on $\cX\ertimes\cY$ {\em does not} imply that it acts by automorphisms on each connected component of $\cX\ertimes\cY$.
For example, as pointed out in Section\nobreakspace \ref {sec:ClassicalExamples}, when
 $\cX$ is an orientable (pre)maniplex of rank 3 and $(\cY,\eta)$ is the snub operator of Figure\nobreakspace \ref {fig:snub}, $\cX\ertimes \cY$ has two connected components, say  the \emph{left snub} (whose flags are those of the form $(x,y)$ with $x\in\cX$ a white flag), and the \emph{right snub} (whose flags are those of the form $(x,y)$ with $x\in\cX$ a black flag). Each component has all the ``rotational'' symmetry of $\cX$ but none of the ``reflection'' symmetries of $\cX$.
However, if there exists a reflection $\tau \in \aut(\cX)$, that is, an automorphism that interchanges white flags with black flags, then $\tau$ acts on $\cX\ertimes \cY$ by swapping the left snub with the right snub.  This means that if we consider rooted premaniplexes $(\cX,x)$ and $(\cY,y)$, then the rooted premaniplex $(\cX\ertimes \cY, (x,y))$ may not have all the symmetries of $\cX$.
In fact we have the following result:

\begin{prop}
Let $(\cX,x)$ be a rooted premaniplex, $((\cY,y),\eta)$ a rooted voltage operator and $\gamma\in\aut(\cX)$. Then $\gamma\in\aut(\cX\ertimes \cY, (x,y))$ if and only if there exists a closed path $W$ based on $y$ such that $x\gamma = \eta(W)x$.
\end{prop}

\begin{proof}
 Assume that $\gamma\in\aut(\cX\ertimes \cY, (x,y))$. This means that $(x,y)\gamma$ is an element of  $(\cX\ertimes \cY, (x,y))$; in other words, $(x,y)\gamma$ and $(x,y)$ are in the same connected component of $\cX\ertimes \cY$. Hence, there exists a monodromy $\omega$ such that $\omega(x,y) = (x,y)\gamma$. Note that on one hand, $(x,y)\gamma = (x\gamma, y)$, while on the other hand $\omega(x,y) = (\eta(\pth[\omega]{y})x, \omega y)$. Thus, $(x\gamma, y) = (\eta(\pth[\omega]{y})x, \omega y)$, which implies that $y=\omega y$ and therefore  $W:= \pth[\omega]{y} \in \fg^y(\cY)$.

 Conversely, suppose that there exists $W\in \fg^y(\cY)$ such that $x\gamma = \eta(W)x$. The fact that  $W\in \fg^y(\cY)$ implies that there exists a monodromy $\omega$ such that $\omega y = y$ and $W=\pth[\omega]{y}$. Hence, $(x,y)\gamma = (x\gamma, y) = (\eta(W)x, y) = (\eta(\pth[\omega]{y})x, y) = (\eta(\pth[\omega]{y})x, \omega y) = \omega(x,y)$, which implies that $(x,y)\gamma$ and $(x,y)$ are in the same connected component of $\cX\ertimes \cY$. Since we already know that $\gamma$ is an automorphism of $\cX\ertimes \cY$, this implies that $\gamma\in\aut(\cX\ertimes \cY, (x,y))$.
\end{proof}

Note that we may use Theorem\nobreakspace \ref {thm:covers} to determine when an operation is not a voltage operation, as we see in the following example.

\begin{eg}\label{eg:2^maniplex}

In this example we will prove that there is no $(n,n+1)$-voltage operator $(\cY,\eta)$ such that $2^\cM = \cM\ertimes \cY$ for every $n$-maniplex $\cM$ (as was hinted in \protect \MakeUppercase {Q}uestion\nobreakspace \ref {Q1}).
    For a $2$-gon $\{2\}$, the polytope $2^{\{2\}}$ is the quadrangular dihedron $\{4,2\}$ (the map on the sphere consisting of two squared faces sharing all the vertices and edges).
    For the square $\{4\}$, we have that $2^{\{4\}}$ is the toroidal map $\{4,4\}_{(4,0)}$ (a $4\times 4$ grid on the torus).
    We can get a $2$-gon from a square by taking its quotient by the cyclic group of order two generated by the half-turn around the center, that is $\{2\}=\{4\}/\bZ_2$.
    However, it is not possible to get the quadrangular dihedron by taking a quotient of $\{4,4\}_{(4,0)}$ by any group of order 2.
    In other words $2^{\{4\}} / \bZ_2 \ncong 2^{\{4\}/\bZ_2}$. Hence, Theorem\nobreakspace \ref {thm:covers} is not satisfied for the operation $2^\cM$.
    Thus, there is no voltage operator $(\cY,\eta)$ such that $2^\cM = \cM \ertimes \cY$ for every maniplex $\cM$.
    \end{eg}

    The previous example uses Theorem\nobreakspace \ref {thm:covers} to prove that there does not exist a voltage operation $(\cY, \eta)$ such that $2^\cM \cong \cM\ertimes \cY$ for all maniplexes $\cM$, in other words, that
    $\cM\mapsto 2^\cM$ is not a voltage operation.
    Observe that this fact can also be established in a more elementary way by simply comparing the number of flags in $\cM \ertimes \cY$ and in $2^{\cM}$.
    We use Theorem\nobreakspace \ref {thm:covers} to emphasize the relation between quotients and voltage operations.
    In fact, Theorem\nobreakspace \ref {thm:covers} characterizes all voltage operations; this shall be shown in Theorem\nobreakspace \ref {thm:OpsAsProducts}.

   It is easy to come up with examples of operations that are not voltage operations, for instance, by defining them differently for different cases. However, the operation $\cM\mapsto 2^\cM$ is of interest because it is a functor between the categories $\pman^n$ and $\pman^{n+1}$,
in other words, if there is a homomorphism $p:\cM \to \cM'$ then there is a natural way to define a homomorphism $\tilde{p}:2^\cM\to 2^{\cM'}$,
as $\tilde{p}(\Phi,v)\mapsto (p(\Phi),u)$ where \[u_i=\sum_{j\in p^{-1}(i)} v_j.\] It is straightforward prove that $\tilde{p}$ is a maniplex homomorphism and that $p\mapsto \tilde{p}$ is indeed a functor.

Let $\cM$ be a premaniplex, let $\Gamma$ be a group of automorphisms of $\cM$ and let $\cX=\cM/\Gamma$.
If $(\cY,\eta)$ is a voltage operation, Theorem\nobreakspace \ref {thm:covers} tells us that $(\cM\ertimes \cY)/\Gamma$ is isomorphic to $\cX \ertimes \cY$.
Hence, there exists a voltage assignment $\xi$ on $\cX$ with voltage group $\Gamma$ such that $\cX^\xi$ is isomorphic to $\cM$, and
 a voltage assignment $\theta:\fg(\cX \ertimes \cY)\to \Gamma$ such that $(\cX \ertimes \cY)^\theta$ is isomorphic to $\cM \ertimes \cY$.
The following theorem tells us how to find $\theta$ in terms of  $\eta$ and $\xi$.

\begin{thm}\label{thm:theta}
Let $(\cX,\xi)$ be a voltage premaniplex with voltage group $\Gamma$ and let $(\cY,\eta)$ be a voltage operator.

Define $\theta = \theta(\eta,\xi):\Pi(\cX \ertimes \cY)\to \Gamma$ as follows:
\begin{equation}\label{eq:theta}
    \theta(\pth[\omega]{(x,y)}):=\xi(\pth[{\eta(\pth[\omega]{y})}]{x}).
\end{equation}

Then $(\cX \ertimes \cY)^\theta$ is isomorphic to $\cX^\xi \ertimes \cY$.
\end{thm}
\begin{proof}
Start by noticing that, by definition, $r_i(x,\gamma) = (r_ix, \xi(^ix)\gamma)$ and similarly $r_i((x,y),\gamma) = (r_i(x,y), \theta(^i(x,y))\gamma)$.
    Now, define the function $$\varphi:(\cX \ertimes \cY)^\theta \to \cX^\xi \ertimes \cY$$ given by
    $$\varphi((x,y),\gamma):=((x,\gamma),y).$$We shall show that $\varphi$ is a premaniplex isomorphism.

    Since both $\cX^\xi$ and $(\cX \ertimes \cY)^\theta$ have $\Gamma$ as voltage group, $\varphi$ is a bijection. We should see next that $\varphi$ preserves $i$\-adjacencies for $i\in\{0,1,\ldots,m-1\}$. In fact:

    \[
    \begin{aligned}
\varphi\left(r_i((x,y),\gamma)\right)
            &= \varphi\left(r_i(x,y),\theta(\pth{(x,y)})\gamma\right)\\
            &= \varphi\left((\eta(\pth{y})x,r_i y),\theta(\pth{(x,y)})\gamma\right)\\
            &= \left((\eta(\pth{y})x,\theta(\pth{(x,y)})\gamma),r_i y\right)\\
            &= \left((\eta(\pth{y})x,\xi(\pth[\eta(\pth{y})]{x})\gamma),r_i y\right)\\
           &= \left(\eta(\pth{y})(x,\gamma), r_i y\right)\\
            &= r_i \left((x,\gamma),y\right)\\
           &= r_i \varphi \left((x,y),\gamma\right).
\end{aligned}
    \]
\end{proof}

\begin{figure}
	\centering
    \begin{scriptsize}
\def\svgwidth{.8\textwidth}
	\begingroup \makeatletter \providecommand\color[2][]{\errmessage{(Inkscape) Color is used for the text in Inkscape, but the package 'color.sty' is not loaded}\renewcommand\color[2][]{}}\providecommand\transparent[1]{\errmessage{(Inkscape) Transparency is used (non-zero) for the text in Inkscape, but the package 'transparent.sty' is not loaded}\renewcommand\transparent[1]{}}\providecommand\rotatebox[2]{#2}\newcommand*\fsize{\dimexpr\f@size pt\relax}\newcommand*\lineheight[1]{\fontsize{\fsize}{#1\fsize}\selectfont}\ifx\svgwidth\undefined \setlength{\unitlength}{1602.27221872bp}\ifx\svgscale\undefined \relax \else \setlength{\unitlength}{\unitlength * \real{\svgscale}}\fi \else \setlength{\unitlength}{\svgwidth}\fi \global\let\svgwidth\undefined \global\let\svgscale\undefined \makeatother \begin{picture}(1,0.53772463)\lineheight{1}\setlength\tabcolsep{0pt}\put(0,0){\includegraphics[width=\unitlength,page=1]{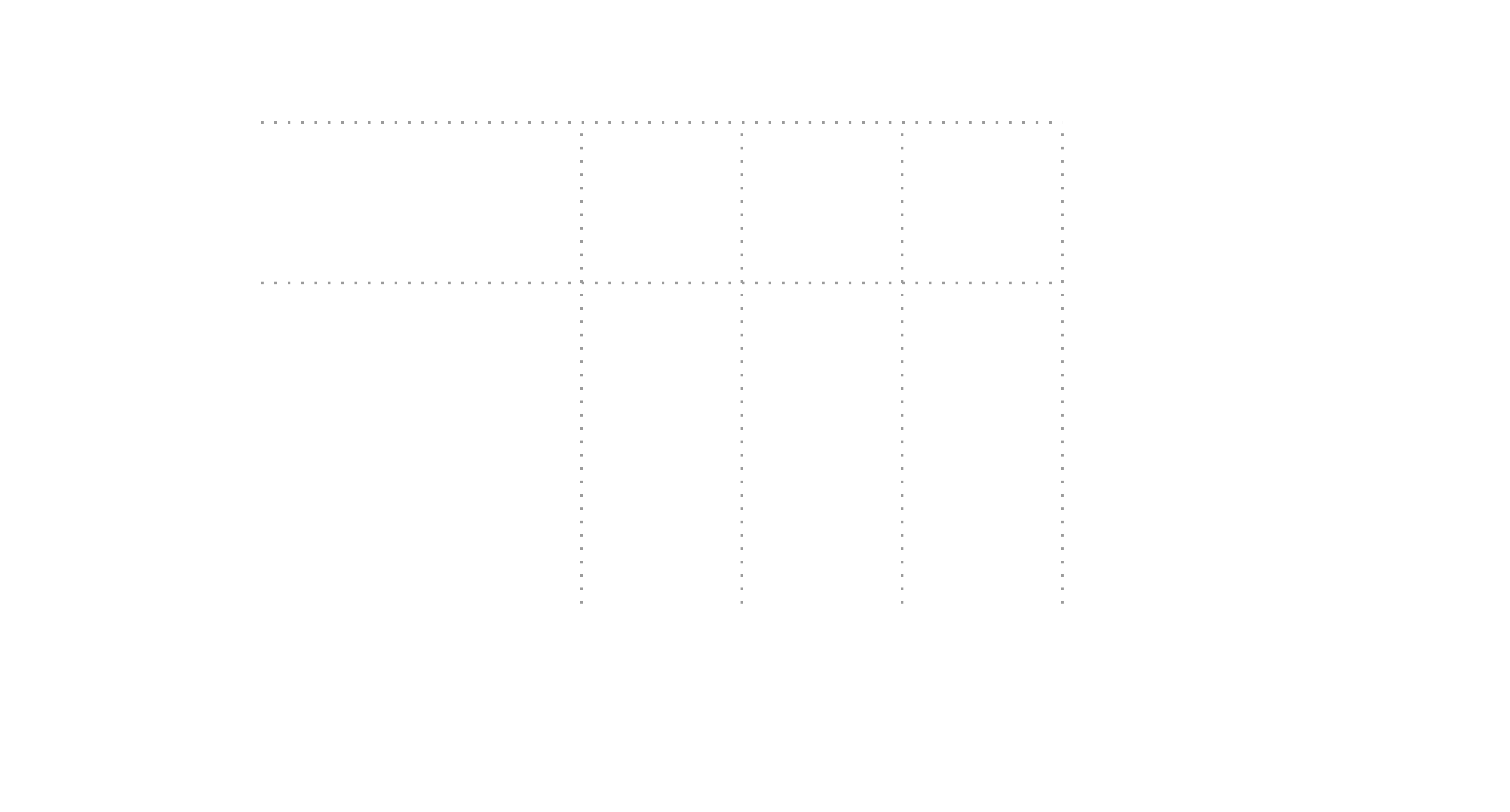}}\put(0.09936723,0.3983948){\color[rgb]{0,0,0}\makebox(0,0)[rt]{\lineheight{1.25}\smash{\begin{tabular}[t]{r}$(\cY, \eta)$\end{tabular}}}}\put(0.54445395,0.00918375){\color[rgb]{0,0,0}\makebox(0,0)[t]{\lineheight{1.25}\smash{\begin{tabular}[t]{c}$(\cX, \xi)$\end{tabular}}}}\put(0.78107944,0.3983948){\color[rgb]{0,0,0}\makebox(0,0)[lt]{\lineheight{1.25}\smash{\begin{tabular}[t]{l}$(\cX \ertimes \cY, \theta)$\end{tabular}}}}\put(0,0){\includegraphics[width=\unitlength,page=2]{theta_pyr.pdf}}\put(0.42940949,0.51838064){\color[rgb]{0,0,0}\makebox(0,0)[t]{\lineheight{1.25}\smash{\begin{tabular}[t]{c}$\rho_1$\end{tabular}}}}\put(0.34095243,0.51838064){\color[rgb]{0,0,0}\makebox(0,0)[t]{\lineheight{1.25}\smash{\begin{tabular}[t]{c}$\rho_0$\end{tabular}}}}\put(0,0){\includegraphics[width=\unitlength,page=3]{theta_pyr.pdf}}\put(0.38518094,0.27954656){\color[rgb]{0,0,0}\makebox(0,0)[t]{\lineheight{1.25}\smash{\begin{tabular}[t]{c}$\rho_1$\end{tabular}}}}\put(0.59747789,0.27954656){\color[rgb]{0,0,0}\makebox(0,0)[t]{\lineheight{1.25}\smash{\begin{tabular}[t]{c}$\rho_1$\end{tabular}}}}\put(0.74785491,0.27954656){\color[rgb]{0,0,0}\makebox(0,0)[t]{\lineheight{1.25}\smash{\begin{tabular}[t]{c}$\rho_1$\end{tabular}}}}\put(0.49132942,0.51838064){\color[rgb]{0,0,0}\makebox(0,0)[t]{\lineheight{1.25}\smash{\begin{tabular}[t]{c}$\rho_0$\end{tabular}}}}\put(0.70362644,0.51838064){\color[rgb]{0,0,0}\makebox(0,0)[t]{\lineheight{1.25}\smash{\begin{tabular}[t]{c}$\rho_0$\end{tabular}}}}\put(0.65939791,0.27954656){\color[rgb]{0,0,0}\makebox(0,0)[t]{\lineheight{1.25}\smash{\begin{tabular}[t]{c}$\rho_0$\end{tabular}}}}\put(0,0){\includegraphics[width=\unitlength,page=4]{theta_pyr.pdf}}\put(0.11980977,0.27954659){\color[rgb]{0,0,0}\makebox(0,0)[t]{\lineheight{1.25}\smash{\begin{tabular}[t]{c}$r_1$\end{tabular}}}}\put(0.22595823,0.27954659){\color[rgb]{0,0,0}\makebox(0,0)[t]{\lineheight{1.25}\smash{\begin{tabular}[t]{c}$r_2$\end{tabular}}}}\put(0,0){\includegraphics[width=\unitlength,page=5]{theta_pyr.pdf}}\put(0.22595823,0.51838063){\color[rgb]{0,0,0}\makebox(0,0)[t]{\lineheight{1.25}\smash{\begin{tabular}[t]{c}$r_1$\end{tabular}}}}\put(0.11980977,0.51838063){\color[rgb]{0,0,0}\makebox(0,0)[t]{\lineheight{1.25}\smash{\begin{tabular}[t]{c}$r_0$\end{tabular}}}}\put(0,0){\includegraphics[width=\unitlength,page=6]{theta_pyr.pdf}}\put(0.17290844,0.4515806){\color[rgb]{0,0,0}\makebox(0,0)[t]{\lineheight{1.25}\smash{\begin{tabular}[t]{c}$x$\end{tabular}}}}\put(0,0){\includegraphics[width=\unitlength,page=7]{theta_pyr.pdf}}\put(0.17290844,0.34543215){\color[rgb]{0,0,0}\makebox(0,0)[t]{\lineheight{1.25}\smash{\begin{tabular}[t]{c}$y$\end{tabular}}}}\put(0,0){\includegraphics[width=\unitlength,page=8]{theta_pyr.pdf}}\put(0.32326096,0.07609543){\color[rgb]{0,0,0}\makebox(0,0)[rt]{\lineheight{1.25}\smash{\begin{tabular}[t]{r}$\rho_1$\end{tabular}}}}\put(0.32326096,0.19108958){\color[rgb]{0,0,0}\makebox(0,0)[rt]{\lineheight{1.25}\smash{\begin{tabular}[t]{r}$\rho_0$\end{tabular}}}}\put(0.49132942,0.20878099){\color[rgb]{0,0,0}\makebox(0,0)[t]{\lineheight{1.25}\smash{\begin{tabular}[t]{c}$\rho_0$\end{tabular}}}}\put(0.59747772,0.05840394){\color[rgb]{0,0,0}\makebox(0,0)[t]{\lineheight{1.25}\smash{\begin{tabular}[t]{c}$\rho_1$\end{tabular}}}}\put(0.76554613,0.07609543){\color[rgb]{0,0,0}\makebox(0,0)[lt]{\lineheight{1.25}\smash{\begin{tabular}[t]{l}$\rho_1$\end{tabular}}}}\put(0.76554613,0.19993529){\color[rgb]{0,0,0}\makebox(0,0)[lt]{\lineheight{1.25}\smash{\begin{tabular}[t]{l}$\rho_0$\end{tabular}}}}\put(0.3852053,0.13313528){\color[rgb]{0,0,0}\makebox(0,0)[t]{\lineheight{1.25}\smash{\begin{tabular}[t]{c}$a$\end{tabular}}}}\put(0.49135382,0.13313528){\color[rgb]{0,0,0}\makebox(0,0)[t]{\lineheight{1.25}\smash{\begin{tabular}[t]{c}$b$\end{tabular}}}}\put(0.59750224,0.13313528){\color[rgb]{0,0,0}\makebox(0,0)[t]{\lineheight{1.25}\smash{\begin{tabular}[t]{c}$c$\end{tabular}}}}\put(0.70365059,0.13313528){\color[rgb]{0,0,0}\makebox(0,0)[t]{\lineheight{1.25}\smash{\begin{tabular}[t]{c}$d$\end{tabular}}}}\end{picture}\endgroup  	\end{scriptsize}
	\caption{A voltage premaniplex for the medial of a pyramid, where the colors red, green and blue for the edges represent $0$-, $1$- and $2$-adjacencies, respectively.}\label{fig:theta_pyr}
\end{figure}

\begin{eg}\label{eg:PyramidMedialVolts}
It is not difficult to see that the antiprism of a $q$-gon can be obtained by taking the medial of the pyramid over a $q$-gon. Thus, we shall use Theorem\nobreakspace \ref {thm:theta} to recover the antiprism of a $q$-gon as a voltage maniplex. To do so, we first construct the pyramid over a $q$-gon as a voltage maniplex and then apply the medial operator (as a voltage operator).

    Note that the automorphism group of a $q$-gonal pyramid coincides with the automorphism group of its base, which is the dihedral group \[ \bD_q = \gen{\rho_0,\rho_1\mid \rho_0^2=\rho_1^2=(\rho_0\rho_1)^q=1},\] where $\rho_0$ is thought as the reflection in a plane orthogonal to the base, that includes the midpoint of an edge $e$ in the base; and $\rho_1$ as the reflection in a plain orthogonal to the base, that includes one of the vertices incident to $e$.
Hence, the $q$-gonal pyramid can be recovered from its symmetry type graph $\cX$ via the voltage premaniplex $(\cX,\xi)$ shown in Figure\nobreakspace \ref {fig:theta_pyr}.
   (As in previous examples, the colors red, green and blue for the edges represent $0$-, $1$- and $2$-adjacencies, respectively.)

   Let $(\cY,\eta)$ be the medial operator. Then, by Theorem\nobreakspace \ref {thm:theta}, $(\cX\ertimes \cY,\theta)$ is the antiprism over a $q$-gon, where $\theta$ is defined as in the theorem.

 Theorem\nobreakspace \ref {thm:theta} gives us a way to construct such antiprism.
When we want to calculate, for example, the endpoint and voltage of the dart $\pth[0]{(a,x)}$ (the red dart that starts at $(a,x)$), we first see that $\eta(\pth[0]{x}) = r_1$,  so we follow the path $\pth[1]{a}$ and we see that it ends in $a$ and has voltage $\xi(\pth[1]{a}) = \rho_1$. Therefore, $(a,x)^0 = (\eta(\pth[0]{x})a,x^i) = (a,x)$ (i.e., the dart is a semiedge) and the dart $\pth[0]{(a,x)}$ has voltage $\theta(\pth[0]{(a,x)}) = \rho_1$.

    Note that $\cX\ertimes \cY$ is the symmetry type graph of the antiprism with respect to the automorphism group of its base, but not with respect to the full automorphism group of the antiprism.
    In fact, the antiprism over any polygon always has extra symmetry, which is induced by the isomorphism between the base of the antiprism and its dual, so antiprisms over polygons have usually four orbits on flags.
    The antiprism over the triangle is in fact an octahedron and therefore regular, so  in that case the actual symmetry type graph is $\1^3$.\end{eg}

An consequence of Theorem\nobreakspace \ref {thm:theta} is the following corollary:

\begin{coro}\label{coro:OperatingRegulars}
  Let $\cX$ be a regular $n$-premaniplex with automorphism group $\langle \rho_{0}, \dots, \rho_{n-1} \rangle$ and let $(\cY,\eta)$ be a $(n,m)$-voltage operator.
  Then $\cX\ertimes \cY$ is isomorphic to the derived graph $\cY^\nu$, where $\nu:\fg(\cY)\to\aut(\cX)$ is the voltage assignment obtained from $\eta$ by replacing each $r_i$ with $\rho_i$.
\end{coro}
\begin{proof}
Since $\cX$ is regular, it is isomorphic to the derived graph $(\1^n)^\xi$ where, if $x$ is the only vertex of $\1^n$, $\xi(\pth{x})=\rho_i$.
Then Theorem\nobreakspace \ref {thm:theta} tells us that $\cX\ertimes \cY$ is isomorphic to $(\1^n\ertimes \cY)^\theta$ with $\theta(\pth[\omega]{(x,y)}) = \xi (\pth[{\eta(\pth[\omega]{y}}]){x})$.
This means precisely that $\theta$ replaces each occurrence of $r_i$ in $\eta$ by the voltage of the semiedge of color $i$ in $\1^n$, but this is exactly $\rho_i$.
By applying the natural isomorphism $\1^n \ertimes \cY\to \cY$ we get the desired result.
\end{proof}

\section{Composition of voltage operations}\label{sec:composition}

The result of applying a voltage operation to a premaniplex is again a premaniplex. Thus, it is natural to think about the composition of two voltage operations. In contrast, the result of applying a voltage operation to a maniplex is not always a maniplex, so one must be careful with this fact when composing operations, as the result of a voltage operation can be disconnected.

It is interesting to note that in fact the composition of two voltage operations can be written as a new voltage operation.
In this section we describe how to do this by using the operator $\theta$ defined in Theorem\nobreakspace \ref {thm:theta}, but instead of using an arbitrary voltage premaniplex $(\cX,\xi)$ and a voltage operator $(\cY,\eta)$ we use two voltage operators $(\cY_1,\eta_1)$ and $(\cY_2,\eta_2)$.

\begin{thm}\label{thm:composition}
    Let $\cX$ be an $n$-premaniplex, $(\cY_1,\eta_1)$ an $(n,m)$-voltage operator and $(\cY_2,\eta_2)$ a $(m,\ell)$-voltage operator.
    Then $$(\cX \ertimes[\eta_1] \cY_1)\ertimes[\eta_2] \cY_2 \cong \cX \ertimes[\theta] (\cY_1\ertimes[\eta_2] \cY_2),$$ where $\theta=\theta(\eta_2,\eta_1)$ is defined as in Equation\nobreakspace \textup {(\ref {eq:theta})}.
\end{thm}
\begin{proof}
    We define $\varphi:(\cX \ertimes[\eta_1] \cY_1)\ertimes[\eta_2] \cY_2 \to \cX \ertimes[\theta] (\cY_1\ertimes[\eta_2] \cY_2)$ in the natural way, that is, $$\varphi\left((x,y_1),y_2\right) = \left(x,(y_1,y_2)\right).$$
    It is clear that $\varphi$ is a bijection, so we only need to prove that it commutes with $r_i$, for all $i\in\{0,1,\ldots,\ell-1\}$:
    \[
    \begin{aligned}
        \varphi\left(r_i\left((x,y_1),y_2\right)\right)
        &= \varphi\left(\eta_2(\pth{y_2})(x,y_1), r_i y_2\right)\\
        &= \varphi\left((\eta_1(\pth[{\eta_2(\pth{y_2})}]{y_1})x,\eta_2(\pth{y_2})y_1), r_i y_2\right)\\
        &= \varphi\left((\theta(\pth{(y_1,y_2)}x,\eta_2(\pth{y_2})y_1), r_i y_2\right)\\
        &= \left(\theta(\pth{(y_1,y_2)}x,(\eta_2(\pth{y_2})y_1, r_i y_2)\right)\\
        &= \left(\theta(\pth{(y_1,y_2)}x,r_i(y_1, y_2)\right)\\
        &= r_i \left( x,(y_1,y_2)\right)\\
        &= r_i \varphi \left((x,y_1),y_2\right)
    \end{aligned}
    \]
\end{proof}

The above theorem can be apply in different contexts, we give some examples here. Let $(\1^m,d)$ be the dual operator and consider a $(n,m)$-voltage operator $(\cY,\eta)$.
If $\cX$ is a premaniplex, then $(\cX\ertimes \cY)\ertimes[d] \1^m$ is the dual of $(\cX\ertimes \cY)$.
Theorem\nobreakspace \ref {thm:composition} tells us that this dual is in fact isomorphic to $\cX\ertimes[\theta](\cY\ertimes[d]\1^m)$. Thus, $(\cY\ertimes[d]\1^m,\theta)$ is the dual of $\cY$, where the voltages of the darts are the same as the ones of the darts of the dual color in $(\cY,\eta)$.
In other words, if  $(\cY^*,\eta^*)$ denotes the voltage operator we get by recoloring the darts of color $i$ in $\cY$ with the color $n-1-i$, then the dual of $\cX\ertimes \cY$ is $\cX\ertimes[\eta^*]\cY^*$.

More generally, Theorem\nobreakspace \ref {thm:composition} lets us define the composition of two operators.
If $(\cY_1,\eta_1)$ is an $(n,m)$-operator and $(\cY_2,\eta_2)$ is an $(m,\ell)$-operator we define the \emph{composition of $(\cY_1,\eta_1)$ with $(\cY_2,\eta_2)$} as $(\cY_1\ertimes[\eta_2] \cY_2,\theta(\eta_1,\eta_2))$ and denote it by $(\cY_1,\eta_1)\circ (\cY_2,\eta_2)$.
Theorem\nobreakspace \ref {thm:composition} tells us that $(\cY_1,\eta_1)\circ (\cY_2,\eta_2)$ is an $(n,\ell)$-operator and that the composition of operators is associative.
This allows us to define a new category:
recall that $\pman^n$ denotes the class of all premaniplexes of rank $n$ and let $\bar{\pman}$ be the category whose objects are the classes $\pman^n$ with $n \geq 1$, and whose arrows are voltage operators.
An $(n,m)$-operator is an arrow from $\pman^n$ to $\pman^m$ and the composition of arrows is defined as above.
The neutral element at the object $\pman^n$ is the arrow $(\1^n,\mixer)$ where $\mixer$ is the mixing voltage, and the isomorphisms are
precisely the voltage operators described in Example\nobreakspace \ref {eg:d-auto}.
It might be also interesting to study the analogous category obtained by considering rooted voltage operations.

Observe that the snub operation seems to act differently on orientable maniplexes than in non-orientable ones.
However, the result of applying the snub operation to a non-orientable maniplex $\cM$ is one of the connected components of doing the same operation to the orientable double cover of $\cM$.
This phenomenon is easy to understand with the following results:

\begin{thm}\label{thm:mixandfix}
    Let $\cX$ be an $n$-premaniplex, $(\cY_1,\mixer)$ be a mix $n$-operator and $(\cY_2,\eta)$ be an $(n,m)$-voltage operator with $\cY_2$ connected. Let $y_1\in\cY_1$ and $y_2\in\cY_2$ be fixed and suppose that $\eta(\fg^{y_2}(\cY_2))$ fixes $y_1$.
    Then, the induced (colored) graph of $(\cX\mixtimes\cY_1)\ertimes \cY_2$ with vertex set $\{((x,y_1),y):x\in\cX, y\in\cY_2\}$
is isomorphic to $\cX\ertimes\cY_2$.
\end{thm}
\begin{proof}
 Because of Theorem\nobreakspace \ref {thm:VoltsEquiv}, we may assume without loss of generality that $\eta(W)\in\eta(\Pi^{y_2}(\cY_2))$ for every path $W\in \fg(\cY_2)$, and thus $\eta(W)$ fixes $y_1$.

  First we notice that the induced (colored) graph of $\cY_1\ertimes\cY_2$ with vertex  set $\{(y_1,y):y\in \cY_2\}$ forms an isomorphic copy of $\cY_2$.
 This is easy to see since the $i$-adjacent flag to $(y_1,y)$ (in $\cY_1\ertimes\cY_2$) is $(\eta(\pth{y})y_1,y^i) = (y_1,y^i)$, because, by assumption, the voltages of paths in $\cY_2$ fix $y_1$.

 Next we use Theorem\nobreakspace \ref {thm:composition} to see $(\cX\mixtimes{\cY_1})\ertimes\cY_2$ as $\cX\ertimes[\theta] (\cY_1\ertimes \cY_2)$ with $\theta=\theta(\eta,\mixer)$. This means that we can consider the $i$-adjacent flag to $((x,y_1),y)$ in $(\cX\mixtimes{\cY_1})\ertimes\cY_2$ as the $i$-adjacent flag to $(x,(y_1,y))$ in $\cX\ertimes[\theta] (\cY_1\ertimes \cY_2)$;  we write $((x,y_1),y)\leftrightarrow (x,(y_1,y))$ to denote that these two points are in correspondence under the isomorphism.
 On one hand, observe that by definition of $\theta$, for any $y\in\cY_2$ and any monodromy $\omega$, we have that $\theta(\pth[\omega]{(y_1,y))}$ is  $\mixer(\pth[{\eta(\pth[\omega]{y}})]{y_1})$. On the other hand, by definition of $\mixer$ we have that $\mixer(\pth[{\eta(\pth[\omega]{y}})]{y_1})=\eta(\pth[\omega]{y})$.

 Hence,
 \[
 \begin{aligned}
  \big((x,y_1),y\big)^i &\leftrightarrow  \big(x,(y_1,y)\big)^i  \\
         &= \big( \theta(^i(y_1,y))x, (y_i,y)^i\big) \\
         &= \big(\mu(^{\eta(^iy)}y_1)x, (y_1,y)^i\big)  \\
         &= \big(\eta(^iy)x , (y_1,y)^i\big)  \\
         &= \big(\eta(^iy)x , (\eta(^iy)y_1, y^i)\big)  \\
         &= \big(\eta(^iy)x , (y_1,y^i)\big),
\end{aligned}
\]

and since $\big(\eta(^iy)x , (y_1,y^i)\big)
         \leftrightarrow  \big((\eta(^iy)x , y_1),y^i\big)$, the theorem follows.
\end{proof}

\begin{coro}\label{coro:mixandfixalot}
  Let $\cX$ be an $n$-premaniplex, $(\cY_1,\mixer)$ be a mix $n$-operator and $(\cY_2,\eta)$ be an $(n,m)$-voltage operator with $\cY_2$ connected. Suppose that $\eta(\fg^{y_2}(\cY_2)$ fixes every vertex of $\cY_1$.
  Then $(\cX\mixtimes \cY_1)\ertimes \cY_2$ has a copy of $\cX\ertimes \cY_2$ for each vertex of $\cY_1$.
\end{coro}

We now have the tools to understand the relation between the snub of a non-orientable 3-maniplex $\cM$ and the snub of its double cover.
In Example\nobreakspace \ref {eg:DoubleCovers} we saw that the orientable double cover of $\cM$ is $\cM\mix \2^3_\emptyset=\cM\mixtimes{\2^3_\emptyset}$.
   Hence, if $(\cY,\eta)$ is the snub operator (see Figure\nobreakspace \ref {fig:snub}),  then Corollary\nobreakspace \ref {coro:mixandfixalot} tells us that $(\cM\mix 2^3_\emptyset)\ertimes\cY$ consists of two copies of $\cM\ertimes \cY$. (Note that the voltages of $\eta$ take values in $\mon^+(\cU^3)$,the group that fixes the vertices of $\2^3_\emptyset$.)
    In other words, the snub $\cM$ is an unrooted snub $\cM\mix \2^3_\emptyset$.

    The following is a similar example. Given a  maniplex $\cM$ with a bipartite 1-skeleton (that is, the graph with the $0$-faces as vertices and the $1$-faces as edges, with the induced incidence), one can find a voltage operator $(\cY,\eta)$ such that $\cM\ertimes\cY$ has two connected components and in each component the vertices of one part of the bipartition are truncated while the ones in the other part remain unchanged.
    If now we apply this same operation to a maniplex $\cM'$ whose 1-skeleton is not bipartite, the result is one of the connected components of $(\cM'\mix \2_{\{1,2,\ldots,n-1\}}^n)\ertimes \cY$.

\section{Voltage operations on the universal maniplex}\label{sec:universal}

Since a voltage operator $(\cY,\eta)$ is a voltage premaniplex, it natural to consider the derived premaniplex $\cY^\eta$. In this section we see that this derived premaniplex is in fact the one obtained from applying the corresponding voltage operation to the universal maniplex $\cU$.
We shall use this fact to prove that Theorem\nobreakspace \ref {thm:covers} characterizes voltage operations.

\begin{thm}\label{thm:derivedAsProduct}
Let $\cU$ be the universal $n$-maniplex and let $(\cY,\eta)$ a $(n,m)$-voltage operator. Then $\cY^\eta$ is isomorphic to $\cU \ertimes \cY$.
\end{thm}
\begin{proof}
Let $\Phi$ be a base flag of $\cU$.
Consider the function $\varphi: \cY^{\eta} \to \cU\ertimes \cY$ given by $\varphi(y,\omega) = (\omega \Phi, y)$.
Recall that any flag of $\cU$ is of the form $\omega \Phi$ for some (unique) $\omega \in \mon(\cU)$, which implies that the function $\varphi$ is bijective.
Finally, observe that
\[\begin{aligned}
    \varphi(r_i(y,\omega))
    &= \varphi(r_iy,\eta(\pth{y})\omega)\\
    &= (\eta(\pth{y})\omega \Phi, r_iy )\\
    &= r_i(\omega \Phi, y)\\
    &= r_i(\varphi(y,\omega))
\end{aligned}
\]
\end{proof}

The proof of Theorem\nobreakspace \ref {thm:derivedAsProduct} is based on the fact that since $\cU$ is regular, its monodromy group acts regularly on its vertices (flags).
We can generalize Theorem\nobreakspace \ref {thm:derivedAsProduct} as follows.
 Let $\cX$ be a regular premaniplex and let $(\cY,\eta)$ be a voltage operator.
If $S_\cX = \{\omega \in \mon(\cU): \omega x= x \ \mathrm{for} \ \mathrm{all} \ x\in \cX\}$, that is, the kernel of the projection $\mon(\cU) \to \mon(\cX)$, then $\mon(\cX) \cong \mon(\cU)/ S_{\cX}$.
 Let $\eta_\cX:\fg(\cY)\to \mon(\cX)$ be the voltage assignment on $\cY$ we get by reducing $\eta$ to $\mon(\cX)$, that is $\eta_\cX(W):=\eta(W) S_\cX$.
 Then by using the exact same argument as in Theorem\nobreakspace \ref {thm:derivedAsProduct} we can see that $\cX\ertimes \cY$ is isomorphic to $\cY^{\eta_\cX}$.

An immediate consequence of Theorems\nobreakspace \ref {thm:covers} and\nobreakspace  \ref {thm:derivedAsProduct} is that operators $(\cY, \eta_{1})$ and $(\cY,\eta_{2})$ with $\eta_{1}$ and $\eta_{2}$ equivalent voltages yield equivalent voltage operations.
More precisely,

\begin{prop}\label{prop:equivalentVolt}
  Let $(\cY,\eta_{1})$ and $(\cY,\eta_{2})$ two $(n,m)$-operators with equivalent voltages $\eta_{1}$ and $\eta_{2}$.
  Then for every $n$-premaniplex $\cX$, \[\cX \ertimes[\eta_{1}] \cY \cong \cX \ertimes[\eta_{2}] \cY.\] Moreover, there is an isomorphism $\cX\ertimes[\eta_{1}] \cY \to \cX \ertimes[\eta_{2}] \cY$ such that the following diagram commutes:
\begin{equation}\label{eq:EquivVoltsOps}
    \begin{tikzcd}[column sep =small]
        \cX\ertimes[\eta_1]\cY \arrow[dashed]{rr} \arrow{rd}{} & & \cX\ertimes[\eta_2]\cY \arrow{dl}{} \\
    & \cY &
\end{tikzcd}
\end{equation}
\end{prop}
\begin{proof}
 From \cite{Hartley_1999_AllPolytopesAre} we know that $\cX \cong \cU/ \Gamma$ for certain $\Gamma \leq \Gamma(\cU)$, now
 \[\begin{aligned}
	\cX \ertimes[\eta_{1}] \cY
	&\cong (\cU / \Gamma) \ertimes[\eta_{1}] \cY \\
	&\cong \left(\cU \ertimes[\eta_{1}] \cY \right) / \Gamma \\
	&\cong \cY^{\eta_{1}} / \Gamma \\
	&\cong \cY^{\eta_{2}} / \Gamma \\
	&\cong \left(\cU \ertimes[\eta_{2}] \cY \right) / \Gamma \\
	&\cong (\cU / \Gamma) \ertimes[\eta_{2}] \cY \\
	&\cong \cX \ertimes[\eta_{2}] \cY.
\end{aligned}\]
If we start with a vertex $(x,y)$ in $\cX\ertimes[\eta_1]\cY$ and apply the natural isomorphisms between consecutive terms in the equation above, we get the following sequence, where $\Phi$ is a base flag of $\cU$:
\[
    \begin{aligned}
        (x,y)&\mapsto (\Psi\Gamma,y)&& \text{for some } \Psi\in\cU\\
        &\mapsto (\Psi,y)\Gamma\\
        &\mapsto (y,\omega)\Gamma&&\text{where } \omega\in\mon(\cU) \text{ is such that } \Psi=\omega \Phi \\
        &\mapsto (y,\tilde{\omega})\Gamma &&\text{for some } \tilde{\omega}\in\mon{\cU} \text{ (since $\eta_1$ is equivalent to $\eta_2$)}\\
        &\mapsto (\tilde{\Psi},y)\Gamma && \text{where } \tilde{\Psi}=\tilde{\omega}\Phi\\
        &\mapsto (\tilde{\Psi}\Gamma,y)\\
        &\mapsto (\tilde{x},y) &&\text{for some } \tilde{x}\in\cX.
    \end{aligned}
\]
We can see that this isomorphism makes the diagram in Equation\nobreakspace \textup {(\ref {eq:EquivVoltsOps})} commute.
\end{proof}

In the proof of Proposition\nobreakspace \ref {prop:equivalentVolt} we strongly used that if $\cX \cong \cU/\Gamma$ for $\Gamma \leq \aut(\cU)$, then $\cX \ertimes \cY \cong (\cU \ertimes \cY)/\Gamma$ for any voltage operator $(\cY,\eta)$ (which is Theorem\nobreakspace \ref {thm:covers} of this paper).
In the following result we will show that this property characterizes all voltage operations.

\begin{thm}\label{thm:OpsAsProducts}
 Let $\oo$ be a mapping that assigns an $m$-premaniplex $\oo(\cX)$ to each $n$-premaniplex $\cX$.
 Assume that there is an action of $\Gamma(\cU)$ on $\oo(\cU)$ such that  $\oo(\cU / \Gamma) \cong \oo(\cU)/\Gamma$ for every $\Gamma \leq \aut(\cU)$, then there exists an $(n,m)$-voltage operator $(\cY, \eta)$ such that \[\oo(\cX) \cong \cX \ertimes \cY\] for every premaniplex $\cX$.
\end{thm}
\begin{proof}
 Let $\1^{n}$ denote the unique $n$-premaniplex with one vertex and define $\cY$ as the premaniplex $\oo(\1^{n})$.
 Observe that
 \[
\oo(\cU)/\aut(\cU) \cong \oo\left( \cU/\aut(\cU) \right) \cong \oo(\1^{n}) = \cY,
 \]
 which implies that there exists a voltage assignment $\eta:\fg(\cY) \to \aut(\cU) \cong \mon(\cU)$ such that $\cY^{\eta} \cong \oo(\cU)$.
 The pair $(\cY, \eta)$ defines a voltage operator and by Theorem\nobreakspace \ref {thm:derivedAsProduct}, $\cU \ertimes \cY  \cong \cY^{\eta} \cong \oo(\cU) $.

 Finally, if $\cX$ is a premaniplex and $\Gamma \leq \aut(\cU)$ is such that $\cX \cong \cU/\Gamma$, then
 \[
	\oo(\cX)
	\cong \oo(\cU/\Gamma)
	\cong \oo(\cU)/ \Gamma
	\cong (\cU \ertimes \cY)/\Gamma
	\cong \left( \cU/\Gamma \right) \ertimes \cY
	\cong \cX \ertimes \cY.
\qedhere
\]
 \end{proof}

We often come across operations that are well defined for some family $F$ of premaniplexes (for example, for convex polytopes) but such that it is not immediately evident how to generalize them for all premaniplexes.
One may ask if it is possible to extend operations of this kind to all premaniplexes (of the given rank) in such a way that the operation is a voltage operation.
As an example, in Section\nobreakspace \ref {sec:ClassicalExamples}, we have found voltage operations that extend the Wythoffian operations, the $k$-bubble and the trapezotope operation to all premaniplexes.
Theorem\nobreakspace \ref {thm:OpsAsProducts} and Corollary\nobreakspace \ref {coro:OperatingRegulars} answer this question and find the corresponding voltage operation, when possible.
The idea is as follows: suppose there is some regular premaniplex $\cP$ and an operation $\oo$ such that we already know $\oo(\cP)$.
According to the proof of Theorem\nobreakspace \ref {thm:OpsAsProducts}, if $\oo$ is indeed a voltage operation $\cX\mapsto \cX\ertimes \cY$, then $\cY$ should be $\oo(\1^n)$, but due to Theorem\nobreakspace \ref {thm:covers} $\cY$ must coincide with $\oo(\cP)/\aut(\cP)$.
Now let $\nu:\fg(\cY)\to\aut(\cP)$ be the voltage assignment such that $\oo(\cP)$ is isomorphic  to $\cY^\nu$.
Corollary\nobreakspace \ref {coro:OperatingRegulars} tells us that $\nu$ is obtained by replacing each occurrence of $r_i$ by $\rho_i$ in the voltage assignment $\eta:\fg(\cY)\to \mon(\cU)$.
So we can recover $\eta$ by replacing every instance of $\rho_i$ by $r_i$ in $\nu$.
Note that, in general, such replacement is not well defined unless $\cP$ is the universal polytope $\cU$, however with some intuition we can find the right way to do it for the natural occurring operations.
However, observe that one can do this without knowing if $\oo$ was in fact a voltage operation.
To know if $\oo$ is a voltage operation or not, one must see that $\oo(\cX)=\cX\ertimes \cY$ for $\cX$ in the family $F$. In fact, if for some $\cX \in F$ we observe that $\oo(\cX)\neq\cX\ertimes \cY$, then $\oo$ cannot be seen as a voltage operation.

 \section{Final remarks and open problems}

We have seen that voltage operations generalize classical operations on maps and polytopes and allow us to define such classical operations on premaniplexes.

One can see that voltage operations naturally generalize to hypertopes (thin residually connected geometries), ``complexes'', and their quotients.
Following \cite{Wilson_2012_ManiplexesPart1}, a {\em complex}  is a properly $n$-edge colored $n$-valent graph (note that a complex is a combinatorial map in the sense of Vince \cite{Vince_1983_CombinatorialMaps}. Complexes generalize (the ``chamber graphs'' of) hypertopes. Thus, one can define a precomplex as a $n$-valent pregraph that has been properly $n$-edge colored.
In this context, the definition of a voltage operator as well as the results in this paper hold if we ask $\eta:\Pi(\cY)\to W^n$, where $W^n$ is the group generated by $n$ involutions with no other relations among the generators.
In particular, this allows us to define operations like the ones given in Sections~\ref{sec:operations} and~\ref{sec:ClassicalExamples} to hypertopes.

Symmetry type graphs and voltage assignments have proven to be strong tools to study polytopes and maniplexes.
A natural (and well-known)  problem that arises when dealing with symmetry type graphs is the following:
\begin{problem}\label{prob:stg}
Given a premaniplex $T$, does there exist a polytope (or maniplex) such that its symmetry type graph (with respect to the full automorphism group) is $T$?
\end{problem}
A particular example of this problem was the question posted in the early 1990's by Schulte and Weiss of whether or not there exist chiral polytopes of all ranks (which was solved by Pellicer in 2010 \cite{Pellicer_2010_ConstructionHigherRank}
Of course, one can generalize \protect \MakeUppercase {P}roblem\nobreakspace \ref {prob:stg} to hypertopes.
\begin{problem}
Given a precomplex $T$, does there exist a hypertope (or complex) such that its symmetry type graph (with respect to the full automorphism group) is $T$?
\end{problem}

One thing that one tries to do when dealing the above question is to use voltage assignments on $T$ to lift to a maniplex $\cM$. In \cite{Mochan_2021_AbstractPolytopesTheir_PhDThesis} the voltage assignments that give a polytope as the derived graph are characterized. Of course the voltage group acts by automorphisms on $\cM$ and the quotient of $\cM$ by it is precisely $T$. However, $\cM$ can have (and often does) symmetries not coming from the voltage group.
In the context of voltage operations, we have shown that given a voltage operator $(\cY, \eta)$ and a premaniplex $\cX$, all automorphisms of $\cX$ act as automorphisms of $\cX\ertimes[\eta]\cY$. However, again, $\cX\ertimes[\eta]\cY$ might have extra symmetry (for example, in the case when one applies the medial operation to a self-dual map).
It is natural to ask when is it true that automorphisms of $\cY$ also lift to $\cX\ertimes[\eta]\cY$, and when all automorphisms of $\cX\ertimes[\eta]\cY$ come either from $\cX$ or from $(\cY, \eta)$. We refer to \cite{HubardMochanMontero_MoreVoltageOperations_preprint} for more details about these questions.
Answering, at least partially, these questions might be of great help \protect \MakeUppercase {P}roblem\nobreakspace \ref {prob:stg} (at least partially).
\begin{problem}
Give necessary conditions on $(\cY,\eta)$ so that one can compute the symmetry type graph of  $\cX\ertimes[\eta]\cY$ with respect to its full automorphism group in terms of $\cX$ and $(\cY,\eta)$.
\end{problem}

If one's interest is on polytopes rather that in (pre)maniplexes, one can ask when a voltage operation preserves polytopality, though we might not care if the result of the operation is not connected. More precisely,
\begin{problem}
Give necessary and sufficient conditions on $(\cY, \eta)$ for $(\cP,\Phi)\ertimes(\cY,y)$ to be a polytope, for all rooted polytopes $(\cP, \Phi)$ and $y\in\cY$.
\end{problem}
Similarly,
\begin{problem}
Give necessary and sufficient conditions on $(\cY, \eta)$ for $({\mathcal H},\Phi)\ertimes(\cY,y)$ to be a hypertope, for all rooted hypertopes $({\mathcal H}, \Phi)$ and $y\in\cY$.
\end{problem}
 
\section*{Acknowledgements}
The authors thank the financial support of CONACyT grant A1-S-21678. 
The third author was supported by the Post Doctoral Scholarship Program at UNAM (DGAPA), Mexico.

\printbibliography
\end{document}